\documentclass[12pt]{article}
\usepackage{amssymb, amsmath, amsthm, amscd}
\usepackage[pdftex]{graphics}
\usepackage[utf8]{inputenc}
\usepackage[all,cmtip]{xy}
\usepackage{enumitem}
\usepackage{setspace}
\usepackage{mathdots}

\usepackage[mathscr]{euscript}
\usepackage[colorlinks=true]{hyperref}
\hypersetup{colorlinks   = true}
\hypersetup{linkcolor=blue}
\usepackage{tikz}

\begin{document}

\newtheorem{The}{Theorem}[section]
\newtheorem{Lem}[The]{Lemma}
\newtheorem{Prop}[The]{Proposition}
\newtheorem{Cor}[The]{Corollary}
\newtheorem{Rem}[The]{Remark}
\newtheorem{Obs}[The]{Observation}
\newtheorem{SConj}[The]{Standard Conjecture}
\newtheorem{Titre}[The]{\!\!\!\! }
\newtheorem{Conj}[The]{Conjecture}
\newtheorem{Question}[The]{Question}
\newtheorem{Prob}[The]{Problem}
\newtheorem{Def}[The]{Definition}
\newtheorem{Not}[The]{Notation}
\newtheorem{Claim}[The]{Claim}
\newtheorem{Conc}[The]{Conclusion}
\newtheorem{Ex}[The]{Example}
\newtheorem{Fact}[The]{Fact}
\newtheorem{Formula}[The]{Formula}
\newtheorem{Formulae}[The]{Formulae}
\newtheorem{The-Def}[The]{Theorem and Definition}
\newtheorem{Prop-Def}[The]{Proposition and Definition}
\newtheorem{Lem-Def}[The]{Lemma and Definition}
\newtheorem{Cor-Def}[The]{Corollary and Definition}
\newtheorem{Conc-Def}[The]{Conclusion and Definition}
\newtheorem{Terminology}[The]{Note on terminology}
\newcommand{\C}{\mathbb{C}}
\newcommand{\R}{\mathbb{R}}
\newcommand{\N}{\mathbb{N}}
\newcommand{\Z}{\mathbb{Z}}
\newcommand{\Q}{\mathbb{Q}}
\newcommand{\Proj}{\mathbb{P}}
\newcommand{\Rc}{\mathcal{R}}
\newcommand{\Oc}{\mathcal{O}}
\newcommand{\Vc}{\mathcal{V}}
\newcommand{\Id}{\operatorname{Id}}
\newcommand{\pr}{\operatorname{pr}}
\newcommand{\rk}{\operatorname{rk}}
\newcommand{\del}{\partial}
\newcommand{\delbar}{\bar{\partial}}
\newcommand{\Cdot}{{\raisebox{-0.7ex}[0pt][0pt]{\scalebox{2.0}{$\cdot$}}}}
\newcommand\nilm{\Gamma\backslash G}
\newcommand\frg{{\mathfrak g}}
\newcommand{\fg}{\mathfrak g}
\newcommand{\Oh}{\mathcal{O}}
\newcommand{\Kur}{\operatorname{Kur}}
\newcommand\gc{\frg_\mathbb{C}}

\begin{center}

{\Large\bf Pluriclosed Star Split Hermitian Metrics}

\end{center}

\begin{center}

{\large Dan Popovici}

\end{center}

\vspace{1ex}

\noindent{\small{\bf Abstract.} We introduce a class of Hermitian metrics, that we call pluriclosed star split, generalising both the astheno-K\"ahler metrics of Jost and Yau and the $(n-2)$-Gauduchon metrics of Fu-Wang-Wu on complex manifolds. They have links with Gauduchon and balanced metrics through the properties of a smooth function associated with any Hermitian metric. After pointing out several examples, we generalise the property to pairs of Hermitian metrics and to triples consisting of a holomorphic map between two complex manifolds and two Hermitian metrics, one on each of these manifolds. Applications include an attack on the Fino-Vezzoni conjecture predicting that any compact complex manifold admitting both SKT and balanced metrics must be K\"ahler, that we answer affirmatively under extra assumptions. We also introduce and study a Laplace-like differential operator of order two acting on the smooth $(1,\,1)$-forms of a Hermitian manifold. We prove it to be elliptic and we point out its links with the pluriclosed star split metrics and pairs defined in the first part of the paper.}

\vspace{2ex}

\section{Introduction}\label{section:Introduction} A central problem in complex geometry aims at classifying compact complex manifolds that share certain properties, either metric or cohomological, with K\"ahler manifolds. In this paper, we pursue the former approach to classification, the one aimed at studying the geometry of these manifolds by means of investigating the types of special Hermitian metrics they support.

Let $X$ be an $n$-dimensional compact complex manifold. Recall that a Hermitian metric on $X$ is defined by any $C^\infty$ positive definite $(1,\,1)$-form $\omega$ on $X$. Any such object can be written locally, in terms of a system of local holomorphic coordinates $z_1,\dots , z_n$ on $X$, as $\omega = \sum_{1\leq j,\,k\leq n}\omega_{j\bar{k}}\,idz_j\wedge d\bar{z}_k$, where the coefficient matrix $(\omega_{j\bar{k}})_{j,\,k}$ is positive definite at every point of the local coordinate domain. Hermitian metrics always exist, but if an extra condition is imposed thereon, the resulting type of metrics need not exist. When they do, they often give useful geometric information on the underlying manifold $X$. Among the best known types of special Hermitian metrics are the ones listed below, of which only the Gauduchon metrics always exist by [Gau77a]. A Hermitian metric $\omega$ is said to be:

\vspace{1ex}

$\cdot$ {\it K\"ahler}, if $d\omega=0$;

\vspace{1ex}

$\cdot$ {\it SKT}, if $\partial\bar\partial\omega=0$ (see [HP96]); 

\vspace{1ex}

$\cdot$ {\it astheno-K\"ahler}, if $\partial\bar\partial\omega^{n-2}=0$ (see [JY93]);

\vspace{1ex}

$\cdot$ {\it $(n-2)$-Gauduchon}, if $\omega\wedge\partial\bar\partial\omega^{n-2}=0$ (see [FWW13]);

\vspace{1ex}

$\cdot${\it balanced}, if $d\omega^{n-1}=0$ (see [Gau77b]);

\vspace{1ex}

$\cdot$ {\it Gauduchon}, if $\partial\bar\partial\omega^{n-1}=0$ (see [Gau77a]).

\vspace{2ex}

An interesting conjecture was proposed by Fino and Vezzoni in [FV15]. It predicts that a K\"ahler metric ought to exist on any compact complex manifold $X$ that carries an SKT metric $\gamma$ and a balanced metric $\omega$. This is currently known to be true only when $\gamma=\omega$ (see [IP13], also [Pop15] for a shorter proof). If this SKT-balanced conjecture is borne out, it will provide an efficient criterion for the K\"ahlerianity of a given manifold $X$.

\vspace{2ex}

Throughout the paper, we will use the convenient notation $$\eta_p:=\frac{\eta^p}{p!}$$ for any positive integer $2\leq p\leq n$ and any (semi-)positive $(1,\,1)$-form $\eta\geq 0$ on $X$.

\vspace{2ex}


In the first part of the paper, motivated by the classification problem for K\"ahler-like compact complex manifolds from a metric perspective and, more specifically, by the above-mentioned SKT-balanced conjecture of Fino and Vezzoni, we introduce a new class of Hermitian metrics starting from the well-known fact that, for any Hermitian metric $\omega$ on an $n$-dimensional complex manifold $X$ and any non-negative integer $p$ such that $2p+1\leq n$, the pointwise map of multiplication of $(p,\,p)$-forms by $\omega_{n-2p}$: \begin{eqnarray*}\omega_{n-2p}\wedge\cdot : \Lambda^{p,\,p}T^\star X\longrightarrow\Lambda^{n-p,\,n-p}T^\star X\end{eqnarray*} is {\bf bijective}. In particular, there exists a unique $(p,\,p)$-form $\rho_\omega^{(p)}$ such that \begin{eqnarray*}i\partial\bar\partial\omega_{n-p-1} = \omega_{n-2p}\wedge\rho_\omega^{(p)}.\end{eqnarray*} Applying the Hodge star operator $\star=\star_\omega$ to $\rho_\omega^{(p)}$, we get the $(n-p,\,n-p)$-form $\star\rho_\omega^{(p)}$ which can be computed explicitly in terms of $\omega$.

In this paper, we concentrate on the case when $p=1$. Setting $\rho_\omega = \rho_\omega^{(1)}$, a computation yields \begin{eqnarray*}\star\rho_\omega = \frac{1}{n-1}\,f_\omega\,\omega_{n-1} - i\partial\bar\partial\omega_{n-2},\end{eqnarray*} where $f_\omega$ is the real-valued $C^\infty$ function on $X$ defined by $f_\omega = (\omega\wedge i\partial\bar\partial\omega_{n-2})/\omega_n$. We say (cf. Definition \ref{Def:pluriclosed-star-split_metrics}) that the Hermitian metric $\omega$ on $X$ is {\bf pluriclosed star split}, respectively {\bf closed star split}, if $\partial\bar\partial(\star\rho_\omega) =0$, respectively if $d(\star\rho_\omega) =0$. 

 The defining property $i\partial\bar\partial\omega_{n-2} = 0$ of astheno-K\"ahler metrics $\omega$ (cf. [JY93]) is equivalent to $\rho_\omega=0$ since the map (\ref{eqn:omega_n-2_map}) is bijective, hence it is also equivalent to $\star\rho_\omega = 0$ since the Hodge star operator $\star = \star_\omega$ is bijective. In particular, every astheno-K\"ahler metric is closed star split, hence also pluriclosed star split. Besides generalising the astheno-K\"ahler metrics of Jost-Yau, the pluriclosed star split metrics also generalise the $(n-2)$-Gauduchon metrics of Fu-Wang-Wu. The defining property $\omega\wedge i\partial\bar\partial\omega_{n-2} =0$ of the latter metrics (cf. [FWW13]) is equivalent to the function $f_\omega$ defined above vanishing identically, a fact that amounts to $\star\rho_\omega = - i\partial\bar\partial\omega_{n-2}$. In particular, $\star\rho_\omega$ is $d$-closed, hence also $\partial\bar\partial$-closed, if $\omega$ is $(n-2)$-Gauduchon. Therefore, we obtain the following implications: \begin{eqnarray*}\omega \hspace{1ex}\mbox{is {\bf astheno-K\"ahler}} \implies \omega \hspace{1ex}\mbox{is}\hspace{1ex} (n-2)\mbox{-{\bf Gauduchon}} \implies & \omega \hspace{1ex}\mbox{is {\bf closed star split}} &  \\
 & \rotatebox{-90}{$\implies$} & \\
 & \omega \hspace{1ex}\mbox{is {\bf pluriclosed star split}}.& \end{eqnarray*}

 On the other hand, using a result from [Gau77a], we notice the following consequence of Proposition \ref{Prop:pluriclosed-star-split_3-choices_f}:

 \begin{Prop}\label{Prop:pluriclosed-star-split_3-choices_f_introd} Let $X$ be a compact connected complex manifold with $\mbox{dim}_\C X = n$ and let $\omega$ be a {\bf pluriclosed star split} Hermitian metric on $X$. Then, its associated function $f_\omega$ satisfies one of the following three conditions: \begin{eqnarray*}f_\omega> 0 \hspace{1ex} \mbox{on} \hspace{1ex} X \hspace{5ex} \mbox{or} \hspace{5ex} f_\omega< 0 \hspace{1ex} \mbox{on} \hspace{1ex} X  \hspace{5ex} \mbox{or} \hspace{5ex} f_\omega\equiv 0.\end{eqnarray*} 

\end{Prop} 

 The case where $f_\omega$ is {\it constant} stands out (cf. Proposition \ref{Prop:pluriclosed-star-split_properties}). 








\vspace{2ex}   

Examples of compact complex manifolds admitting pluriclosed star split Gauduchon (even balanced) metrics $\omega$ include (cf. $\S.$\ref{subsection:examples_pluriclosed-star-split-metrics}): the Iwasawa manifold $I^{(3)}$ (for which $f_\omega\equiv 1$) and its small deformations; the Nakamura manifolds (for which $f_\omega\equiv 2$); the $5$-dimensional Iwasawa manifold $I^{(5)}$ (for which $f_\omega\equiv 3$); the $3$-dimensional Calabi-Eckmann manifold $(S^3\times S^3,\,J_{CE})$ and its small deformations. The constant values the function $f_\omega$ assumes in these cases seem to point to some intrinsic features of these manifolds that are yet to be singled out. For example, in the cases of $I^{(3)}$, $I^{(5)}$ and the Nakamura manifolds, the constant $f_\omega$ equals the number of holomorphic $(1,\,0)$-forms that are not $d$-closed among the canonical forms that determine the cohomology of these manifolds. While the actual value of $f_\omega$ may change when $\omega$ changes, for example as described in (\ref{eqn:rescaling_omega_f-omega}), it is worth wondering if and when the sign of $f_\omega$ depends only on the complex structure of $X$. Lemma 3.2 in [FU13]\footnote{The author is grateful to L. Ugarte for pointing out this result and this reference to him.} shows that this is, indeed, the case for nilmanifolds of complex dimension $3$: the sign of $f_\omega$ ($=$ the sign of $\gamma_1(\Omega)$ in the notation of that paper) remains the same as $\omega$ ranges over all the invariant Hermitian metrics on such a manifold. 

\vspace{2ex}

As a consequence of the proof of Proposition \ref{Prop:astheno-balanced-K} and other observations, we get the following information about the relations between the balanced condition and the sign of the function $f_\omega$ of a pluriclosed star split metric $\omega$.

\begin{The}\label{The:bal-pluri_star-split_f-not-neg} (i)\, Let $\omega$ be a Hermitian metric on a compact connected complex manifold $X$ with $\mbox{dim}_\C X = n\geq 3$. If $\omega$ is {\bf balanced} and {\bf pluriclosed star split}, the (necessarily constant) function $f_\omega$ is {\bf non-negative}.

\vspace{1ex}
  
  (ii)\, There are examples of compact Hermitian manifolds $(X,\,\omega)$ such that the metric $\omega$ is {\bf pluriclosed star split} and {\bf non-balanced}, while the function $f_\omega$ is a {\bf negative constant}.

\end{The}  

\vspace{2ex}

One of the results we get in $\S.$\ref{subsection:gen-theory_one-metric} is the following characterisation of the K\"ahlerianity of a Hermitian metric $\omega$ in terms of the associated function $f_\omega$. It is a consequence of Proposition \ref{Prop:balanced_f-omega_del-omega}.

\begin{Prop}\label{Prop:balanced_f-omega_del-omega_introd} Suppose there exists a {\bf balanced} metric $\omega$ on a compact complex manifold $X$ with $\mbox{dim}_\C X =n\geq 3$.

  Then, $\omega$ is K\"ahler if and only if $\int_Xf_\omega\,\omega_n=0$. If $f_\omega$ is constant, $\omega$ is K\"ahler if and only if $f_\omega=0$.

\end{Prop}

Note that the condition $\int_Xf_\omega\,\omega_n=0$ is equivalent to $f_\omega$ being $L^2_\omega$-orthogonal to the constant functions on $X$, in other words $L^2_\omega$-orthogonal to $\ker\Delta_\omega$ since the Laplacian $\Delta_\omega:=\Lambda_\omega(i\partial\bar\partial)$ acting on functions is elliptic and $X$ is compact. Thus, the condition $\int_Xf_\omega\,\omega_n=0$ is equivalent to $f_\omega$ lying in the image of $\Delta_\omega^\star$, the adjoint of $\Delta_\omega$. On the other hand, when $\omega$ is balanced, $\Delta_\omega = \Delta_\omega^\star$, so in that case $\omega$ is K\"ahler if and only if $f_\omega\in\mbox{Im}\,\Delta_\omega$.

\vspace{3ex}

In $\S.$\ref{section:def_two-metrics}, we generalise the pluriclosed star split and the closed star split conditions to pairs $(\omega,\,\gamma)$ of Hermitian metrics by performing the division of $i\partial\bar\partial\omega_{n-2}$ by $\gamma_{n-2}$ (rather than $\omega_{n-2}$) to obtain a unique smooth $(1,\,1)$-form $\rho_{\omega,\,\gamma}$ satisfying \begin{eqnarray*}i\partial\bar\partial\omega_{n-2} = \gamma_{n-2}\wedge \rho_{\omega,\,\gamma}.\end{eqnarray*} Taking the Hodge star operator $\star = \star_\gamma$ induced (again) by $\gamma$, we obtain the $(n-1,\,n-1)$-form $\star_\gamma\rho_{\omega,\,\gamma}$ on the given $n$-dimensional complex manifold $X$. It is expressed in terms of $\omega$ and $\gamma$ by means of a $C^\infty$ function $f_{\omega,\,\gamma}$ that has similar properties to those of $f_\omega$. 

As in the case of one metric, we say (cf. Definition \ref{Def:pluriclosed-star-split_pairs}) that the pair $(\omega,\,\gamma)$ of Hermitian metrics is {\bf pluriclosed star split}, respectively {\bf closed star split}, if $\partial\bar\partial(\star_\gamma\rho_{\omega,\,\gamma}) = 0$, respectively if $d(\star_\gamma\rho_{\omega,\,\gamma}) = 0$.

Thanks to the extra flexibility they afford, (pluriclosed star split) pairs of Hermitian metrics seem better suited to investigating the SKT-balanced conjecture of Fino-Vezzoni than single metrics. We get the following positive answer to this conjecture under the extra {\it pluriclosed star split} assumption on the pair $(\omega,\,\gamma)$ and the extra semi-definiteness assumption on the $(1,\,1)$-form $\rho_{\omega,\,\gamma}$ that plays in this context a role analogous to that of the curvature.

\begin{The}\label{The:SKT_pss_rho-pos_K} Let $\omega$ and $\gamma$ be Hermitian metrics on a compact connected complex manifold $X$ with $\mbox{dim}_\C X =n\geq 3$. Suppose that $\gamma$ is {\bf SKT}, the pair $(\omega,\,\gamma)$ is {\bf pluriclosed star split} and the associated real $(1,\,1)$-form $\rho_{\omega,\,\gamma}$ is either {\bf positive semi-definite} or {\bf negative semi-definite}. 

Then, $\omega$ is {\bf astheno-K\"ahler}. If, moreover, $\omega$ is {\bf balanced}, then $\omega$ is {\bf K\"ahler}. 

\end{The}

\vspace{2ex}

In $\S.$\ref{section:maps_pluriclosed}, we allow ourselves an even greater flexibility by extending the main construction of this paper to the case where the two Hermitian metrics $\omega$ and $\gamma$ exist on two possibly different complex manifolds $Y$, respectively $X$, with $\mbox{dim} X\leq\mbox{dim} Y$, between which holomorphic maps, supposed non-degenerate at some point, are considered: \begin{eqnarray*}\phi:(X,\,\gamma)\longrightarrow(Y,\,\omega).\end{eqnarray*}

We extend the {\it pluriclosed star split} and {\it closed star split} conditions to the triple $(\phi,\,\omega,\,\gamma)$ by means of the same construction as in the case of pairs applied to the degenerate metric $\widetilde\omega:=\phi^\star\omega$ and the genuine metric $\gamma$ on $X$. In particular, we associate a unique smooth $(1,\,1)$-form $\rho_{\phi,\,\omega,\,\gamma}$ to any such triple $(\phi,\,\omega,\,\gamma)$.   

When the manifolds $X$ and $Y$ coincide, we obtain the following generalisation of Theorem \ref{The:SKT_pss_rho-pos_K} that highlights the role of the automorphism group $Aut(X)$ of a given compact complex manifold $X$ and may prove useful in a future attack on the SKT-balanced conjecture.

\begin{The}\label{The:SKT-bal_map_pss_pos} Let $X$ be a compact complex manifold with $\mbox{dim}_\C X =n\geq 3$. Suppose there exists an {\bf SKT} metric $\gamma$ and a {\bf balanced} metric $\omega$ on $X$.

If there exists a biholomorphism $\phi:X\longrightarrow X$ such that the triple $(\phi,\,\omega,\,\gamma)$ is {\bf pluriclosed star split} and the $(1,\,1)$-form $\rho_{\phi,\,\omega,\,\gamma}$ is either {\bf positive semi-definite} or {\bf negative semi-definite} on $X$, the metric $\omega$ is {\bf K\"ahler}.

\end{The}   

It turns out (cf. Corollary \ref{Cor:subgroups}) that the $\gamma$-isometries of $X$ (i.e. the biholomorphisms $\phi:X\longrightarrow X$ satisfying $\phi^\star\gamma = \gamma$) for which the triple $(\phi,\,\omega,\,\gamma)$ is pluriclosed star split form a subgroup of the automorphism group $Aut(X)$ of $X$. Theorem \ref{The:SKT-bal_map_pss_pos} is a special case of Theorem \ref{The:map_two-metrics_positivity}. The latter can be seen in the context of Siu's Theorems 1 and 5 on rigidity in [Siu80] -- see comments just after the statement of Theorem \ref{The:map_two-metrics_positivity}.  

\vspace{3ex}

The second part of the paper ($\S.$\ref{section:further-applications}) is motivated by a desire to extend in the long run the resolution of certain geometric PDEs, such as the celebrated Monge-Amp\`ere equation, from the case where the solution is a function (i.e. a differential form of degree $0$) to the case where it is a form of positive degree. This would have numerous geometric appplications. One such equation was proposed in [DP21, $\S.5$, $(\star)$, p. 30]. In the same vein, functionals acting on positive-degreed differential forms were introduced in [DP22], where the first variation of several operators depending on Hermitian metrics and featuring frequently in complex geometric problems was also computed.

Continuing this effort, we introduce in Definition \ref{Def:P_omega_def} a differential operator \begin{eqnarray*}P_\omega:C^\infty_{1,\,1}(X,\,\C)\longrightarrow C^\infty_{1,\,1}(X,\,\C)\end{eqnarray*} of order two that seems to be a natural analogue acting on smooth $(1,\,1)$-forms of the standard Laplacian $\Delta_\omega:C^\infty(X,\,\C)\longrightarrow C^\infty(X,\,\C)$ acting on smooth functions $\varphi$ by $\Delta_\omega(\varphi) = \Lambda_\omega(i\partial\bar\partial\varphi)$. Once a Hermitian metric $\omega$ has been fixed on an $n$-dimensional complex manifold $X$, we define \begin{eqnarray*}P_\omega(\alpha) = (\omega_{n-2}\wedge\cdot)^{-1}(i\partial\bar\partial\alpha\wedge\omega_{n-3}), \hspace{5ex} \alpha\in C^\infty_{1,\,1}(X,\,\C),\end{eqnarray*} where $(\omega_{n-2}\wedge\cdot)^{-1}:C^\infty_{n-1,\,n-1}(X,\,\C)\longrightarrow C^\infty_{1,\,1}(X,\,\C)$ is the inverse of the bijection that multiplies $(1,\,1)$-forms by $\omega_{n-2}$. 

Lemma \ref{Lem:star-rho-P_integral-link} reveals a natural link, in the form of an integral equality, between the operator $P_\omega$ and the $(n-1,\,n-1)$-form $\star_\gamma\rho_{\omega,\,\gamma}$ that we associate with every pair $(\omega,\,\gamma)$ of Hermitian metrics on $X$ via the construction leading to {\it pluriclosed star split} pairs of metrics. 

On the other hand, $P_\omega$ can be explicitly computed as (cf. Lemma \ref{Lem:P-omega_trace-trace-square}): \begin{eqnarray*}P_\omega(\alpha) = \Lambda_\omega(i\partial\bar\partial\alpha) - \frac{1}{2(n-1)}\,\Lambda_\omega^2(i\partial\bar\partial\alpha)\,\omega,\end{eqnarray*} for every smooth $(1,\,1)$-form $\alpha$. 

We then render $P_\omega$ elliptic by adding to it certain terms that we compute in $\S.$\ref{subsubsection:computation_P_omega_alpha-alpha} and $\S.$\ref{subsubsection:R-omega}. We get an elliptic second-order differential operator $Q_\omega:C^\infty_{1,\,1}(X,\,\C)\longrightarrow C^\infty_{1,\,1}(X,\,\C)$ that differs from the usual $\bar\partial$-Laplacian $\Delta''_\omega = \bar\partial\bar\partial^\star + \bar\partial^\star\bar\partial$ multiplied by $-1$ by lower-order terms (cf. Corollary \ref{Cor:Q_omega_elliptic}): \begin{eqnarray*}Q_\omega = -\Delta''_\omega + l.o.t.\end{eqnarray*} It is a key feature of this construction that the above-mentioned integral equality of Lemma \ref{Lem:star-rho-P_integral-link} satisfied by $P_\omega$ continues to be satisfied by the elliptic operator $Q_\omega$ when the metric $\omega$ is balanced (cf. Lemma \ref{Lem:star-rho-Q_integral-link}). In particular, $Q_\omega$ retains a link with the pluriclosed star split pairs of metrics introduced in the first part of the paper.

\vspace{2ex}

\noindent {\bf Acknowledgments.} The author is grateful to S. Dinew and L. Ugarte for useful comments on this text, as well as to the referee for their careful reading of the preprint and for interesting suggestions.

\section{Pluriclosed star split and closed star split metrics}\label{section:def_one-metric} Let $X$ be an $n$-dimensional complex manifold with $n\geq 3$. For any integer $p=1,\dots , n$ and any $(1,\,1)$-form $\alpha$ on $X$, we will use the following notation throughout the text: $\alpha_p:=\alpha^p/p!$.

\subsection{Definitions and general properties}\label{subsection:gen-theory_one-metric} Fix an arbitrary Hermitian metric $\omega$ on $X$. Since the pointwise map \begin{eqnarray}\label{eqn:omega_n-2_map}\omega_{n-2}\wedge\cdot :\Lambda^{1,\,1}T^\star X\longrightarrow\Lambda^{n-1,\,n-1}T^\star X\end{eqnarray} is bijective, there exists a unique $C^\infty$ $(1,\,1)$-form $\rho_\omega$ on $X$ such that \begin{eqnarray}\label{eqn:rho_omega_def}i\partial\bar\partial\omega_{n-2} = \omega_{n-2}\wedge\rho_\omega.\end{eqnarray}

The $(1,\,1)$-form $\rho_\omega$ can be computed explicitly in the following way. The standard Lefschetz decomposition theorem ensures that $\rho_\omega$ can be written in a unique way as \begin{eqnarray*}\rho_\omega = \rho_{prim} + g\,\omega,\end{eqnarray*} where $\rho_{prim}$ is a primitive (w.r.t. $\omega$) $(1,\,1)$-form and $g$ is a function. Taking $\Lambda_\omega$, we get $\Lambda_\omega(\rho_{prim})=0$ (since $\rho_{prim}$ is primitive), hence $g=(1/n)\,\Lambda_\omega(\rho_\omega)$, so \begin{eqnarray}\label{eqn:Lefschetz-decomp_rho}\rho_\omega = \rho_{prim} + \frac{1}{n}\,\Lambda_\omega(\rho_\omega)\,\omega.\end{eqnarray}

Taking the Hodge star operator $\star = \star_\omega$ in (\ref{eqn:Lefschetz-decomp_rho}) and using the following standard formula (cf. e.g. [Voi02, Proposition 6.29, p. 150]) for this operator acting on {\it primitive} forms $v$ of arbitrary bidegree $(p, \, q)$: \begin{eqnarray}\label{eqn:prim-form-star-formula-gen}\star\, v = (-1)^{k(k+1)/2}\, i^{p-q}\, \omega_{n-p-q}\wedge v, \hspace{2ex} \mbox{where}\,\, k:=p+q,\end{eqnarray} we get: \begin{eqnarray}\label{eqn:Lefschetz-decomp_rho_app1}\star\rho_\omega = -\omega_{n-2}\wedge\rho_{prim} + \frac{1}{n}\,\Lambda_\omega(\rho_\omega)\,\omega_{n-1},\end{eqnarray} since it is well known that $\star\omega = \omega_{n-1}$.

 On the other hand, after multiplying (\ref{eqn:Lefschetz-decomp_rho}) by $\omega_{n-2}$, we get \begin{eqnarray}\label{eqn:Lefschetz-decomp_rho_app2}\omega_{n-2}\wedge\rho_\omega = \omega_{n-2}\wedge\rho_{prim} + \frac{n-1}{n}\,\Lambda_\omega(\rho_\omega)\,\omega_{n-1}.\end{eqnarray} So, adding (\ref{eqn:Lefschetz-decomp_rho_app1}) and (\ref{eqn:Lefschetz-decomp_rho_app2}) up and using (\ref{eqn:rho_omega_def}), we get \begin{eqnarray}\label{eqn:Lefschetz-decomp_rho_app3}\star\rho_\omega = \Lambda_\omega(\rho_\omega)\,\omega_{n-1} - i\partial\bar\partial\omega_{n-2}.\end{eqnarray}

Meanwhile, multiplying (\ref{eqn:rho_omega_def}) by $1/(n-1)\,\omega$, we get the first equality below: \begin{eqnarray*}\frac{1}{n-1}\omega\wedge i\partial\bar\partial\omega_{n-2} = \omega_{n-1}\wedge\rho_\omega = \Lambda_\omega(\rho_\omega)\,\omega_n.\end{eqnarray*} Together with (\ref{eqn:Lefschetz-decomp_rho_app3}), this proves the following

\begin{Lem}\label{Lem:star-rho_omega_formula} The $(1,\,1)$-form $\rho_\omega$ uniquely determined by an arbitrary Hermitian metric $\omega$ on an $n$-dimensional complex manifold $X$ via property (\ref{eqn:rho_omega_def}) is given by the formula \begin{eqnarray}\label{eqn:star-rho_omega_formula}\star\rho_\omega = \frac{1}{n-1}\,\frac{\omega\wedge i\partial\bar\partial\omega_{n-2}}{\omega_n}\,\omega_{n-1} - i\partial\bar\partial\omega_{n-2}.\end{eqnarray}

\end{Lem}

Applying $\star$ again to (\ref{eqn:star-rho_omega_formula}), we get: \begin{eqnarray}\label{eqn:rho_omega_formula}\rho_\omega = \frac{1}{n-1}\,\frac{\omega\wedge i\partial\bar\partial\omega_{n-2}}{\omega_n}\,\omega + i\bar\partial^\star\partial^\star\omega_2.\end{eqnarray}

In particular, we see that with each Hermitian metric $\omega$ on $X$ we can associate the $C^\infty$ function $f_\omega:X\to\R$ defined by \begin{eqnarray}\label{eqn:function_f_omega}f_\omega:=\frac{\omega\wedge i\partial\bar\partial\omega_{n-2}}{\omega_n}.\end{eqnarray} The above computations also show that $f_\omega = (n-1)\,\Lambda_\omega(\rho_\omega)$.

\vspace{2ex}

 Rescaling $\omega$ by a positive constant $\lambda$ has the following obvious effect on the associated function: \begin{eqnarray}\label{eqn:rescaling_omega_f-omega}f_{\lambda\omega} = \frac{1}{\lambda}\,f_\omega,  \hspace{5ex} \lambda\in(0,\,\infty).\end{eqnarray} 

Note that the function $f_\omega$ vanishes identically if and only if the metric $\omega$ is $(n-2)$-Gauduchon in the sense of [FWW13]. 

\begin{Def}\label{Def:pluriclosed-star-split_metrics} Let $\omega$ be a Hermitian metric on an $n$-dimensional complex manifold $X$, let $\star=\star_\omega$ be the Hodge star operator induced by $\omega$ and let $\rho_\omega$ be the $(1,\,1)$-form on $X$ uniquely determined by $\omega$ via property (\ref{eqn:rho_omega_def}).

\vspace{1ex}

 (a)\, The metric  $\omega$ is said to be {\bf pluriclosed star split} if $\partial\bar\partial(\star\rho_\omega) = 0$. 

\vspace{1ex}

 (b)\, The metric $\omega$ is said to be {\bf closed star split} if $d(\star\rho_\omega) = 0$. 

\end{Def}

Thanks to formula (\ref{eqn:star-rho_omega_formula}), we see that the {\bf pluriclosed star split} condition on $\omega$ is equivalent to \begin{eqnarray}\label{eqn:pluriclosed-star-split_metrics_bis}\partial\bar\partial(f_\omega\,\omega_{n-1}) = 0,\end{eqnarray} while the {\bf closed star split} condition on $\omega$ is equivalent to \begin{eqnarray}\label{eqn:closed-star-split_metrics_bis}d(f_\omega\,\omega_{n-1}) = 0.\end{eqnarray}

The {\it closed star split} condition obviously implies the {\it pluriclosed star split} one. Both are implied by the $(n-2)$-Gauduchon condition of Fu-Wang-Wu [FWW13].

\vspace{3ex}

$\bullet$ Let us now recall a few facts from [Gau77a]. Consider the Laplace-type operator on functions: $$\Delta_\omega:=i\Lambda_\omega\bar\partial\partial:C^\infty(X,\,\C)\longrightarrow C^\infty(X,\,\C).$$ Its adjoint is the operator $\Delta_\omega^\star:C^\infty(X,\,\C)\longrightarrow C^\infty(X,\,\C)$ given by \begin{eqnarray}\label{eqn:Laplace-adjoint_functions}\Delta_\omega^\star(f) = i\star_\omega\bar\partial\partial(f\,\omega_{n-1}),\end{eqnarray} where $\star=\star_\omega$ is the Hodge star operator induced by $\omega$. This follows at once from the formulae $\partial^\star = -\star\bar\partial\star$, $\bar\partial^\star = -\star\partial\star$ and $\star\omega = \omega_{n-1}$. Now, Gauduchon proves, as a consequence of his Lemmas 1 and 2 of [Gau77a], the following key result.

\begin{Prop}([Gau77a])\label{Prop:existence_Gauduchon_conformal_3-4} Let $X$ be a compact connected complex manifold and let $\omega$ be a Hermitian metric on $X$. Let $f_0:X\longrightarrow\R$ be a $C^\infty$ function such that $f_0\in\ker\Delta^{\star}_{\omega}$. Then $$f_0> 0 \hspace{1ex} \mbox{on} \hspace{1ex} X \hspace{5ex} \mbox{or} \hspace{5ex} f_0< 0 \hspace{1ex} \mbox{on} \hspace{1ex} X  \hspace{5ex} \mbox{or} \hspace{5ex} f_0\equiv 0.$$

\end{Prop}

In our context, an immediate consequence of (\ref{eqn:pluriclosed-star-split_metrics_bis}) and (\ref{eqn:Laplace-adjoint_functions}) is the following

\begin{Prop}\label{Prop:pluriclosed-star-split_3-choices_f} Let $X$ be a compact connected complex manifold with $\mbox{dim}_\C X = n$ and let $\omega$ be a Hermitian metric on $X$. The following equivalence holds: \begin{eqnarray*}\omega \hspace{1ex}\mbox{is {\bf pluriclosed star split}} \iff \Delta_\omega^\star(f_\omega) =0.\end{eqnarray*}

\end{Prop} 

Proposition \ref{Prop:pluriclosed-star-split_3-choices_f_introd} stated in the introduction follows as a corollary of Propositions \ref{Prop:existence_Gauduchon_conformal_3-4} and \ref{Prop:pluriclosed-star-split_3-choices_f}.

\vspace{3ex}

We notice the following properties of pluriclosed star split metrics in relation to balanced and Gauduchon metrics. Part (iii) implies part (ii), but we state both of them since the proofs we give are slightly different.

\begin{Prop}\label{Prop:pluriclosed-star-split_properties} Let $(X,\,\omega)$ be a complex Hermitian manifold with $\mbox{dim}_\C X=n$. 

  (i)\, If the function $f_\omega$ is a {\bf non-zero constant}, the metric $\omega$ is {\bf pluriclosed star split} if and only if it is {\bf Gauduchon}.

\vspace{1ex}

(ii)\, Suppose $X$ is {\bf compact} and connected. If the metric $\omega$ is {\bf balanced}, then it is {\bf pluriclosed star split} if and only if the function $f_\omega$ is {\bf constant}.

\vspace{1ex}

(iii)\, Suppose $X$ is {\bf compact} and connected. If the metric $\omega$ is {\bf Gauduchon}, then it is {\bf pluriclosed star split} if and only if the function $f_\omega$ is {\bf constant}.

\end{Prop}

\noindent {\it Proof.} Part (i) follows obviously from (\ref{eqn:pluriclosed-star-split_metrics_bis}). To prove part (ii), we note that $\omega$ being balanced implies the first of the following equalities: \begin{eqnarray*}i\partial\bar\partial(f_\omega\,\omega_{n-1}) = i\partial\bar\partial f_\omega\wedge\omega_{n-1} = \Lambda_\omega(i\partial\bar\partial f_\omega)\,\omega_n = -\Delta_\omega(f_\omega)\,\omega_n.\end{eqnarray*} Therefore, thanks to (\ref{eqn:pluriclosed-star-split_metrics_bis}), $\omega$ is pluriclosed star split if and only if $\Delta_\omega(f_\omega) = 0$ on $X$. By the maximum principle, this is equivalent to $f_\omega$ being constant since $X$ is compact and the Laplacian $-\Delta_\omega= \Lambda_\omega(i\partial\bar\partial)$ is elliptic of order two with no zero-th order terms.

We now prove part (iii). Suppose the metric $\omega$ is Gauduchon.

``$\implies$'' If $\omega$ is pluriclosed star split, then $f_\omega>0$ on $X$, or $f_\omega<0$ on $X$, or 
$f_\omega\equiv 0$. In the last case, $f_\omega$ is indeed constant. In the case where $f_\omega>0$ on $X$, the pluriclosed star split condition (\ref{eqn:pluriclosed-star-split_metrics_bis}) translates to the metric $(f_\omega)^{1/(n-1)}\,\omega$ being Gauduchon (and also, obviously, conformally equivalent to $\omega$). Since the metric $\omega$ is already Gauduchon, this implies that the function $(f_\omega)^{1/(n-1)}$, hence also $f_\omega$, is constant. Indeed, Gauduchon's main theorem (th\'eor\`eme 1) in [Gau77a] stipulates that in every conformal class of Hermitian metrics on $X$ there exists a unique (up to multiplicative positive constants) Gauduchon metric. In the case where $f_\omega<0$ on $X$, $-f_\omega>0$ and the above argument can be repeated with $-f_\omega$ in place of $f_\omega$.   

``$\Longleftarrow$'' If the function $f_\omega$ is constant, then either $f_\omega$ is a non-zero constant, in which case $\omega$ is pluriclosed star split by (i); or $f_\omega$ vanishes identically, in which case (\ref{eqn:pluriclosed-star-split_metrics_bis}) holds in an obvious way, so $\omega$ is again pluriclosed star split.  \hfill $\Box$

\vspace{3ex}

Formula (\ref{eqn:star-rho_omega_formula}) shows that if the metric $\omega$ is both {\it Gauduchon} and {\it pluriclosed star split}, then \begin{eqnarray}\label{eqn:Gauduchon-pluriclosed-star-split_A}[\star\rho_\omega]_A = \frac{1}{n-1}\,f_\omega\,[\omega_{n-1}]_A,\end{eqnarray} while if $\omega$ is both a {\it balanced} and a {\it closed star split} metric, then \begin{eqnarray}\label{eqn:Gauduchon-pluriclosed-star-split_BC}[\star\rho_\omega]_{BC} = \frac{1}{n-1}\,f_\omega\,[\omega_{n-1}]_{BC}.\end{eqnarray} In both of these cases, $f_\omega$ is constant thanks to Proposition \ref{Prop:pluriclosed-star-split_properties}.

\vspace{3ex}

$\bullet$ The following immediate observation will come in handy later on.

\begin{Lem}\label{Lem:star-trace-11n-1n-1} Let $\gamma$ be a positive definite $(1,\,1)$-form and $\Gamma$ a real $(n-1,\,n-1)$-form on a complex manifold $X$ with $\mbox{dim}_\C X=n$. The following equality holds: \begin{eqnarray}\label{eqn:star-trace-11n-1n-1}\gamma\wedge\Gamma = \star_\gamma\Gamma\wedge\gamma_{n-1},\end{eqnarray} where $\star_\gamma$ is the Hodge star operator induced by the Hermitian metric $\gamma$.

\end{Lem}

\noindent {\it Proof.} Using the properties of $\star_\gamma$, we have: \begin{eqnarray*}\gamma\wedge\Gamma & = & \gamma\wedge\star_\gamma\star_\gamma\Gamma = \langle\gamma,\,\star_\gamma\Gamma\rangle_\gamma\,\gamma_n = \langle\star_\gamma\star_\gamma\gamma,\,\star_\gamma\Gamma\rangle_\gamma\,\gamma_n \\
 & \stackrel{(a)}{=} & \langle\gamma_{n-1},\,\Gamma\rangle_\gamma\,\gamma_n = \gamma_{n-1}\wedge\star_\gamma\Gamma,\end{eqnarray*} where (a) follows from the fact that $\star_\gamma$ is an isometry for the pointwise inner product induced by $\gamma$ and from the standard identity $\star_\gamma\gamma = \gamma_{n-1}$.  \hfill $\Box$

\vspace{2ex}

Independently, we make the following observation whose proof is very similar to the one of Proposition 1.1 in [Pop15]. The latter, weaker, part appeared in [MT01] with a different proof. 

\begin{Prop}\label{Prop:astheno-balanced-K} Let $\omega$ be a Hermitian metric on a compact complex manifold $X$ with $\mbox{dim}_\C X =n\geq 3$. If $\omega$ is both {\bf balanced} and {\bf $(n-2)$-Gauduchon}, it is {\bf K\"ahler}. 

In particular, if $\omega$ is both {\bf balanced} and {\bf astheno-K\"ahler}, it is {\bf K\"ahler}.  

\end{Prop}

\noindent {\it Proof.} Taking $\gamma=\omega$ and $\Gamma= i\partial\bar\partial\omega_{n-2}$ in Lemma \ref{Lem:star-trace-11n-1n-1}, we get the first equality below: \begin{eqnarray}\label{eqn:astheno_equiv_1}\nonumber\omega\wedge i\partial\bar\partial\omega_{n-2} & = & \star(i\partial\bar\partial\omega_{n-2})\wedge\omega_{n-1} = -i\,(-\star\bar\partial\star)(\star\partial\omega_{n-2})\wedge\omega_{n-1} \\
& = & -i\,\partial^\star\bigg(\star(\omega_{n-3}\wedge\partial\omega)\bigg)\wedge\omega_{n-1},\end{eqnarray} where we also used the standard identities $\star\star = \pm\mbox{Id}$ on forms of even, resp. odd, degree, and $\partial^\star = -\star\bar\partial\star$.

On the other hand, the balanced hypothesis on $\omega$ is well known to be equivalent to $\partial\omega$ being primitive. Indeed, $\partial\omega_{n-1} = \omega_{n-2}\wedge\partial\omega$ and, since $\partial\omega$ is a $3$-form, the condition $\omega_{n-2}\wedge\partial\omega = 0$ expresses the primitivity of $\partial\omega$, by definition. Thus, using the standard formula (\ref{eqn:prim-form-star-formula-gen}) for primitive forms of bidegree $(p,\,q)=(2,\,1)$, we get: $\star(\partial\omega) = i\,\omega_{n-3}\wedge\partial\omega$, hence \begin{eqnarray}\label{eqn:balanced_equiv_1}\star(\omega_{n-3}\wedge\partial\omega) = i\,\partial\omega\end{eqnarray} when $\omega$ is balanced. 

Putting (\ref{eqn:astheno_equiv_1}) and (\ref{eqn:balanced_equiv_1}) together, we infer that \begin{eqnarray}\label{eqn:balanced_conseq_proof-astheno}\omega\wedge i\partial\bar\partial\omega_{n-2} = \partial^\star\partial\omega\wedge\omega_{n-1} = \Lambda_\omega(\partial^\star\partial\omega)\,\omega_n\end{eqnarray} whenever $\omega$ is balanced. In particular, if besides being balanced, $\omega$ is also $(n-2)$-Gauduchon (i.e. $\omega\wedge i\partial\bar\partial\omega_{n-2}=0$), then $\Lambda_\omega(\partial^\star\partial\omega)=0$. Taking the $L^2_\omega$ inner product with the constant function $1$, we get in that case: \begin{eqnarray*}0 = \langle\langle\Lambda_\omega(\partial^\star\partial\omega),\,1\rangle\rangle_\omega = \langle\langle\partial^\star\partial\omega,\,\omega\rangle\rangle_\omega = ||\partial\omega||^2_\omega.\end{eqnarray*} This implies $\partial\omega=0$, hence $\omega$ is K\"ahler.  \hfill $\Box$

\vspace{3ex}

This discussion is summed up in the following

\begin{Cor}\label{Cor:astheno-balanced-K-rho} Suppose $\omega$ is a {\bf balanced} metric on a compact complex manifold $X$ with $\mbox{dim}_\C X =n\geq 3$. The following equivalences hold: \begin{eqnarray*}\label{eqn:astheno-balanced-K-rho_0}\omega \hspace{1ex}\mbox{is K\"ahler}\hspace{1ex} \iff \omega \hspace{1ex}\mbox{is}\hspace{1ex} (n-2)\mbox{-Gauduchon}\end{eqnarray*}\begin{eqnarray*}\label{eqn:astheno-balanced-K-rho_1}\omega \hspace{1ex}\mbox{is K\"ahler}\hspace{1ex} \iff \omega \hspace{1ex}\mbox{is astheno-K\"ahler}\end{eqnarray*} \begin{eqnarray*}\label{eqn:astheno-balanced-K-rho_2}\omega \hspace{1ex}\mbox{is K\"ahler}\hspace{1ex} \iff \rho_\omega=0\end{eqnarray*} \begin{eqnarray*}\label{eqn:astheno-balanced-K-rho_3}\omega \hspace{1ex}\mbox{is K\"ahler}\hspace{1ex} \iff \star\rho_\omega=0.\end{eqnarray*}

\end{Cor}

\vspace{3ex}

$\bullet$ {\bf Proof of (i) of Theorem \ref{The:bal-pluri_star-split_f-not-neg}.} It is easily obtained from the proof of Proposition \ref{Prop:astheno-balanced-K} in the following way. Under the assumptions of (i) of Theorem \ref{The:bal-pluri_star-split_f-not-neg}, the function $f_\omega$ is necessarily constant by (ii) of Proposition \ref{Prop:pluriclosed-star-split_properties}.

If $f_\omega$ were a negative constant $c$, we would have $\Lambda_\omega(\partial^\star\partial\omega) = c<0$ by (\ref{eqn:balanced_conseq_proof-astheno}), hence also \begin{eqnarray*}c\,\int_X\omega_n = \langle\langle\Lambda_\omega(\partial^\star\partial\omega),\,1\rangle\rangle_\omega = \langle\langle\partial^\star\partial\omega,\,\omega\rangle\rangle_\omega = ||\partial\omega||^2_\omega,\end{eqnarray*} which is impossible if $c<0$. \hfill $\Box$

\vspace{2ex}

However, it may happen that $\omega$ be balanced and pluriclosed star split with the function $f_\omega$ a positive constant. Examples include the standard metric $\omega$ on the Iwasawa manifold, where $f_\omega\equiv 1$ (see Proposition \ref{Prop:Iwasawa_examp}), any Nakamura manifold, where $f_\omega\equiv 2$ (see Proposition \ref{Prop:Nakamura_examp}).

\vspace{2ex}

$\bullet$ {\bf Proof of (ii) of Theorem \ref{The:bal-pluri_star-split_f-not-neg}.} That the function $f_\omega$ may be a negative constant when the pluriclosed star split metric $\omega$ is not balanced is shown by the examples of the small deformations $X_t$ of the Calabi-Eckmann manifold $X_0=(S^3\times S^3,\,J_{CE})$ corresponding to complex numbers $t$ close enough to $0$ for which $\mbox{Im}\,(t)<0$ (see Proposition \ref{Prop::small-def-S3S3_examp}).  \hfill $\Box$

\vspace{3ex}

$\bullet$ Let us now return to the function $f_\omega:X\to\R$ associated via (\ref{eqn:function_f_omega}) with an arbitrary Hermitian metric $\omega$. We get at once that \begin{eqnarray*}\omega\wedge\star\rho_\omega = \frac{1}{n-1}\,\omega\wedge i\partial\bar\partial\omega_{n-2} =  \frac{1}{n-1}\,f_\omega\,\omega_n.\end{eqnarray*} Hence, integrating and applying Stokes, we get: \begin{eqnarray*}\int_X\omega\wedge\star\rho_\omega = -\frac{1}{n-1}\,\int_Xi\partial\omega\wedge\bar\partial\omega\wedge\omega_{n-3}.\end{eqnarray*} Meanwhile, if $\omega$ is {\it balanced}, $\bar\partial\omega$ is primitive (as seen above for $\partial\omega$), so the standard formula (\ref{eqn:prim-form-star-formula-gen}) for primitive forms of bidegree $(p,\,q)=(1,\,2)$ gives: $\star(\bar\partial\omega) = -i\,\omega_{n-3}\wedge\bar\partial\omega$. Together with the last integral formula, this leads to \begin{eqnarray*}\int_X\omega\wedge\star\rho_\omega = \frac{1}{n-1}\,\int_X\partial\omega\wedge\star(\bar\partial\omega) = \frac{1}{n-1}\,||\partial\omega||^2_\omega,\end{eqnarray*} if $\omega$ is {\it balanced}, where $||\,\,||_\omega$ is the $L^2_\omega$-norm.

We have thus proved the following

\begin{Prop}\label{Prop:balanced_f-omega_del-omega} Suppose $\omega$ is a {\bf balanced} metric on a compact complex manifold $X$ with $\mbox{dim}_\C X =n\geq 3$. The following equality holds: \begin{eqnarray}\label{eqn:balanced_f-omega_del-omega}\int_Xf_\omega\,\omega_n = ||\partial\omega||^2_\omega.\end{eqnarray} 

In particular, $\omega$ is K\"ahler if and only if $\int_Xf_\omega\,\omega_n=0$. If $f_\omega$ is constant, then $\omega$ is K\"ahler if and only if $f_\omega=0$.

\end{Prop}

Proposition \ref{Prop:balanced_f-omega_del-omega_introd} is an immediate corollary of Proposition \ref{Prop:balanced_f-omega_del-omega}.

\subsection{Examples of pluriclosed star split and closed star split Hermitian metrics}\label{subsection:examples_pluriclosed-star-split-metrics} All the examples we now point out are Gauduchon metrics $\omega$ with $f_\omega$ constant (see (i) of Proposition \ref{Prop:pluriclosed-star-split_properties}). However, the value of the constant will vary from case to case.    

\subsubsection{Example of the Iwasawa manifold}\label{subsubsection:Iwasawa_examp} For the Iwasawa manifold $X$, $\mbox{dim}_\C X=3$ and the cohomology of $X$ is completely determined by three holomorphic $(1,\,0)$-forms $\alpha$, $\beta$, $\gamma$ (so, $\bar\partial\alpha = \bar\partial\beta = \bar\partial\gamma = 0$) that satisfy the structure equations: \begin{eqnarray}\label{eqn:Iwasawa_structure-eq}\partial\alpha = \partial\beta = 0  \hspace{3ex}\mbox{and} \hspace{3ex} \partial\gamma = -\alpha\wedge\beta\neq 0.\end{eqnarray}

We will show that the standard Hermitian metric \begin{eqnarray*}\omega = i\alpha\wedge\bar\alpha + i\beta\wedge\bar\beta + i\gamma\wedge\bar\gamma,\end{eqnarray*} which is known to be balanced, is also pluriclosed star split. We know from (ii) of Proposition \ref{Prop:pluriclosed-star-split_properties} that the function $f_\omega$ will then be constant. We will compute this constant.

Since $n-2=1$, we have to compute $i\partial\bar\partial\omega$. We get $\bar\partial\omega = i\gamma\wedge\bar\alpha\wedge\bar\beta$, $i\partial\bar\partial\omega = i\alpha\wedge\bar\alpha\wedge i\beta\wedge\bar\beta$ and \begin{eqnarray*}\omega\wedge(i\partial\bar\partial\omega) = i\alpha\wedge\bar\alpha\wedge i\beta\wedge\bar\beta \wedge i\gamma\wedge\bar\gamma = \omega_3.\end{eqnarray*} In particular, $f_\omega = \omega\wedge(i\partial\bar\partial\omega)/\omega_3 =1$.

On the other hand, \begin{eqnarray*}\omega_2 = i\alpha\wedge\bar\alpha\wedge i\beta\wedge\bar\beta + i\alpha\wedge\bar\alpha\wedge i\gamma\wedge\bar\gamma + i\beta\wedge\bar\beta\wedge i\gamma\wedge\bar\gamma,\end{eqnarray*} (so, we see that $\partial\omega_2 = 0$, which means that $\omega$ is balanced, as is well known) and we get: \begin{eqnarray*}\star\rho_\omega = \frac{1}{2}\,f_\omega\,\omega_2 - i\partial\bar\partial\omega = \frac{1}{2}\,\bigg(i\alpha\wedge\bar\alpha\wedge i\gamma\wedge\bar\gamma + i\beta\wedge\bar\beta\wedge i\gamma\wedge\bar\gamma - i\alpha\wedge\bar\alpha\wedge i\beta\wedge\bar\beta\bigg).\end{eqnarray*} In particular, we see that $\partial(\star\rho_\omega) = 0$, hence also $\bar\partial(\star\rho_\omega) = 0$ because $\star\rho_\omega$ is a real form. This, of course, implies $\partial\bar\partial(\star\rho_\omega) = 0$, but we even have $d(\star\rho_\omega) = 0$. Meanwhile, from $\rho_\omega = \star\star\rho_\omega$ and the above formula for $\star\rho_\omega$, we get: \begin{eqnarray*}\rho_\omega = \frac{1}{2}\,\bigg(i\alpha\wedge\bar\alpha + i\beta\wedge\bar\beta - i\gamma\wedge\bar\gamma\bigg).\end{eqnarray*}
 
The conclusions of this computation are summed up in the following

\begin{Prop}\label{Prop:Iwasawa_examp} The balanced metric $\omega = i\alpha\wedge\bar\alpha + i\beta\wedge\bar\beta + i\gamma\wedge\bar\gamma$ on the {\bf Iwasawa manifold} is {\bf closed star split}, hence also {\bf pluriclosed star split}, and its associated function $f_\omega$ is {\bf constant} equal to $1$.

Moreover, the eigenvalues of the real $(2,\,2)$-form $\star\rho_\omega$ with respect to $\omega$ are $1, 1, -1$.

\end{Prop}

\subsubsection{Example of the Nakamura manifolds}\label{subsubsection:Nakamura_examp} For the Nakamura manifolds $X$, $\mbox{dim}_\C X=3$ and the cohomology of $X$ is completely determined by three holomorphic $(1,\,0)$-forms $\varphi_1$, $\varphi_2$, $\varphi_3$ (so, $\bar\partial\varphi_1 = \bar\partial\varphi_2 = \bar\partial\varphi_3 = 0$) that satisfy the structure equations: \begin{eqnarray}\label{eqn:Nakamura_structure-eq}\partial\varphi_1 = 0,  \hspace{3ex} \partial\varphi_2 = \varphi_1\wedge\varphi_2 \hspace{3ex}\mbox{and} \hspace{3ex} \partial\varphi_3 = -\varphi_1\wedge\varphi_3.\end{eqnarray}

We will show that the standard Hermitian metric \begin{eqnarray*}\omega = i\varphi_1\wedge\bar\varphi_1 + i\varphi_2\wedge\bar\varphi_2 + i\varphi_3\wedge\bar\varphi_3,\end{eqnarray*} which is known to be balanced, is also pluriclosed star split. In particular, the function $f_\omega$ will be constant (see (ii) of Proposition \ref{Prop:pluriclosed-star-split_properties}) and we will compute this constant.

We get: \begin{eqnarray*}\bar\partial\omega & = & -i\varphi_2\wedge(\bar\varphi_1\wedge\bar\varphi_2) -i\varphi_3\wedge(-\bar\varphi_1\wedge\bar\varphi_3) = \bar\varphi_1\wedge i\varphi_2\wedge\bar\varphi_2 - \bar\varphi_1\wedge i\varphi_3\wedge\bar\varphi_3\end{eqnarray*} and \begin{eqnarray*}i\partial\bar\partial\omega & = & -i\bar\varphi_1\wedge i(\varphi_1\wedge\varphi_2)\wedge\bar\varphi_2 + i\bar\varphi_1\wedge i(-\varphi_1\wedge\varphi_3)\wedge\bar\varphi_3\\
 & = & i\varphi_1\wedge\bar\varphi_1\wedge i\varphi_2\wedge\bar\varphi_2 + i\varphi_1\wedge\bar\varphi_1\wedge i\varphi_3\wedge\bar\varphi_3.\end{eqnarray*} Hence \begin{eqnarray*}\omega\wedge(i\partial\bar\partial\omega) = 2\,\omega_3.\end{eqnarray*} In particular, $f_\omega = \omega\wedge(i\partial\bar\partial\omega)/\omega_3 =2$.

On the other hand, \begin{eqnarray*}\omega_2 = i\varphi_1\wedge\bar\varphi_1\wedge i\varphi_2\wedge\bar\varphi_2 + i\varphi_1\wedge\bar\varphi_1\wedge i\varphi_3\wedge\bar\varphi_3 + i\varphi_2\wedge\bar\varphi_2\wedge i\varphi_3\wedge\bar\varphi_3,\end{eqnarray*} (so, we see that $\partial\omega_2 = i(\varphi_1\wedge\varphi_2)\wedge\bar\varphi_2\wedge i\varphi_3\wedge\bar\varphi_3 + i\varphi_2\wedge\bar\varphi_2\wedge i(-\varphi_1\wedge\varphi_3)\wedge\bar\varphi_3 = 0$, which means that $\omega$ is balanced, as is well known) and we get: \begin{eqnarray*}\star\rho_\omega = \frac{1}{2}\,f_\omega\,\omega_2 - i\partial\bar\partial\omega = i\varphi_2\wedge\bar\varphi_2\wedge i\varphi_3\wedge\bar\varphi_3.\end{eqnarray*}

 In particular, we see that \begin{eqnarray*}\partial(\star\rho_\omega) = i(\varphi_1\wedge\varphi_2)\wedge\bar\varphi_2\wedge i\varphi_3\wedge\bar\varphi_3 + i\varphi_2\wedge\bar\varphi_2\wedge i(-\varphi_1\wedge\varphi_3)\wedge\bar\varphi_3 = 0,\end{eqnarray*} hence also $\bar\partial(\star\rho_\omega) = 0$ because $\star\rho_\omega$ is a real form. This, of course, implies $\partial\bar\partial(\star\rho_\omega) = 0$, but we even have $d(\star\rho_\omega) = 0$. Meanwhile, from $\rho_\omega = \star\star\rho_\omega$ and the above formula for $\star\rho_\omega$, we get: \begin{eqnarray*}\rho_\omega = i\varphi_1\wedge\bar\varphi_1.\end{eqnarray*}

The conclusions of this computation are summed up in the following

\begin{Prop}\label{Prop:Nakamura_examp} The balanced metric $\omega = i\varphi_1\wedge\bar\varphi_1 + i\varphi_2\wedge\bar\varphi_2 + i\varphi_3\wedge\bar\varphi_3$ on any {\bf Nakamura manifold} is {\bf closed star split}, hence also {\bf pluriclosed star split}, and its associated function $f_\omega$ is {\bf constant} equal to $2$.

Moreover, the eigenvalues of the real $(2,\,2)$-form $\star\rho_\omega$ with respect to $\omega$ are $1, 0, 0$.

\end{Prop}

\subsubsection{Example of the small deformations of the Iwasawa manifold}\label{subsubsection:small-def-Iwasawa_examp} The small deformations $X_t$ lying in Nakamura's class $(i)$ of the Iwasawa manifold $X$ have the same properties as $X$, so we have implicitly discussed them in $\S.$\ref{subsubsection:Iwasawa_examp}. 

We now concentrate on the small deformations $X_t$ lying in one of Nakamura's classes $(ii)$ or $(iii)$. Their cohomology is completely determined by three smooth $(1,\,0)$-forms $\alpha_t$, $\beta_t$, $\gamma_t$ satisfying the {\bf structure equations} (cf. [Ang11, $\S.4.3$]): 
\begin{eqnarray*}\label{eq:structure}
d\alpha_t & = & d\beta_t = 0, \\
\partial_t\gamma_t&=& \sigma_{12}(t)\,\alpha_t\wedge\beta_t, \\
\bar\partial_t\gamma_t&=& \sigma_{1\bar{1}}(t)\,\alpha_t\wedge\bar\alpha_t + \sigma_{1\bar{2}}(t)\,\alpha_t\wedge\bar\beta_t + \sigma_{2\bar{1}}(t)\,\beta_t\wedge\bar\alpha_t + \sigma_{2\bar{2}}(t)\,\beta_t\wedge\bar\beta_t,
  \end{eqnarray*}
where $\sigma_{12}$ and $\sigma_{i\bar{j}}$ are $C^{\infty}$ functions of $t\in\Delta$ that depend only on $t$ (so $\sigma_{12}(t)$ and $\sigma_{i\bar{j}}(t)$ are complex numbers for every fixed $t$ in the parameter space $\Delta$) and satisfy $\sigma_{12}(0) = -1$ and $\sigma_{i\bar{j}}(0) = 0$ for all $i, j$.

Now, for every $t\in\Delta$ close to $0$, the $J_t$-$(1,\,1)$-form 
\begin{equation*}\label{eqn:omega_t_def}\omega_t:=i\alpha_t\wedge\overline{\alpha}_t + i\beta_t\wedge\overline{\beta}_t + i\gamma_t\wedge\overline{\gamma}_t\end{equation*} 
\noindent is positive definite, hence it defines a Hermitian metric on $X_t$ that varies in a $C^{\infty}$ way with $t$. 

Computing, we get: \begin{eqnarray*}\bar\partial_t\omega_t & = & i\bar\partial_t\gamma_t\wedge\bar\gamma_t - i\gamma_t\wedge\overline{\partial_t\gamma_t} \\
 & = & \sigma_{1\bar{1}}(t)\,i\alpha_t\wedge\bar\alpha_t\wedge\bar\gamma_t + \sigma_{1\bar{2}}(t)\,i\alpha_t\wedge\bar\beta_t\wedge\bar\gamma_t + \sigma_{2\bar{1}}(t)\,i\beta_t\wedge\bar\alpha_t\wedge\bar\gamma_t + \sigma_{2\bar{2}}(t)\,i\beta_t\wedge\bar\beta_t\wedge\bar\gamma_t\\
& - & \overline{\sigma_{12}(t)}\,i\gamma_t\wedge\bar\alpha_t\wedge\bar\beta_t,\end{eqnarray*} hence \begin{eqnarray*}i\partial_t\bar\partial_t\omega_t = A(t)\,i\alpha_t\wedge\bar\alpha_t\wedge i\beta_t\wedge\bar\beta_t,\end{eqnarray*} where \begin{eqnarray}\label{eqn:At}A(t):=|\sigma_{12}(t)|^2 + |\sigma_{2\bar{1}}(t)|^2 + |\sigma_{1\bar{2}}(t)|^2 - 2\,\mbox{Re}\,\bigg(\sigma_{1\bar{1}}(t)\,\overline{\sigma_{2\bar{2}}(t)}\bigg).\end{eqnarray} 

Thus, we get: \begin{eqnarray*}\omega_t\wedge i\partial_t\bar\partial_t\omega_t = A(t)\,i\alpha_t\wedge\bar\alpha_t\wedge i\beta_t\wedge\bar\beta_t\wedge i\gamma_t\wedge\bar\gamma_t = A(t)\,(\omega_t)_3.\end{eqnarray*} In particular, $f_{\omega_t} = \omega_t\wedge i\partial_t\bar\partial_t\omega_t/(\omega_t)_3 = A(t)$, which is a real constant for every fixed $t$.

On the other hand, we get: \begin{eqnarray*}(\omega_t)_2 = i\alpha_t\wedge\bar\alpha_t\wedge i\beta_t\wedge\bar\beta_t + i\alpha_t\wedge\bar\alpha_t\wedge i\gamma_t\wedge\bar\gamma_t + i\beta_t\wedge\bar\beta_t\wedge i\gamma_t\wedge\bar\gamma_t,\end{eqnarray*} so we see that $\partial_t\bar\partial_t(\omega_t)_2 = 0$, proving that $\omega_t$ is a Gauduchon metric for every $t$ close to $0$. Thanks to (i) of Proposition \ref{Prop:pluriclosed-star-split_properties}, we conclude that $\omega_t$ is pluriclosed star split.

Explicitly, we further have: \begin{eqnarray*}\star\rho_{\omega_t} = \frac{A(t)}{2}\,(\omega_t)_2 - i\partial_t\bar\partial_t\omega_t = \frac{A(t)}{2}\,\bigg(i\alpha_t\wedge\bar\alpha_t\wedge i\gamma_t\wedge\bar\gamma_t + i\beta_t\wedge\bar\beta_t\wedge i\gamma_t\wedge\bar\gamma_t - i\alpha_t\wedge\bar\alpha_t\wedge i\beta_t\wedge\bar\beta_t\bigg).\end{eqnarray*} 

We point out that, unlike in the cases of the Iwasawa and Nakamura manifolds, the form $\star\rho_{\omega_t}$ is not $d$-closed. To see this, we notice that $i\alpha_t\wedge\bar\alpha_t\wedge i\beta_t\wedge\bar\beta_t$ is $d$-closed and that \begin{eqnarray*}\bar\partial_t(i\alpha_t\wedge\bar\alpha_t\wedge i\gamma_t\wedge\bar\gamma_t) & = & \sigma_{2\bar{2}}(t)\,i\alpha_t\wedge\bar\alpha_t\wedge i\beta_t\wedge\bar\beta_t\wedge\bar\gamma_t \\
 \bar\partial_t(i\beta_t\wedge\bar\beta_t\wedge i\gamma_t\wedge\bar\gamma_t) & = & \sigma_{1\bar{1}}(t)\,i\alpha_t\wedge\bar\alpha_t\wedge i\beta_t\wedge\bar\beta_t\wedge\bar\gamma_t.\end{eqnarray*} Hence \begin{eqnarray}\label{eqn:dbar_t-rho_I-def}\bar\partial_t(\star\rho_{\omega_t}) = \frac{A(t)}{2}\,\bigg(\sigma_{1\bar{1}}(t) + \sigma_{2\bar{2}}(t)\bigg)\,i\alpha_t\wedge\bar\alpha_t\wedge i\beta_t\wedge\bar\beta_t\wedge\bar\gamma_t\end{eqnarray} for all $t$ sufficiently close to $0$ in one of Nakamura's classes (ii) or (iii). In particular, we see that there is no reason for $\bar\partial_t(\star\rho_{\omega_t})$ to vanish when $t\neq 0$ since $\sigma_{1\bar{1}}(t) + \sigma_{2\bar{2}}(t)$ need not vanish in that case.

The conclusions of this computation are summed up in the following

\begin{Prop}\label{Prop:small-def-Iwasawa_examp} For every $t\neq 0$ sufficiently close to $0$ in one of Nakamura's classes (ii) or (iii), the Gauduchon metric $\omega_t = i\alpha_t\wedge\bar\alpha_t + i\beta_t\wedge\bar\beta_t + i\gamma_t\wedge\bar\gamma_t$ on the corresponding {\bf small deformation} $X_t$ of the Iwasawa manifold $X=X_0$ is {\bf pluriclosed star split} and its associated function $f_{\omega_t}$ is {\bf constant} equal to $A(t)$ defined in (\ref{eqn:At}).

Moreover, the $\partial\bar\partial$-closed real $(2,\,2)$-form $\star\rho_{\omega_t}$ need not be $d$-closed and its eigenvalues with respect to $\omega$ are $1, 1, -1$.

\end{Prop}

 \subsubsection{Example of the $5$-dimensional Iwasawa manifold}\label{subsubsection:5-Iwasawa_examp} The cohomology of the $5$-dimensional Iwasawa manifold $X=I^{(5)}$ is completely determined by five holomorphic $(1,\,0)$-forms $\varphi_j$, with $j=1,\dots , 5$, (so $\bar\partial\varphi_j = 0$) that satisfy the structure equations: \begin{eqnarray}\label{eqn:5-Iwasawa_structure-eq}\partial\varphi_1 = \partial\varphi_2 = 0,  \hspace{3ex} \partial\varphi_3 = \varphi_1\wedge\varphi_2, \hspace{3ex} \partial\varphi_4 = \varphi_1\wedge\varphi_3, \hspace{3ex} \partial\varphi_5 = \varphi_2\wedge\varphi_3.\end{eqnarray} 

The standard Hermitian metric on $X$ is defined by $\omega = \sum_{j=1}^5 i\varphi_j\wedge\bar\varphi_j$. Since $n=5$, we need to compute $i\partial\bar\partial\omega_3$. We get successively: \begin{eqnarray*}\omega_3 = \sum\limits_{1\leq j<k<l\leq 5}i\varphi_j\wedge\bar\varphi_j\wedge i\varphi_k\wedge\bar\varphi_k\wedge i\varphi_l\wedge\bar\varphi_l,\end{eqnarray*} \begin{eqnarray*}\bar\partial\omega_3 & = & -i\varphi_1\wedge\bar\varphi_1\wedge i\varphi_4\wedge\bar\varphi_4\wedge i\varphi_5\wedge\bar\varphi_2\wedge\bar\varphi_3 -i\varphi_2\wedge\bar\varphi_2\wedge i\varphi_4\wedge\bar\varphi_1\wedge\bar\varphi_3\wedge i\varphi_5\wedge\bar\varphi_5 \\
 & - & i\varphi_3\wedge\bar\varphi_1\wedge\bar\varphi_2\wedge i\varphi_4\wedge\bar\varphi_4\wedge i\varphi_5\wedge\bar\varphi_5\end{eqnarray*} and \begin{eqnarray*}i\partial\bar\partial\omega_3 = \widehat{i\varphi_3\wedge\bar\varphi_3} + \widehat{i\varphi_4\wedge\bar\varphi_4} + \widehat{i\varphi_5\wedge\bar\varphi_5},\end{eqnarray*} where $\widehat{i\varphi_j\wedge\bar\varphi_j}$ stands for the product of all the $(1,\,1)$-forms $\widehat{i\varphi_k\wedge\bar\varphi_k}$ with $k\in\{1,\dots 5\}\setminus\{j\}$.

From this, we get: $\omega\wedge i\partial\bar\partial\omega_3 = 3\,\omega_5$, hence $f_\omega = \omega\wedge i\partial\bar\partial\omega_3/\omega_5 = 3$.

Meanwhile, \begin{eqnarray*}\omega_4 = \sum\limits_{j=1}^5\widehat{i\varphi_j\wedge\bar\varphi_j},\end{eqnarray*} so we see that $\bar\partial\omega_4 = 0$ (which means that the metric $\omega$ is balanced) and we get: \begin{eqnarray*}\star\rho_\omega & = & \frac{1}{4}\,f_\omega\,\omega_4 - i\partial\bar\partial\omega_3 = \frac{3}{4}\,\omega_4 - i\partial\bar\partial\omega_3\\
 & = & \frac{3}{4}\,\bigg(\widehat{i\varphi_1\wedge\bar\varphi_1} + \widehat{i\varphi_2\wedge\bar\varphi_2}\bigg) - \frac{1}{4}\,\bigg(\widehat{i\varphi_3\wedge\bar\varphi_3} + \widehat{i\varphi_4\wedge\bar\varphi_4} + \widehat{i\varphi_5\wedge\bar\varphi_5}\bigg).\end{eqnarray*} In particular, $\bar\partial(\star\rho_\omega) = 0$.

The conclusions of this computation are summed up in the following

\begin{Prop}\label{Prop:5-Iwasawa_examp} The balanced metric $\omega = \sum_{j=1}^5 i\varphi_j\wedge\bar\varphi_j$ on the {\bf $5$-dimensional Iwasawa manifold} $I^{(5)}$ is {\bf closed star split}, hence also {\bf pluriclosed star split}, and its associated function $f_\omega$ is {\bf constant} equal to $3$.

Moreover, the eigenvalues of the real $(4,\,4)$-form $\star\rho_\omega$ with respect to $\omega$ are $3/4$, $3/4$, $-1/4$, $-1/4$, $-1/4$.

\end{Prop}

\subsubsection{Example of the Calabi-Eckmann manifold $(S^3\times S^3,\, J_{CE})$ and its small deformations}\label{subsubsection:small-def-S3S3_examp} As is well known, Calabi and Eckmann defined in [CE53] a complex structure $J_{CE}$ on the $C^\infty$ manifold $S^3\times S^3$, where $S^3$ is the $3$-sphere. In [TT17], Tardini and Tomassini studied a one-dimensional deformation $(X_t,\,J_t)$, with $t$ varying in a small disc $D$ about $0$ in $\C$, of the Calabi-Eckmann complex manifold $X=X_0= (S^3\times S^3,\, J_{CE})$, with $\mbox{dim}_\C X =3$, whose cohomology is determined by three smooth $1$-forms $\varphi_t^1$, $\varphi_t^2$, $\varphi_t^3$ satisfying the following structure equations: \begin{align*}\partial_t\varphi_t^1 & = \frac{\bar{t} + 1}{1-|t|^2}\,i\varphi_t^1\wedge\varphi_t^3, &  \bar\partial_t\varphi_t^1 & = \frac{t + 1}{1-|t|^2}\,i\varphi_t^1\wedge\overline\varphi_t^3 &      \\
\partial_t\varphi_t^2 & =  \frac{1 - \bar{t}}{1-|t|^2}\,\varphi_t^2\wedge\varphi_t^3, & \bar\partial_t\varphi_t^2 & = \frac{t - 1}{1-|t|^2}\,\varphi_t^2\wedge\overline\varphi_t^3  & \\
\partial_t\varphi_t^3 & =  0, &  \bar\partial_t\varphi_t^3 & = (t-1)\,i\varphi_t^1\wedge\overline\varphi_t^1 + (t + 1)\,\varphi_t^2\wedge\overline\varphi_t^2.  \end{align*}

For $t\in D$, we consider the Hermitian metric $\omega_t=\frac{1}{2}\,\sum_{j=1}^3 i\varphi_t^j\wedge\overline\varphi_t^j$ on $X_t$. Computing, we get \begin{eqnarray*}\bar\partial_t\omega_t = \frac{i}{2}(t-1)\,i\varphi_t^1\wedge\overline\varphi_t^1\wedge\overline\varphi_t^3 +  \frac{i}{2}(t+1)\,\varphi_t^2\wedge\overline\varphi_t^2\wedge\overline\varphi_t^3\end{eqnarray*} and \begin{eqnarray*}\partial_t\bar\partial_t\omega_t = \frac{1}{2}\,\bigg((t+1)(\bar{t}-1) - (t-1)(\bar{t}+1)\bigg)\,i\varphi_t^1\wedge\overline\varphi_t^1\wedge\varphi_t^2\wedge i\overline\varphi_t^2.\end{eqnarray*} Thus, \begin{eqnarray*}i\partial_t\bar\partial_t\omega_t = 2\,\mbox{Im}\,(t)\,i\varphi_t^1\wedge\overline\varphi_t^1\wedge\varphi_t^2\wedge i\overline\varphi_t^2,\end{eqnarray*} hence \begin{eqnarray*}\omega_t\wedge i\partial_t\bar\partial_t\omega_t = \mbox{Im}\,(t)\,i\varphi_t^1\wedge\overline\varphi_t^1\wedge\varphi_t^2\wedge i\overline\varphi_t^2\wedge \varphi_t^3\wedge i\overline\varphi_t^3 = 8\,\mbox{Im}\,(t)\,(\omega_t)_3.\end{eqnarray*} In particular, $f_{\omega_t} = \omega_t\wedge i\partial_t\bar\partial_t\omega_t/(\omega_t)_3 = 8\,\mbox{Im}\,(t)$, which is a real constant for every fixed $t$.

On the other hand, we have \begin{eqnarray*}(\omega_t)_2 = \frac{1}{4}\,\sum\limits_{j=1}^3\widehat{i\varphi_t^j\wedge\overline\varphi_t^j},\end{eqnarray*} where $\widehat{i\varphi_j\wedge\bar\varphi_j}$ stands for the product of all the $(1,\,1)$-forms $\widehat{i\varphi_k\wedge\bar\varphi_k}$ with $k\in\{1,\,2,\,3\}\setminus\{j\}$. Hence, \begin{eqnarray*}\star\rho_{\omega_t} & = & \frac{1}{2}\,f_{\omega_t}\,(\omega_t)_2 - i\partial_t\bar\partial_t\omega_t =  \mbox{Im}\,(t)\,\sum\limits_{j=1}^3\widehat{i\varphi_t^j\wedge\overline\varphi_t^j} - 2\,\mbox{Im}\,(t)\,\widehat{i\varphi_t^3\wedge\overline\varphi_t^3}  \\
 & = & \mbox{Im}\,(t)\,\bigg(\widehat{i\varphi_t^1\wedge\overline\varphi_t^1} + \widehat{i\varphi_t^2\wedge\overline\varphi_t^2} - \widehat{i\varphi_t^3\wedge\overline\varphi_t^3}\bigg).\end{eqnarray*} In particular, whenever $t$ is real, $\rho_{\omega_t}=0$, or equivalently $\partial_t\bar\partial_t\omega_t=0$, which means that the metric $\omega_t$ on $X_t$ is SKT.

To compute $\partial_t(\star\rho_{\omega_t})$, we first compute separately the following quantities: \begin{eqnarray*}\partial_t(\widehat{i\varphi_t^1\wedge\overline\varphi_t^1}) & = & (1-\bar{t})\,i\varphi_t^1\wedge\overline\varphi_t^1\wedge i\varphi_t^2\wedge\overline\varphi_t^2\wedge \varphi_t^3 \\
\partial_t(\widehat{i\varphi_t^2\wedge\overline\varphi_t^2}) & = & (1+\bar{t})\,i\varphi_t^1\wedge\overline\varphi_t^1\wedge i\varphi_t^2\wedge\overline\varphi_t^2\wedge \varphi_t^3 \\
\partial_t(\widehat{i\varphi_t^3\wedge\overline\varphi_t^3}) & = & 0.\end{eqnarray*} From this, we get: \begin{eqnarray*}\partial_t(\star\rho_{\omega_t}) = 2\,\mbox{Im}\,(t)\,i\varphi_t^1\wedge\overline\varphi_t^1\wedge i\varphi_t^2\wedge\overline\varphi_t^2\wedge \varphi_t^3.\end{eqnarray*} In particular, $\partial_t(\star\rho_{\omega_t})\neq 0$ whenever $t$ is not real. 

Finally, taking $\bar\partial_t$ in the expression for $\partial_t(\star\rho_{\omega_t})$, we get \begin{eqnarray*}\bar\partial_t\partial_t(\star\rho_{\omega_t}) & = & 2\,\mbox{Im}\,(t)\,i\varphi_t^1\wedge\overline\varphi_t^1\wedge i\varphi_t^2\wedge\overline\varphi_t^2\wedge\bar\partial_t\varphi_t^3 \\
& = & 2\,\mbox{Im}\,(t)\,i\varphi_t^1\wedge\overline\varphi_t^1\wedge i\varphi_t^2\wedge\overline\varphi_t^2\wedge\bigg((t-1)\,i\varphi_t^1\wedge\overline\varphi_t^1 + (t + 1)\,\varphi_t^2\wedge\overline\varphi_t^2\bigg) = 0\end{eqnarray*} for every $t\in D$. This means that the metric $\omega_t$ is pluriclosed star split, hence also Gauduchon since $f_{\omega_t}$ is constant on $X_t$ (see (i) of Proposition \ref{Prop:pluriclosed-star-split_properties}).

The conclusions of this computation are summed up in the following

\begin{Prop}\label{Prop::small-def-S3S3_examp} For every $t\in D$, the Hermitian metric $\omega_t$ on the corresponding {\bf small deformation} $X_t$ of the {\bf Calabi-Eckmann manifold} $X=X_0 = (S^3\times S^3,\, J_{CE})$ is {\bf pluriclosed star split} and its associated function $f_{\omega_t}$ is {\bf constant} equal to $8\,\mbox{Im}\,(t)$.

The $\partial\bar\partial$-closed real $(2,\,2)$-form $\star\rho_{\omega_t}$ is $d$-closed if and only if $t\in\R$, in which case $\rho_{\omega_t} = 0$ and the metric $\omega_t$ on $X_t$ is SKT. The eigenvalues of $\star\rho_{\omega_t}$ with respect to $\omega$ are $4\mbox{Im}\,(t), 4\mbox{Im}\,(t), -4\mbox{Im}\,(t)$.

\end{Prop}

\subsubsection{Example of a Hermitian metric $\omega$ whose function $f_\omega$ is non-constant}\label{subsubsection:examp_non-constant_f-omega}

Since $f_\omega$ is constant in all of the above examples, we will now try to impress upon the reader the relative ease with which one can produce examples of metrics whose associated function is non-constant by varying a metric $\omega$ with constant function $f_\omega$ in its conformal class. If the original $\omega$ is a Gauduchon metric and $g$ is a non-constant smooth positive function on the manifold, the new metric $g\omega$ is no longer Gauduchon due to the uniqueness, up to multiplicative positive constants, of a Gauduchon metric in its conformal class proved in Gauduchon's foundational work [Gau77a].  

We start by noticing a general formula describing the variation of $f_\omega$ under conformal changes of a balanced $\omega$. For the sake of convenience, we will confine ourselves to the $3$-dimensional case.

\begin{Lem}\label{Lem:f_omega-changes_conformal-var} Let $\omega$ be a {\bf balanced} Hermitian metric on a $3$-dimensional complex manifold $X$.

  Then, for any $C^\infty$ function $g:X\longrightarrow (0,\,+\infty)$, one has: \begin{eqnarray}\label{eqn:f_omega-changes_conformal-var}f_{g\omega} = \frac{1}{g}\,f_\omega - \frac{2}{g^2}\,\Delta_\omega(g),\end{eqnarray} where $\Delta_\omega$ is the standard Laplacian defined on the functions $h$ on $X$ by $\Delta_\omega(h)=-\Lambda_\omega(i\partial\bar\partial h)$.

\end{Lem}  

\noindent {\it Proof.} In the case $n=3$, formula (\ref{eqn:function_f_omega}) spells: \begin{eqnarray*}\label{eqn:function_f_omega_proof-conformal_1}f_{g\omega}=\frac{\omega\wedge i\partial\bar\partial(g\omega)}{g^2\omega_3}.\end{eqnarray*}

On the other hand, we get: \begin{eqnarray*}i\partial\bar\partial(g\omega) = gi\partial\bar\partial\omega + i\partial\bar\partial g\wedge\omega + i\partial g\wedge\bar\partial\omega - i\bar\partial g\wedge\partial\omega.\end{eqnarray*}

Hence, we get: \begin{eqnarray*}f_{g\omega} & = & \frac{1}{g}\,f_\omega + \frac{2}{g^2}\,\frac{i\partial\bar\partial g\wedge\omega_2}{\omega_3} + \frac{1}{g^2}\,\frac{i\partial g\wedge\omega\wedge\bar\partial\omega}{\omega_3} - \frac{1}{g^2}\,\frac{i\bar\partial g\wedge\omega\wedge\partial\omega}{\omega_3}\\
  & = & \frac{1}{g}\,f_\omega + \frac{2}{g^2}\,\Lambda_\omega(i\partial\bar\partial g),\end{eqnarray*} where the balanced hypothesis on $\omega$ was used to infer the latter equality in each of the pairs of equalities $\omega\wedge\bar\partial\omega = (1/2)\,\bar\partial\omega^2 = 0$ and $\omega\wedge\partial\omega = (1/2)\,\partial\omega^2 = 0$.  \hfill $\Box$

\vspace{2ex}

We are now in a position to describe our example. We choose $X$ to be the {\it Iwasawa manifold} equipped with the {\it balanced} metric $\omega$ described in $\S$\ref{subsubsection:Iwasawa_examp}. We saw in Proposition \ref{Prop:Iwasawa_examp} that $f_\omega$ is the constant function $1$ on $X$.

We will multiply $\omega$ by the non-constant positive $C^\infty$ function $g:X\longrightarrow(0,\,+\infty)$ induced by the function (denoted by the same symbol): \begin{eqnarray}\label{eqn:g_def_non-constant-function}g:\C^3\longrightarrow(0,\,+\infty), \hspace{3ex} g(z_1,\,z_2,\,z_3)=e^{\sin(2\pi\,\mbox{Re}\,z_1)}.\end{eqnarray} Recall that $\C^3$ coincides, as a complex manifold, with the Heisenberg group $G$. On the other hand, the action on $G$ that defines the Iwasawa manifold $X$ acts on the first coordinate as the sum, namely $z_1\in\C$ is mapped to $z_1 + (p_1 + iq_1)$ with $p_1$ and $q_1$ ranging over $\Z$. This implies that the above function $g$ defined on $\C^3$ is {\it constant} along the orbits of the action. Hence, $g$ passes to the quotient and defines a function on $X$.\footnote{The author is grateful to L. Ugarte for pointing out this function to him.}

In order to apply formula (\ref{eqn:f_omega-changes_conformal-var}) to $g\omega$, we compute $\Delta_\omega(g)$. We get successively: \begin{eqnarray*}\bar\partial g & = & \pi g\,\cos(2\pi\,\mbox{Re}\,z_1)\,d\bar{z}_1 \\
  i\partial\bar\partial g & = & \pi^2 g\,\bigg(\cos(2\pi\,\mbox{Re}\,z_1) - \sin(2\pi\,\mbox{Re}\,z_1)\bigg)\,idz_1\wedge d\bar{z}_1  \\
  \Delta_\omega(g) = -\Lambda_\omega(i\partial\bar\partial g) & = & - \pi^2 g\,\bigg(\cos(2\pi\,\mbox{Re}\,z_1) - \sin(2\pi\,\mbox{Re}\,z_1)\bigg).\end{eqnarray*}

Thus, formula (\ref{eqn:f_omega-changes_conformal-var}) and $f_\omega \equiv 1$ yield: \begin{eqnarray*}\label{eqn:f_omega-changes_conformal-var_example}f_{g\omega} & = & \frac{1}{g} + \frac{2\pi^2}{g}\,\bigg(\cos(2\pi\,\mbox{Re}\,z_1) - \sin(2\pi\,\mbox{Re}\,z_1)\bigg)\\
  & = & \frac{1 + 2\pi^2\,\bigg(\cos(2\pi\,\mbox{Re}\,z_1) - \sin(2\pi\,\mbox{Re}\,z_1)\bigg)}{e^{\sin(2\pi\,\mbox{Re}\,z_1)}}.\end{eqnarray*} This function is readily seen to be non-constant. For example, the above formula shows that:

\vspace{1ex}

(i)\, if $\mbox{Re}\,z_1=0$, then $f_{g\omega}(z_1,\,z_2,\,z_3) = 1 + 2\pi^2$;

\vspace{1ex}

(ii)\, if $\mbox{Re}\,z_1=\frac{1}{4}$, then $f_{g\omega}(z_1,\,z_2,\,z_3) = \frac{1 - 2\pi^2}{e}$,

\vspace{1ex}

\noindent so $f_{g\omega}$ assumes different values at the above two sets of points.

\vspace{2ex}

We have thus proved the following

\begin{Prop}\label{Prop:non-constant_omega-function} Let $(X,\,\omega)$ be the Iwasawa manifold equipped with the balanced metric described in $\S$\ref{subsubsection:Iwasawa_examp}. Let $g$ be the $C^\infty$ function $g:X\longrightarrow(0,\,+\infty)$ induced by the function defined in (\ref{eqn:g_def_non-constant-function}).

  Then, the function $f_{g\omega}$ associated with the Hermitian metric $g\omega$ on $X$ is {\bf non-constant}. 

\end{Prop}

\section{The pluriclosed star split condition for pairs of metrics}\label{section:def_two-metrics} Let $X$ be an $n$-dimensional complex manifold and let $\omega$ and $\gamma$ be arbitrary Hermitian metrics on $X$. Since the pointwise map \begin{eqnarray}\label{eqn:gamma_n-2_map}\gamma_{n-2}\wedge\cdot :\Lambda^{1,\,1}T^\star X\longrightarrow\Lambda^{n-1,\,n-1}T^\star X\end{eqnarray} is bijective, there exists a unique $C^\infty$ $(1,\,1)$-form $\rho_{\omega,\,\gamma}$ on $X$ such that \begin{eqnarray}\label{eqn:rho_omega-gamma_def}i\partial\bar\partial\omega_{n-2} = \gamma_{n-2}\wedge\rho_{\omega,\,\gamma}.\end{eqnarray}

In particular, $\rho_{\omega,\,\gamma} = 0$ if and only if $\omega$ is astheno-K\"ahler. Together with Proposition \ref{Prop:astheno-balanced-K}, this implies the following

\begin{Obs}\label{Obs:balanced_rho-zero_K} Let $\omega$ be a Hermitian metric on a compact complex manifold $X$ with $\mbox{dim}_\C X =n\geq 3$. If $\omega$ is {\bf balanced}, then $\omega$ is {\bf K\"ahler} if and only if $\rho_{\omega,\,\gamma}=0$ for some (hence every) Hermitian metric $\gamma$.  

\end{Obs}

The form $\rho_{\omega,\,\gamma}$ can be computed in the same way as $\rho_\omega$ was computed in $\S.$\ref{subsection:gen-theory_one-metric}. Without repeating all the steps of that computation, we only mention that the Lefschetz decomposition with respect to $\gamma$ spells: \begin{eqnarray*}\rho_{\omega,\,\gamma} = \rho_{\omega,\,\gamma,\,prim} + \frac{1}{n}\,\Lambda_\gamma(\rho_{\omega,\,\gamma})\,\gamma,\end{eqnarray*} where $\rho_{\omega,\,\gamma,\,prim}$ is a $(1,\,1)$-form that is primitive with respect to $\gamma$, namely $\Lambda_\gamma(\rho_{\omega,\,\gamma,\,prim})=0$ or, equivalently, $\rho_{\omega,\,\gamma}\wedge\gamma_{n-1}=0$. Then, using the Hodge star operator $\star_\gamma$ induced by $\gamma$ as we did for $\omega$ in $\S.$\ref{subsection:gen-theory_one-metric}, we get the following analogue of Lemma \ref{Lem:star-rho_omega_formula}.

\begin{Lem}\label{Lem:star-rho_omega-gamma_formula} The $(1,\,1)$-form $\rho_{\omega,\,\gamma}$ uniquely determined by an arbitrary pair $(\omega,\,\gamma)$ of Hermitian metrics on an $n$-dimensional complex manifold $X$ via property (\ref{eqn:rho_omega-gamma_def}) is given by the formula \begin{eqnarray}\label{eqn:star-rho_omega-gamma_formula}\star_\gamma\rho_{\omega,\,\gamma} = \frac{1}{n-1}\,\frac{\gamma\wedge i\partial\bar\partial\omega_{n-2}}{\gamma_n}\,\gamma_{n-1} - i\partial\bar\partial\omega_{n-2}.\end{eqnarray}

\end{Lem}

In particular, we see that with each pair $(\omega,\,\gamma)$ of Hermitian metrics on $X$ we can associate the $C^\infty$ function $f_{\omega,\,\gamma}:X\to\R$ defined by \begin{eqnarray}\label{eqn:function_f_omega-gamma}f_{\omega,\,\gamma}:=\frac{\gamma\wedge i\partial\bar\partial\omega_{n-2}}{\gamma_n}.\end{eqnarray} As with $f_\omega$ in $\S.$\ref{subsection:gen-theory_one-metric}, we also get $f_{\omega,\,\gamma} = (n-1)\,\Lambda_\gamma(\rho_{\omega,\,\gamma})$. This immediately yields the following

\begin{Lem}\label{Lem:f-vanishing_equiv} Let $\omega$ and $\gamma$ be Hermitian metrics on a compact connected complex manifold $X$ with $\mbox{dim}_\C X =n\geq 3$.

The function $f_{\omega,\,\gamma}$ vanishes identically if and only if $\Lambda_\gamma(\rho_{\omega,\,\gamma}) = 0$.

\end{Lem}

Note that this statement also follows from the following considerations.

\vspace{2ex}

On the one hand, we always have the following equivalences: \begin{eqnarray}\label{eqn:equiv_f-omega-gamma_id-zero}f_{\omega,\,\gamma}\equiv 0 \iff \gamma\wedge i\partial\bar\partial\omega_{n-2} = 0 \iff \Lambda_\gamma\bigg(\star_\gamma(i\partial\bar\partial\omega_{n-2})\bigg) = 0,\end{eqnarray} where the former follows from (\ref{eqn:function_f_omega-gamma}) and the latter from Lemma \ref{Lem:star-trace-11n-1n-1}. 

On the other hand, (iii) of (\ref{eqn:standard_comm_Lambda-L-powers}) implies the first equality below: \begin{eqnarray}\label{eqn:equiv_f-omega-gamma_id-zero_bis}\Lambda_\gamma\bigg(\star_\gamma(i\partial\bar\partial\omega_{n-2})\bigg) = \star_\gamma(\gamma\wedge i\partial\bar\partial\omega_{n-2}) = (n-1)\,\star_\gamma(\gamma_{n-1}\wedge\rho_{\omega,\,\gamma}) = (n-1)\,\Lambda_\gamma(\rho_{\omega,\,\gamma}),\end{eqnarray} where the second equality follows from (\ref{eqn:rho_omega-gamma_def}) and the third from $\gamma_{n-1}\wedge\rho_{\omega,\,\gamma} = \Lambda_\gamma(\rho_{\omega,\,\gamma})\,\gamma_n$ and from $\star_\gamma\gamma_n = 1$.

Thus, putting (\ref{eqn:equiv_f-omega-gamma_id-zero}) and (\ref{eqn:equiv_f-omega-gamma_id-zero_bis}) together, we get again Lemma \ref{Lem:f-vanishing_equiv}.

\begin{Def}\label{Def:pluriclosed-star-split_pairs} Let $(\omega,\,\gamma)$ be a pair of Hermitian metrics on an $n$-dimensional complex manifold $X$, let $\star_\gamma$ be the Hodge star operator induced by $\gamma$ and let $\rho_{\omega,\,\gamma}$ be the $(1,\,1)$-form on $X$ uniquely determined by $\omega$ and $\gamma$ via property (\ref{eqn:rho_omega-gamma_def}).

\vspace{1ex}

(a)\, The pair $(\omega,\,\gamma)$ is said to be {\bf pluriclosed star split} if $\partial\bar\partial(\star_\gamma\rho_{\omega,\,\gamma}) = 0$. 

\vspace{1ex}

(b)\, The pair $(\omega,\,\gamma)$ is said to be {\bf closed star split} if $d(\star_\gamma\rho_{\omega,\,\gamma}) = 0$.

\end{Def}

Thanks to formula (\ref{eqn:star-rho_omega-gamma_formula}), we see that the {\bf pluriclosed star split} condition on the pair $(\omega,\,\gamma)$ is equivalent to \begin{eqnarray}\label{eqn:pluriclosed-star-split_pairs_bis}\partial\bar\partial(f_{\omega,\,\gamma}\,\gamma_{n-1}) = 0,\end{eqnarray} while the {\bf closed star split} condition on $(\omega,\,\gamma)$ is equivalent to \begin{eqnarray}\label{eqn:closed-star-split_pairs_bis}d(f_{\omega,\,\gamma}\,\gamma_{n-1}) = 0.\end{eqnarray}

In particular, if the function $f_{\omega,\,\gamma}$ vanishes identically, the pair $(\omega,\,\gamma)$ is closed star split, hence also pluriclosed star split. See Corollary \ref{Cor:SKT_pss_f-omega-gamma_zero} for a partial converse. 

\vspace{2ex}

In the present two-metric context, the analogue of Proposition \ref{Prop:pluriclosed-star-split_3-choices_f} is the following

\begin{Prop}\label{Prop:pluriclosed-star-split-pair_3-choices_f} Let $X$ be a compact connected complex manifold with $\mbox{dim}_\C X = n$ and let $\omega$ and $\gamma$ be Hermitian metrics on $X$. The following equivalence holds: \begin{eqnarray*}(\omega,\,\gamma) \hspace{1ex}\mbox{is {\bf pluriclosed star split}} \iff \Delta_\gamma^\star(f_{\omega,\,\gamma}) =0,\end{eqnarray*} where $\Delta_\gamma^\star$ is the adjoint, given by formula (\ref{eqn:Laplace-adjoint_functions}), of the Laplace-type operator on functions: $\Delta_\gamma:=i\Lambda_\gamma\bar\partial\partial:C^\infty(X,\,\C)\longrightarrow C^\infty(X,\,\C)$.

Moreover, any of these equivalent conditions implies that \begin{eqnarray*}f_{\omega,\,\gamma}> 0 \hspace{1ex} \mbox{on} \hspace{1ex} X \hspace{5ex} \mbox{or} \hspace{5ex} f_{\omega,\,\gamma}< 0 \hspace{1ex} \mbox{on} \hspace{1ex} X  \hspace{5ex} \mbox{or} \hspace{5ex} f_{\omega,\,\gamma}\equiv 0.\end{eqnarray*} 

\end{Prop}

\vspace{3ex}

The following statement is the analogue of Proposition \ref{Prop:pluriclosed-star-split_properties}.

\begin{Prop}\label{Prop:pluriclosed-star-split-pair_properties} Let $\omega$ and $\gamma$ be Hermitian metrics on a complex manifold $X$ with $\mbox{dim}_\C X=n$. 

 \vspace{1ex}

(i)\, If the function $f_{\omega,\,\gamma}$ is a {\bf non-zero constant}, the pair $(\omega,\,\gamma)$ is {\bf pluriclosed star split} if and only if the metric $\gamma$ is {\bf Gauduchon}.

\vspace{1ex}

(ii)\, Suppose $X$ is {\bf compact} and connected. If the metric $\gamma$ is {\bf balanced}, the pair $(\omega,\,\gamma)$ is {\bf pluriclosed star split} if and only if the function $f_{\omega,\,\gamma}$ is {\bf constant}.

\vspace{1ex}

(iii)\, Suppose $X$ is {\bf compact} and connected. If the metric $\gamma$ is {\bf Gauduchon}, the pair $(\omega,\,\gamma)$ is {\bf pluriclosed star split} if and only if the function $f_{\omega,\,\gamma}$ is {\bf constant}.

\end{Prop}

\noindent {\it Proof.} It is the analogue of the proof of Proposition \ref{Prop:pluriclosed-star-split_properties} with $\omega$ replaced by $\gamma$ at the obvious places. For example, when $\gamma$ is balanced, the pair $(\omega,\,\gamma)$ is pluriclosed star split if and only if $\Delta_\gamma(f_{\omega,\,\gamma})=0$ on $X$, a condition that is equivalent to $f_{\omega,\,\gamma}$ being constant. \hfill $\Box$

\vspace{3ex}

The following statement is the analogue of Proposition \ref{Prop:balanced_f-omega_del-omega}.

\begin{Prop}\label{Prop:balanced-SKT_f-omega-gamma_del} Let $\omega$ and $\gamma$ be Hermitian metrics on a compact complex manifold $X$ with $\mbox{dim}_\C X =n\geq 3$.

\vspace{1ex}

(i)\, If $\gamma$ is {\bf SKT}, then $\int_Xf_{\omega,\,\gamma}\,\gamma_n = 0$.

 \vspace{1ex}

 (ii)\, If $\omega$ is {\bf balanced}, then \begin{eqnarray}\label{eqn:balanced_f-omega_del-omega}\int_Xf_{\omega,\,\gamma}\,\gamma_n = \langle\langle\partial\gamma,\,\partial\omega\rangle\rangle_\omega.\end{eqnarray}

\end{Prop}

\noindent {\it Proof.} From (\ref{eqn:star-rho_omega-gamma_formula}), we get \begin{eqnarray*}\gamma\wedge\star_\gamma\rho_{\omega,\,\gamma} = \frac{1}{n-1}\,\gamma\wedge i\partial\bar\partial\omega_{n-2} =  \frac{1}{n-1}\,f_{\omega,\,\gamma}\,\gamma_n.\end{eqnarray*} Hence, integrating and applying Stokes twice, we get: \begin{eqnarray*}\int_Xf_{\omega,\,\gamma}\,\gamma_n = i\,\int_X\omega_{n-2}\wedge\partial\bar\partial\gamma.\end{eqnarray*} Since $\partial\bar\partial\gamma = 0$ when $\gamma$ is SKT, this proves part (i).

On the other hand, integrating as above but applying Stokes only once, we get: \begin{eqnarray*}\int_Xf_{\omega,\,\gamma}\,\gamma_n = \int_X i\,\bar\partial\omega_{n-2}\wedge\partial\gamma = - i\,\int_X\partial\gamma\wedge\omega_{n-3}\wedge\bar\partial\omega.\end{eqnarray*} As we saw in the proof of Proposition \ref{Prop:balanced_f-omega_del-omega}, if $\omega$ is {\it balanced}, $\bar\partial\omega$ is primitive, so $\star(\bar\partial\omega) = -i\,\omega_{n-3}\wedge\bar\partial\omega$. Thus, the last integral formula becomes: \begin{eqnarray*}\int_Xf_{\omega,\,\gamma}\,\gamma_n = \int_X\partial\gamma\wedge\star(\bar\partial\omega) = \langle\langle\partial\gamma,\,\partial\omega\rangle\rangle_\omega,\end{eqnarray*} as claimed under (ii).  \hfill $\Box$

\begin{Cor}\label{Cor:SKT_pss_f-omega-gamma_zero} Let $\omega$ and $\gamma$ be Hermitian metrics on a compact connected complex manifold $X$ with $\mbox{dim}_\C X =n\geq 3$.

  If $\gamma$ is {\bf SKT} and the pair $(\omega,\,\gamma)$ is {\bf pluriclosed star split}, the function $f_{\omega,\,\gamma}$ vanishes identically.

\end{Cor}

\noindent {\it Proof.} Since $\gamma$ is SKT, $\int_Xf_{\omega,\,\gamma}\,\gamma_n = 0$ by (i) of Proposition \ref{Prop:balanced-SKT_f-omega-gamma_del}. This implies that $f_{\omega,\,\gamma}\equiv 0$ thanks to the pair $(\omega,\,\gamma)$ being pluriclosed star split and to Proposition \ref{Prop:pluriclosed-star-split-pair_3-choices_f}.  \hfill $\Box$

\vspace{2ex}

Putting Corollary \ref{Cor:SKT_pss_f-omega-gamma_zero} and Lemma \ref{Lem:f-vanishing_equiv} together, we can now prove another result stated in the introduction.

\vspace{2ex}

\noindent {\bf Proof of Theorem \ref{The:SKT_pss_rho-pos_K}.} By Corollary \ref{Cor:SKT_pss_f-omega-gamma_zero}, $f_{\omega,\,\gamma}$ vanishes identically. By Lemma \ref{Lem:f-vanishing_equiv}, this amounts to $\Lambda_\gamma(\rho_{\omega,\,\gamma}) = 0$. Thus, the trace of $\rho_{\omega,\,\gamma}$ w.r.t. $\gamma$ vanishes. Meanwhile, the hypothesis also ensures that the eigenvalues of $\rho_{\omega,\,\gamma}$ w.r.t. $\gamma$ are either all non-negative or all non-positive. Therefore, all the eigenvalues must vanish, which means that $\rho_{\omega,\,\gamma} = 0$. Thanks to (\ref{eqn:rho_omega-gamma_def}), this implies $i\partial\bar\partial\omega_{n-2} = 0$, which means that $\omega$ is astheno-K\"ahler.

The last statement follows from Proposition \ref{Prop:astheno-balanced-K}.  \hfill $\Box$

\section{Maps and the pluriclosed star split condition}\label{section:maps_pluriclosed} In this section, we give the first applications of the constructions performed in the previous two sections adapted to the context of holomorphic maps between complex Hermitian manifolds.

Let $X$, $Y$ be complex manifolds of respective dimensions $n$ and $m$ with $n\leq m$, let $\gamma$, $\omega$ be Hermitian metrics on $X$, respectively $Y$ and let \begin{eqnarray*}\phi:(X,\,\gamma)\longrightarrow(Y,\,\omega)\end{eqnarray*} be a holomorphic map supposed to be non-degenerate at some point $x\in X$. By this we mean that the differential map $d_x\phi:T^{1,\,0}_xX\longrightarrow T^{1,\,0}_{\phi(x)}Y$ at $x$ is of maximal rank (i.e. rank $n$). Since the points where $\phi$ degenerates form a (possibly empty) analytic subset $\Sigma\subset X$, the map $\phi$ is non-degenerate at least almost everywhere on $X$, namely everywhere on $X\setminus\Sigma$.

The pullback form $$\widetilde\omega:=\phi^\star\omega\geq 0$$ is a $C^\infty$ positive semidefinite $(1,\,1)$-form on $X$ that is positive definite on $X\setminus\Sigma$. Thus, it is a degenerate metric on $X$ and a genuine Hermitian metric on $X\setminus\Sigma$. We can rerun the construction in $\S.$\ref{section:def_two-metrics} with $\widetilde\omega$ and $\gamma$ in place of $\omega$ and $\gamma$.  

Thus, there exists exists a unique $C^\infty$ real $(1,\,1)$-form $\rho_{\phi,\,\omega,\,\gamma}$ on $X$ such that \begin{eqnarray}\label{eqn:rho_phi-omega-gamma_def}i\partial\bar\partial\widetilde\omega_{n-2} = \gamma_{n-2}\wedge\rho_{\phi,\,\omega,\,\gamma}.\end{eqnarray} Intuitively, $\rho_{\phi,\,\omega,\,\gamma}$ can be thought of as a kind of {\it curvature} form for the triple $(\phi,\,\omega,\,\gamma)$. As in $\S.$\ref{section:def_one-metric} and $\S.$\ref{section:def_two-metrics}, (\ref{eqn:rho_phi-omega-gamma_def}) implies that \begin{eqnarray}\label{eqn:star-rho_phi-omega-gamma_formula}\star_\gamma\rho_{\phi,\,\omega,\,\gamma} = \frac{1}{n-1}\,f_{\phi,\,\omega,\,\gamma}\,\gamma_{n-1} - i\partial\bar\partial\widetilde\omega_{n-2},\end{eqnarray} where $f_{\phi,\,\omega,\,\gamma}:X\to\R$ is the $C^\infty$ function \begin{eqnarray}\label{eqn:function_f_phi-omega-gamma}f_{\phi,\,\omega,\,\gamma}:=\frac{\gamma\wedge i\partial\bar\partial\widetilde\omega_{n-2}}{\gamma_n} = (n-1)\,\Lambda_\gamma(\rho_{\phi,\,\omega,\,\gamma}).\end{eqnarray}

\begin{Def}\label{Def:pluriclosed-star-split_triples} A triple $(\phi,\,\omega,\,\gamma)$ of a holomorphic map $\phi:(X,\,\gamma)\longrightarrow(Y,\,\omega)$ that is non-degenerate at some point $x\in X$ and Hermitian metrics $\omega$ and $\gamma$ on the complex manifolds $Y$, respectively $X$, is said to be {\bf pluriclosed star split} if $\partial\bar\partial(\star_\gamma\rho_{\phi,\,\omega,\,\gamma}) = 0$.

  The triple $(\phi,\,\omega,\,\gamma)$ is said to be {\bf closed star split} if $d(\star_\gamma\rho_{\phi,\,\omega,\,\gamma}) = 0$.

\end{Def}

\vspace{2ex}

As in $\S.$\ref{section:def_two-metrics}, especially in Proposition \ref{Prop:pluriclosed-star-split-pair_3-choices_f}, we have the following equivalences when $X$ is compact: \begin{eqnarray}\label{eqn:pluriclosed-star-split-triple_f}(\phi,\,\omega,\,\gamma) \hspace{1ex} \mbox{is {\bf pluriclosed star split}} \iff \partial\bar\partial(f_{\phi,\,\omega,\,\gamma}\,\gamma_{n-1}) = 0 \iff \Delta_\gamma^\star(f_{\phi,\,\omega,\,\gamma}) = 0.\end{eqnarray} Thus, the pluriclosed star split condition on the triple $(\phi,\,\omega,\,\gamma)$ is a kind of harmonicity condition on the function $f_{\phi,\,\omega,\,\gamma}$, hence also on the map $\phi$.

When $X$ is compact, any of the equivalent conditions (\ref{eqn:pluriclosed-star-split-triple_f}) implies that \begin{eqnarray*}f_{\phi,\,\omega,\,\gamma}> 0 \hspace{1ex} \mbox{on} \hspace{1ex} X \hspace{5ex} \mbox{or} \hspace{5ex} f_{\phi,\,\omega,\,\gamma}< 0 \hspace{1ex} \mbox{on} \hspace{1ex} X  \hspace{5ex} \mbox{or} \hspace{5ex} f_{\phi,\,\omega,\,\gamma}\equiv 0.\end{eqnarray*}

\vspace{2ex}

Combining the arguments in Propositions \ref{Prop:astheno-balanced-K}, \ref{Prop:balanced-SKT_f-omega-gamma_del} and in Theorem \ref{The:SKT_pss_rho-pos_K}, we get the following result where the map $\phi$ is supposed to be non-degenerate at every point of $X$, so $\widetilde\omega$ defines a Hermitian metric on the whole of $X$.

\begin{The}\label{The:map_two-metrics_positivity} Let $\phi:X\longrightarrow Y$ be a surjective holomorphic submersion between compact complex manifolds of the same dimension $n\geq 3$. Suppose there exist an {\bf SKT} metric $\gamma$ on $X$ and a {\bf balanced} metric $\omega$ on $Y$ such that the triple $(\phi,\,\omega,\,\gamma)$ is {\bf pluriclosed star split}.

  If the $(1,\,1)$-form $\rho_{\phi,\,\omega,\,\gamma}$ is either {\bf positive semi-definite} or {\bf negative semi-definite} on $X$, the metric $\omega$ of $Y$ is {\bf K\"ahler}.

\end{The}

It may be useful to compare this result with Siu's rigidity results in [Siu80, Theorems 1 and 5] to which it is, in a certain sense, complementary. In Siu's case, the map $\phi$ (denoted by $f$ there) is supposed to be harmonic (a role played by our pluriclosed star split hypothesis) and to satisfy a certain non-degeneracy assumption, while the curvature tensor of $Y$ (denoted by $M$ there) is supposed to be strongly negative or negative of a certain order (a role played by the semi-definiteness assumption on $\rho_{\phi,\,\omega,\,\gamma}$ in our case). Both manifolds are supposed to be K\"ahler in [Siu80] and the conclusion is that the map $\phi$ is holomorphic or conjugate holomorphic. In our case, the map $\phi$ is supposed to be holomorphic, but neither of the manifolds $X$ and $Y$ is supposed to be K\"ahler. We obtain the K\"ahlerianity of $Y$ as the conclusion of our result.

\vspace{2ex}

\noindent {\it Proof of Theorem \ref{The:map_two-metrics_positivity}.} Since $\phi$ is non-degenerate everywhere, $\widetilde\omega$ is a Hermitian metric on $X$, so the pluriclosed star split hypothesis on the triple $(\phi,\,\omega,\,\gamma)$ is equivalent to the pair $(\widetilde\omega,\,\gamma)$ being pluriclosed star split. In particular, since $\gamma$ is SKT, Corollary \ref{Cor:SKT_pss_f-omega-gamma_zero} ensures that $f_{\phi,\,\omega,\,\gamma} = f_{\widetilde\omega,\,\gamma}$ vanishes identically on $X$. As in the proof of Theorem \ref{The:SKT_pss_rho-pos_K}, this implies, together with the semi-definiteness assumption on $\rho_{\phi,\,\omega,\,\gamma} = \rho_{\widetilde\omega,\,\gamma}$, that $\rho_{\widetilde\omega,\,\gamma} = 0$, a fact that amounts to the metric $\widetilde\omega$ of $X$ being astheno-K\"ahler.

On the other hand, the metric $\omega$ of $Y$ is balanced, hence so is the metric $\widetilde\omega$ of $X$ because $d\widetilde\omega^{n-1} = d\phi^\star(\omega^{n-1}) = \phi^\star(d\omega^{n-1}) = 0$.

Now, Proposition \ref{Prop:astheno-balanced-K} ensures that the metric $\widetilde\omega$ of $X$ is K\"ahler. This means that $\phi^\star(d\omega) = 0$ on $X$, hence $s\,d\omega = \phi_\star\phi^\star(d\omega) = 0$ on $Y$, where $s\in\N^\star$ is the number of elements in the fibre $\phi^{-1}(y)$ for any $y\in Y$. Thus, $d\omega=0$, so $\omega$ is a K\"ahler metric on $Y$. \hfill $\Box$

\vspace{3ex}

By taking $X=Y$ in Theorem \ref{The:map_two-metrics_positivity}, we get Theorem \ref{The:SKT-bal_map_pss_pos} stated in the introduction.

\vspace{3ex}

To further stress the possible relevance of the automorphism group $Aut(X)$ to the SKT-balanced conjecture and the relations among the various notions of {\it pluriclosed star split} objects introduced above, we now digress briefly, starting from the following very simple (and probably known) observation dealing with {\it $\gamma$-isometries} (i.e. automorphisms $\phi$ of $X$ that preserve a given Hermitian metric $\gamma$ on $X$ in the sense that $\phi^\star\gamma = \gamma$).

\begin{Lem}\label{Lem:star_Phi-star_commutation_isometry} Let $(X,\,\gamma)$ be a Hermitian complex manifold and let $\phi:X\longrightarrow X$ be a biholomorphism such that $\phi^\star\gamma = \gamma$. Then \begin{eqnarray}\label{eqn:star_Phi-star_commutation_isometry}\star_\gamma\circ\phi^\star = \phi^\star\circ\star_\gamma,\end{eqnarray} where $\star_\gamma$ is the Hodge star operator defined by the metric $\gamma$ and $\phi^\star : C^\infty_{p,\,q}(X,\,\C)\longrightarrow C^\infty_{p,\,q}(X,\,\C)$ is the pullback map under $\phi$ for smooth differential forms of any bidegree $(p,\,q)$ on $X$.

\end{Lem}

\noindent {\it Proof.} Recall that $\star_\gamma$ is defined by the pointwise identity $u\wedge\star_\gamma\bar{v} = \langle u,\,v\rangle_\gamma\,dV_\gamma$ required to hold for all $(p,\,q)$-forms $u,v$ on $X$, where $\langle\cdot,\,\cdot\rangle_\gamma$ is the pointwise inner product and $dV_\gamma=\gamma^n/n!$ is the volume form defined by $\gamma$, $n$ being the complex dimension of $X$. Thus, we get: \begin{eqnarray*}\phi^\star u\wedge\star_\gamma\overline{\phi^\star v} & = & \phi^\star u\wedge\star_{\phi^\star\gamma}\overline{\phi^\star v} = \langle \phi^\star u,\,\phi^\star v\rangle_{\phi^\star\gamma}\,dV_{\phi^\star\gamma} = \bigg(\langle u,\,v\rangle_\gamma\circ\phi\bigg)\,\phi^\star(dV_\gamma) \\
  & = & \phi^\star\bigg(\langle u,\,v\rangle_\gamma\,dV_\gamma\bigg) = \phi^\star\bigg(u\wedge\star_\gamma\bar{v}\bigg) = \phi^\star u\wedge\overline{\phi^\star(\star_\gamma v)},\end{eqnarray*} where for the first equality we used the hypothesis $\phi^\star\gamma = \gamma$. This proves that $\star_\gamma(\phi^\star v) = \phi^\star(\star_\gamma v)$ for every form $v$, as claimed.  \hfill $\Box$

\vspace{3ex}

This leads to the following

\begin{Prop}\label{Prop:pluriclosed-star-split_various-relations} Let $X$ be a complex manifold equipped with Hermitian metrics $\omega$ and $\gamma$.

\vspace{1ex}  

(a) Let $\phi:X\longrightarrow X$ be a biholomorphism such that $\phi^\star\gamma = \gamma$. If the pair $(\omega,\,\gamma)$ is {\bf pluriclosed star split}, the triple $(\phi,\,\omega,\,\gamma)$ is {\bf pluriclosed star split} and we have $\rho_{\phi,\,\omega,\,\gamma}:=\phi^\star\rho_{\omega,\,\gamma}$.

\vspace{1ex}

(b) Let $\phi, \psi:X\longrightarrow X$ be biholomorphisms such that $\phi^\star\gamma = \gamma$ and $\psi^\star\gamma = \gamma$. If the triples $(\phi,\,\omega,\,\gamma)$ and $(\psi,\,\omega,\,\gamma)$ are {\bf pluriclosed star split}, the triples $(\phi\circ\psi,\,\omega,\,\gamma)$ and $(\psi\circ\phi,\,\omega,\,\gamma)$ are again {\bf pluriclosed star split} and we have $\rho_{\phi\circ\psi,\,\omega,\,\gamma}:=\psi^\star\rho_{\phi,\,\omega,\,\gamma}$ and $\rho_{\psi\circ\phi,\,\omega,\,\gamma}:=\phi^\star\rho_{\psi,\,\omega,\,\gamma}$.

\end{Prop}

\noindent {\it Proof.} Let $n=\mbox{dim}_\C X$.

\vspace{1ex}

(a) The pluriclosed star split hypothesis on $(\omega,\,\gamma)$ means that \begin{eqnarray*}i\partial\bar\partial\omega_{n-2} = \gamma_{n-2}\wedge\rho_{\omega,\,\gamma} \hspace{5ex} \mbox{with}\hspace{2ex} \partial\bar\partial(\star_\gamma\rho_{\omega,\,\gamma}) = 0.\end{eqnarray*} Appying $\phi^\star$ to the first equality above and using the $\gamma$-isometry hypothesis, we get the latter equality below, where the former equality follows from $\phi$ being holomorphic: \begin{eqnarray*} i\partial\bar\partial(\phi^\star\omega)_{n-2} = \phi^\star(i\partial\bar\partial\omega_{n-2}) = \gamma_{n-2}\wedge\phi^\star\rho_{\omega,\,\gamma}.\end{eqnarray*} Then, putting $\rho_{\phi,\,\omega,\,\gamma}:=\phi^\star\rho_{\omega,\,\gamma}$, we get: \begin{eqnarray*}\partial\bar\partial(\star_\gamma\rho_{\phi,\,\omega,\,\gamma}) = \partial\bar\partial\bigg(\phi^\star(\star_\gamma\rho_{\omega,\,\gamma})\bigg) = \phi^\star\bigg(\partial\bar\partial(\star_\gamma\rho_{\omega,\,\gamma})\bigg) = 0,\end{eqnarray*} where the first equality follows from Lemma \ref{Lem:star_Phi-star_commutation_isometry}, the second one follows from $\phi$ being holomorphic and the third one follows from the hypothesis. This proves the contention.

\vspace{1ex}

(b) The pluriclosed star split hypothesis on the triples $(\phi,\,\omega,\,\gamma)$ and $(\psi,\,\omega,\,\gamma)$ means that \begin{eqnarray*}i\partial\bar\partial(\phi^\star\omega)_{n-2} & = & \gamma_{n-2}\wedge\rho_{\phi,\,\omega,\,\gamma} \hspace{5ex} \mbox{with}\hspace{2ex} \partial\bar\partial(\star_\gamma\rho_{\phi,\,\omega,\,\gamma}) = 0 \\
  i\partial\bar\partial(\psi^\star\omega)_{n-2} & = & \gamma_{n-2}\wedge\rho_{\psi,\,\omega,\,\gamma} \hspace{5ex} \mbox{with}\hspace{2ex} \partial\bar\partial(\star_\gamma\rho_{\psi,\,\omega,\,\gamma}) = 0.\end{eqnarray*}

Thus, we get: \begin{eqnarray*}i\partial\bar\partial((\phi\circ\psi)^\star\omega)_{n-2} = \psi^\star\bigg(i\partial\bar\partial(\phi^\star\omega)_{n-2}\bigg) =  \psi^\star\bigg(\gamma_{n-2}\wedge\rho_{\phi,\,\omega,\,\gamma}\bigg) = \gamma_{n-2}\wedge\psi^\star\rho_{\phi,\,\omega,\,\gamma},\end{eqnarray*} where for the last equality we used the $\gamma$-isometry hypothesis on $\psi$. Now, putting $\rho_{\phi\circ\psi,\,\omega,\,\gamma}:=\psi^\star\rho_{\phi,\,\omega,\,\gamma}$, we get: \begin{eqnarray*}\partial\bar\partial(\star_\gamma\rho_{\phi\circ\psi,\,\omega,\,\gamma}) = \partial\bar\partial\bigg(\psi^\star(\star_\gamma\rho_{\phi,\,\omega,\,\gamma})\bigg) = \psi^\star\bigg(\partial\bar\partial(\star_\gamma\rho_{\phi,\,\omega,\,\gamma})\bigg) = 0,\end{eqnarray*} where for the first equality we used the $\gamma$-isometry hypothesis on $\psi$, for the second equality we used the holomorphicity of $\psi$ and the third equality follows from the triple $(\phi,\,\omega,\,\gamma)$ being pluriclosed star split.

This proves that the triple $(\phi\circ\psi,\,\omega,\,\gamma)$ is  pluriclosed star split and its associated $(1,\,1)$-form $\rho_{\phi\circ\psi,\,\omega,\,\gamma}$ satisfies the claimed equality.

The analogous claim on the triple $(\psi\circ\phi,\,\omega,\,\gamma)$ can be proved in the same way.  \hfill $\Box$

\vspace{3ex}

Using the notation: \begin{eqnarray*}Aut_\gamma(X):=\bigg\{\phi:X\longrightarrow X\,\mid\,\phi\hspace{1ex}\mbox{is a biholomorphism such that}\hspace{1ex} \phi^\star\gamma = \gamma\bigg\},\end{eqnarray*} \begin{eqnarray*}Aut_{\omega,\,\gamma}(X):=\bigg\{\phi:X\longrightarrow X\,\mid\,\phi\in Aut_\gamma(X) \hspace{1ex}\mbox{and}\hspace{1ex} (\phi,\,\omega,\,\gamma)\hspace{1ex}\mbox{is pluriclosed star split}\bigg\},\end{eqnarray*} some of the above results can be restated as

\begin{Cor}\label{Cor:subgroups} Let $X$ be a complex manifold.

\vspace{1ex}

(a) For any Hermitian metrics $\omega$ and $\gamma$ on $X$, $Aut_{\omega,\,\gamma}(X)$ is a subgroup of $Aut(X)$.

\vspace{1ex}

(b) For any {\bf pluriclosed star split} pair $(\omega,\,\gamma)$ of Hermitian metrics on $X$, we have \begin{eqnarray*}Aut_\gamma(X) = Aut_{\omega,\,\gamma}(X)\subset Aut(X).\end{eqnarray*}

\end{Cor}

 It will probably be interesting to further study the subgroup $Aut_{\omega,\,\gamma}(X)$ of automorphisms at least in the case where $X$ is compact, including as a tool to tackle the SKT-balanced conjecture.

\section{Further applications}\label{section:further-applications} The context will be mainly the one in $\S.$\ref{section:def_one-metric}, occasionally the one in $\S.$\ref{section:def_two-metrics}.

\subsection{Two types of operators associated with a Hermitian metric}\label{subsection:two-types_operators} We begin with a simple general observation that will be used in what follows. It generalises Lemma \ref{Lem:star-trace-11n-1n-1}.

\begin{Lem}\label{Lem:star_prod_star} Let $\alpha$ and $\beta$ be differential forms on an $n$-dimensional complex manifold $X$ equipped with an arbitrary Hermitian metric $\omega$ such that $\mbox{deg}\,\alpha + \mbox{deg}\,\beta = 2n$. The following equality holds: \begin{eqnarray}\label{eqn:star_prod_star}\alpha\wedge\beta = \star\alpha\wedge\star\beta,\end{eqnarray} where $\star=\star_\omega$ is the Hodge star operator induced by $\omega$. 

\end{Lem} 

\noindent {\it Proof.} Let $k$ be the degree of $\alpha$. We have: \begin{eqnarray*}\alpha\wedge\beta & = & (-1)^k\alpha\wedge\star\overline{\star\bar\beta} = (-1)^k\langle\alpha,\,\star\bar\beta\rangle\,dV_\omega = \langle\star\star\alpha,\,\star\bar\beta\rangle\,dV_\omega = \langle\star\alpha,\,\bar\beta\rangle\,dV_\omega = \star\alpha\wedge\star\beta,\end{eqnarray*} where we have used: (i)\, the definition of $\star$ requiring $u\wedge\star\bar{v} = \langle u,\,v\rangle\,dV_\omega$; (ii)\, the property $\star\star = \pm\,\mbox{Id}$ according to whether this is evaluated on even-degreed, resp. odd-degreed, forms; (iii)\, the fact that $\star$ is an isometry for the pointwise inner product $\langle\,\,,\,\,\rangle$; (iv)\,$(-1)^k = (-1)^{2n-k}$.   \hfill $\Box$

\vspace{2ex}

$\bullet$ Next, we introduce two linear operators defined pointwise that will be involved in subsequent definitions.

\begin{Def}\label{Def:T-S_def} Let $X$ be an $n$-dimensional complex manifold equipped with an arbitrary Hermitian metric $\omega$ and let \begin{eqnarray*}(\omega_{n-2}\wedge\cdot)^{-1}:\Lambda^{n-1,\,n-1}T^\star X\longrightarrow\Lambda^{1,\,1}T^\star X\end{eqnarray*} be the operator of division by $\omega_{n-2}$ of $(n-1,\,n-1)$-forms, namely the inverse of the bijection (\ref{eqn:omega_n-2_map}).

We consider the following $\omega$-dependent $\C$-linear operators: \begin{eqnarray}\label{eqn:T_omega_def}T_\omega:\Lambda^{1,\,1}T^\star X\longrightarrow\Lambda^{1,\,1}T^\star X,  \hspace{5ex} T_\omega = (\omega_{n-2}\wedge\cdot)^{-1}\circ\star_\omega,\end{eqnarray} and \begin{eqnarray}\label{eqn:S_omega_def}S_\omega:\Lambda^{n-1,\,n-1}T^\star X\longrightarrow\Lambda^{n-1,\,n-1}T^\star X,  \hspace{5ex} S_\omega = \star_\omega\circ(\omega_{n-2}\wedge\cdot)^{-1}.\end{eqnarray}

\end{Def}  

Both $T_\omega$ and $S_\omega$ are bijections. One notices the following properties right away: \begin{eqnarray}\label{eqn:S_composed-star_T}S_\omega\circ\star_\omega = \star_\omega\circ T_\omega:\Lambda^{1,\,1}T^\star X\longrightarrow\Lambda^{n-1,\,n-1}T^\star X\end{eqnarray} and \begin{eqnarray}\label{eqn:star_composed_S}\star_\omega\circ S_\omega = T_\omega\circ\star_\omega = (\omega_{n-2}\wedge\cdot)^{-1}:\Lambda^{n-1,\,n-1}T^\star X\longrightarrow\Lambda^{1,\,1}T^\star X.\end{eqnarray}

\vspace{2ex}

Resolving the division by $\omega_{n-2}$, we get the explicit expressions in the following

\begin{Lem}\label{Lem:T-omega_S-omega_explicit-formulae} (a)\, For every $(1,\,1)$-form $\alpha$, we have: \begin{eqnarray*}\label{eqn:T-omega_explicit-formula}T_\omega(\alpha) = -\alpha + \frac{1}{n-1}\,(\Lambda_\omega\alpha)\,\omega.\end{eqnarray*}

\noindent In particular, $T_\omega(\omega) = \frac{1}{n-1}\,\omega$. 

\vspace{1ex}

(b)\, For every $(n-1,\,n-1)$-form $\Omega$, we have: \begin{eqnarray*}\label{eqn:S-omega_explicit-formula}S_\omega(\Omega) = -\Omega + \frac{1}{n-1}\,\Lambda_\omega(\star_\omega\Omega)\,\omega_{n-1}.\end{eqnarray*}

\noindent In particular, $S_\omega(\omega_{n-1}) = \frac{1}{n-1}\,\omega_{n-1}$.

\end{Lem}

\noindent {\it Proof.} (a)\, Let $\alpha$ be a $(1,\,1)$-form. Taking $\star=\star_\omega$ in the Lefschetz decomposition $\alpha = \alpha_{prim} + \frac{1}{n}\,(\Lambda_\omega\alpha)\,\omega$, we get \begin{eqnarray*}\star\alpha = -\alpha_{prim}\wedge\omega_{n-2} + \bigg(\frac{1}{n(n-1)}\,(\Lambda_\omega\alpha)\,\omega\bigg)\wedge\omega_{n-2},\end{eqnarray*} where we also used the standard formula (\ref{eqn:prim-form-star-formula-gen}) for the primitive $(1,\,1)$-form $\alpha_{prim}$ and the standard equality $\star\omega = \omega_{n-1}$. We infer that \begin{eqnarray*}T_\omega(\alpha) = (\omega_{n-2}\wedge\cdot)^{-1}(\star\alpha) = -\bigg(\alpha_{prim} + \frac{1}{n}\,(\Lambda_\omega\alpha)\,\omega\bigg) + \frac{1}{n}\,\bigg(1+\frac{1}{n-1}\bigg)\,(\Lambda_\omega\alpha)\,\omega,\end{eqnarray*} which is nothing but the formula claimed under (a). 

\vspace{1ex}

(b)\, Let $\Omega$ be an $(n-1,\,n-1)$-form. There exists a unique $(1,\,1)$-form $\alpha$ such that $\Omega = \star\alpha$. Then, $\alpha = \star\Omega$ and formula (\ref{eqn:S_composed-star_T}) yields the first equality below: \begin{eqnarray*}S_\omega(\Omega) = \star(T_\omega(\alpha)) = -\star\alpha + \frac{1}{n-1}(\Lambda_\omega\alpha)\,\omega_{n-1},\end{eqnarray*} where the second equality follows from the formula for $T_\omega(\alpha)$ proved under (a). This proves the contention. \hfill $\Box$

\vspace{2ex}

$\bullet$ In the same vein, we introduce a differential operator of order two acting on $(1,\,1)$-forms. It can be seen as a higher-degree analogue of the standard Laplacian $\varphi\mapsto\Lambda_\omega(i\partial\bar\partial\varphi)$ acting on functions.

\begin{Def}\label{Def:P_omega_def} Let $X$ be an $n$-dimensional complex manifold equipped with a Hermitian metric $\omega$. We consider the following operator: \begin{eqnarray}\label{eqn:P_omega_def}P_\omega:C^\infty_{1,\,1}(X,\,\C)\longrightarrow C^\infty_{1,\,1}(X,\,\C), \hspace{5ex} P_\omega(\alpha) = (\omega_{n-2}\wedge\cdot)^{-1}(i\partial\bar\partial\alpha\wedge\omega_{n-3}).\end{eqnarray}

\end{Def}

Thus, the definition of $P_\omega$ is equivalent to the equality \begin{eqnarray}\label{eqn:P_omega_def_again}i\partial\bar\partial\alpha\wedge\omega_{n-3} = P_\omega(\alpha)\wedge\omega_{n-2}\end{eqnarray} holding for every smooth $(1,\,1)$-form $\alpha$. 

We now observe a link with the discussion in $\S.$\ref{section:def_two-metrics}, especially with the real $(1,\,1)$-form $\rho_{\omega,\,\gamma}$ associated with a pair of Hermitian metrics $(\omega,\gamma)$ via (\ref{eqn:rho_omega-gamma_def}).

\begin{Lem}\label{Lem:star-rho-P_integral-link} Let $X$ be an $n$-dimensional compact complex manifold and let $\omega$, $\gamma$ be two Hermitian metrics on $X$. Then, for every smooth $(1,\,1)$-form $\eta$ on $X$, we have: \begin{eqnarray}\label{eqn:star-rho-P_integral-link}\int\limits_X\eta\wedge\star_\gamma\rho_{\omega,\,\gamma} = \frac{n-1}{n-2}\int\limits_X P_\omega\bigg(T_\gamma(\eta)\bigg)\wedge\omega_{n-1},\end{eqnarray} where $\rho_{\omega,\,\gamma} = (\gamma_{n-2}\wedge\cdot)^{-1}(i\partial\bar\partial\omega_{n-2})$.

\end{Lem}

\noindent {\it Proof.} We have: \begin{eqnarray*}\int\limits_X\eta\wedge\star_\gamma\rho_{\omega,\,\gamma} & \stackrel{(i)}{=} & \int\limits_X\star_\gamma\star_\gamma\eta\wedge\star_\gamma\rho_{\omega,\,\gamma} \stackrel{(ii)}{=} \int\limits_X\star_\gamma\eta\wedge\rho_{\omega,\,\gamma} \stackrel{(iii)}{=} \int\limits_XT_\gamma(\eta)\wedge\gamma_{n-2}\wedge\rho_{\omega,\,\gamma} \\
& \stackrel{(iv)}{=} & \int\limits_XT_\gamma(\eta)\wedge i\partial\bar\partial\omega_{n-2} \stackrel{(v)}{=} \int\limits_Xi\partial\bar\partial T_\gamma(\eta)\wedge\omega_{n-2},\end{eqnarray*} where (i) follows from $\star_\gamma\star_\gamma = \mbox{Id}$ on even-degreed forms, (ii) follows from (\ref{eqn:star_prod_star}), (iii) follows from $T_\gamma(\eta) = (\gamma_{n-2}\wedge\cdot)^{-1}(\star_\gamma\eta)$, (iv) follows from the definition (\ref{eqn:rho_omega-gamma_def}) of $\rho_{\omega,\,\gamma}$ and (v) follows from the Stokes theorem. 

Now, writing \begin{eqnarray*}i\partial\bar\partial T_\gamma(\eta)\wedge\omega_{n-2} = \frac{1}{n-2}\,i\partial\bar\partial T_\gamma(\eta)\wedge\omega_{n-3}\wedge\omega = \frac{n-1}{n-2}\,P_\omega\bigg(T_\gamma(\eta)\bigg)\wedge\omega_{n-1},\end{eqnarray*} where the last equality follows from (\ref{eqn:P_omega_def_again}), the above sequence of equalities (i)-(v) can be continued and yields (\ref{eqn:star-rho-P_integral-link}). \hfill $\Box$

\begin{Cor}\label{Cor:pluriclosed-star-split_P_omega} Let $X$ be an $n$-dimensional compact complex manifold and let $\omega$, $\gamma$ be Hermitian metrics on $X$ such that the pair $(\omega,\,\gamma)$ is {\bf pluriclosed star split}. Then \begin{eqnarray}\label{eqn:pluriclosed-star-split_P_omega}\int\limits_X P_\omega\bigg(T_\gamma(i\partial\bar\partial\varphi)\bigg)\wedge\omega_{n-1} = 0\end{eqnarray} for every $C^\infty$ function $\varphi:X\longrightarrow\C$.

\end{Cor}

\noindent {\it Proof.} The statement follows by taking $\eta=i\partial\bar\partial\varphi$ in (\ref{eqn:star-rho-P_integral-link}) and using the Stokes theorem and Definition \ref{Def:pluriclosed-star-split_pairs} to get $\int_X i\partial\bar\partial\varphi\wedge\star_\gamma\rho_{\omega,\,\gamma} = \int_X\varphi\,i\partial\bar\partial(\star_\gamma\rho_{\omega,\,\gamma}) = 0$.  \hfill $\Box$

\vspace{3ex}

Note that the quantity featuring under the integral of (\ref{eqn:pluriclosed-star-split_P_omega}) can be easily transformed. Taking $\alpha = i\partial\bar\partial\varphi$ in (a) of Lemma \ref{Lem:T-omega_S-omega_explicit-formulae}, we get \begin{eqnarray*}T_\gamma(i\partial\bar\partial\varphi) = -i\partial\bar\partial\varphi + \frac{1}{n-1}\,\Delta_\gamma(\varphi)\,\gamma.\end{eqnarray*} Since $P_\omega(i\partial\bar\partial\varphi) = 0$, we get \begin{eqnarray*}P_\omega\bigg(T_\gamma(i\partial\bar\partial\varphi)\bigg) = \frac{1}{n-1}\,P_\omega\bigg(\Delta_\gamma(\varphi)\,\gamma\bigg).\end{eqnarray*}

\vspace{2ex}

If we choose $\varphi = f_{\omega,\,\gamma}$ in Corollary \ref{Cor:pluriclosed-star-split_P_omega}, we can think of the $(1,\,1)$-form $\Theta_{\omega,\,\gamma}:=P_\omega(T_\gamma(i\partial\bar\partial f_{\omega,\,\gamma}))$ as a kind of curvature form for the pair of metrics $(\omega,\,\gamma)$. We get the following

\begin{Cor}\label{Cor:pluriclosed-star-split_P_omega_vanishing} Let $X$ be an $n$-dimensional compact complex manifold and let $\omega$, $\gamma$ be Hermitian metrics on $X$ such that the pair $(\omega,\,\gamma)$ is {\bf pluriclosed star split}.

  If $P_\omega(T_\gamma(i\partial\bar\partial f_{\omega,\,\gamma}))\geq 0$ on $X$ or $P_\omega(T_\gamma(i\partial\bar\partial f_{\omega,\,\gamma}))\leq 0$ on $X$, then
$P_\omega(T_\gamma(i\partial\bar\partial f_{\omega,\,\gamma})) = 0$ on $X$.

\end{Cor}

\noindent {\it Proof.} Suppose that $\Theta_{\omega,\,\gamma}: = P_\omega(T_\gamma(i\partial\bar\partial f_{\omega,\,\gamma}))\geq 0$ on $X$. Then \begin{eqnarray*}\Theta_{\omega,\,\gamma}\wedge\omega_{n-1} = \Lambda_\omega(\Theta_{\omega,\,\gamma})\,\omega_n\geq 0  \hspace{3ex} \mbox{on}\hspace{1ex} X.\end{eqnarray*} Since $\int_X\Theta_{\omega,\,\gamma}\wedge\omega_{n-1} = 0$, by Corollary \ref{Cor:pluriclosed-star-split_P_omega}, we must have $\Lambda_\omega(\Theta_{\omega,\,\gamma}) = 0$. Thus, at every point of $X$, the sum of the eigenvalues of $\Theta_{\omega,\,\gamma}$ with respect to $\omega$ vanishes.

Meanwhile, each of these eigenvalues is non-negative because $\Theta_{\omega,\,\gamma}\geq 0$. Therefore, all the eigenvalues must vanish, so $\Theta_{\omega,\,\gamma} = 0$.

The case where $\Theta_{\omega,\,\gamma}\leq 0$ on $X$ is similar. \hfill $\Box$

\subsection{Properties of the differential operator $P_\omega:C^\infty_{1,\,1}(X,\,\C)\longrightarrow C^\infty_{1,\,1}(X,\,\C)$}\label{subsection:properties_P-omega_11} We start with the following observation.

\begin{Lem}\label{Lem:P-omega_trace-trace-square} Let $X$ be an $n$-dimensional complex manifold equipped with a Hermitian metric $\omega$. Suppose that $n\geq 3$. The second-order differential operator $P_\omega:C^\infty_{1,\,1}(X,\,\C)\longrightarrow C^\infty_{1,\,1}(X,\,\C)$ defined in (\ref{eqn:P_omega_def}) is given by the formula: \begin{eqnarray}\label{eqn:P-omega_trace-trace-square}P_\omega(\alpha) = \Lambda_\omega(i\partial\bar\partial\alpha) - \frac{1}{2(n-1)}\,\Lambda_\omega^2(i\partial\bar\partial\alpha)\,\omega,  \hspace{5ex} \alpha\in C^\infty_{1,\,1}(X,\,\C).\end{eqnarray}

\end{Lem}

\noindent {\it Proof.} We will prove, more generally, that for every form $\Gamma\in\Lambda^{2,\,2}T^\star X$, the pointwise formula holds: \begin{eqnarray}\label{eqn:Gamma-omega_n-3_trace-trace-square}(\omega_{n-2}\wedge\cdot)^{-1}(\Gamma\wedge\omega_{n-3}) = \Lambda_\omega(\Gamma) - \frac{1}{2(n-1)}\,\Lambda_\omega^2(\Gamma)\,\omega.\end{eqnarray} Then, (\ref{eqn:P-omega_trace-trace-square}) will follow from (\ref{eqn:Gamma-omega_n-3_trace-trace-square}) by taking $\Gamma=i\partial\bar\partial\alpha$.

Thus, we fix an arbitrary $\Gamma\in\Lambda^{2,\,2}T^\star X$. If $n\geq 4$, the Lefschetz decomposition of $\Gamma$ spells \begin{eqnarray}\label{eqn:Lefschetz_22_proof_1}\Gamma = \Gamma_{prim} + \omega\wedge V_{prim} + f\omega_2 = \Gamma_{prim} + \omega\wedge V,\end{eqnarray} where $\Gamma_{prim}$ is an $\omega$-primitive $(2,\,2)$-form, $V_{prim}$ is  an $\omega$-primitive $(1,\,1)$-form, $f$ is a function and we put $V:=V_{prim} + (1/2)f\omega$. Multiplying the last equality by $\omega_{n-3}$, we get: \begin{eqnarray*}\Gamma\wedge\omega_{n-3} = \Gamma_{prim}\wedge\omega_{n-3} + (n-2)\,\omega_{n-2}\wedge V = (n-2)\,\omega_{n-2}\wedge V,\end{eqnarray*} where the last equality follows from the previous one since $\Gamma_{prim}\wedge\omega_{n-3}=0$ thanks to $\Gamma_{prim}$ being a primitive $4$-form. Thus \begin{eqnarray}\label{eqn:Gamma-omega_n-3_trace-trace-square_proof_1}(\omega_{n-2}\wedge\cdot)^{-1}(\Gamma\wedge\omega_{n-3}) = (n-2)\,V,\end{eqnarray} so we are reduced to computing $V$.

If $n=3$, the pointwise multiplication map $\omega\wedge\cdot:\Lambda^{1,\,1}T^\star X\longrightarrow\Lambda^{2,\,2}T^\star$ is bijective, so (\ref{eqn:Lefschetz_22_proof_1}) holds with $\Gamma_{prim}=0$ and a uniquely determined $(1,\,1)$-form $V$. In particular, (\ref{eqn:Gamma-omega_n-3_trace-trace-square_proof_1}) holds as well.

On the other hand, we have $\Lambda_\omega(\Gamma) = \Lambda_\omega(\omega\wedge V) = [\Lambda_\omega,\,L_\omega](V) + \Lambda_\omega(V)\,\omega = (n-2)\,V + \Lambda_\omega(V)\,\omega$. So, taking $\Lambda_\omega$ again, we get: $\Lambda_\omega^2(\Gamma) = (n-2)\,\Lambda_\omega(V) + n\,\Lambda_\omega(V)$, hence \begin{eqnarray*}\Lambda_\omega(V) = \frac{1}{2(n-1)}\,\Lambda_\omega^2(\Gamma),  \hspace{3ex}\mbox{hence also}\hspace{3ex} \Lambda_\omega(\Gamma) = (n-2)\,V + \frac{1}{2(n-1)}\,\Lambda_\omega^2(\Gamma)\,\omega.\end{eqnarray*} Together with (\ref{eqn:Gamma-omega_n-3_trace-trace-square_proof_1}), the last equality proves (\ref{eqn:Gamma-omega_n-3_trace-trace-square}).  \hfill $\Box$

\vspace{2ex}

We note that, if one works with $V_{prim}$ and $f$ in the above proof, one gets \begin{eqnarray*}(\omega_{n-2}\wedge\cdot)^{-1}(\Gamma\wedge\omega_{n-3}) = \Lambda_\omega(\Gamma) - \frac{1}{n-1}\,\frac{\Gamma\wedge\omega_{n-2}}{\omega_n}\,\omega.\end{eqnarray*} Comparing with (\ref{eqn:Gamma-omega_n-3_trace-trace-square}), this implies \begin{eqnarray}\label{eqn:Gamma-omega_n-3_trace-trace-square_proof_2}\frac{1}{2}\,\Lambda_\omega^2(\Gamma) = \frac{\Gamma\wedge\omega_{n-2}}{\omega_n},  \hspace{5ex}  \Gamma\in\Lambda^{2,\,2}T^\star X.\end{eqnarray}

\noindent Together with (\ref{eqn:P-omega_trace-trace-square}), this implies the following property of $P_\omega$: \begin{eqnarray}\label{eqn:P_omega_alpha-wedge_omega_n-1}P_\omega(\alpha)\wedge\omega_{n-1} = \frac{n-2}{n-1}\,i\partial\bar\partial\alpha\wedge\omega_{n-2},  \hspace{5ex} \alpha\in C^\infty_{1,\,1}(X,\,\C),\end{eqnarray} which is also an immediate consequence of (\ref{eqn:P_omega_def_again}). 

\begin{Cor}\label{Cor:P-omega_trace-trace-square} Under the assumptions of Lemma \ref{Lem:P-omega_trace-trace-square}, the following equality holds: \begin{eqnarray}\label{eqn:P-omega_trace-trace-square_Cor}\Lambda_\omega(P_\omega(\alpha)) = \frac{n-2}{2(n-1)}\,\Lambda_\omega^2(i\partial\bar\partial\alpha),  \hspace{5ex} \alpha\in C^\infty_{1,\,1}(X,\,\C).\end{eqnarray}

\end{Cor}

\noindent {\it Proof.} This follows at once by taking $\Lambda_\omega$ in (\ref{eqn:P-omega_trace-trace-square}) and using the fact that $\Lambda_\omega(\omega)=n$. \hfill $\Box$


\subsubsection{Computation of $\langle\langle P_\omega(\alpha),\,\beta\rangle\rangle$ for arbitrary smooth $(1,\,1)$-forms $\alpha$, $\beta$}\label{subsubsection:computation_P_omega_alpha-alpha} From (\ref{eqn:P-omega_trace-trace-square}), we get \begin{eqnarray}\label{eqn:P-omega_alpha,alpha}\langle\langle P_\omega(\alpha),\,\beta\rangle\rangle = \langle\langle i\partial\bar\partial\alpha,\,\omega\wedge\beta\rangle\rangle - \frac{1}{2(n-1)}\,\langle\langle\Lambda_\omega^2(i\partial\bar\partial\alpha),\,\Lambda_\omega\beta\rangle\rangle,  \hspace{5ex} \alpha,\beta\in C^\infty_{1,\,1}(X,\,\C).\end{eqnarray}

The computation of the first term on the right-hand side of (\ref{eqn:P-omega_alpha,alpha}) will yield the following conclusion.

\begin{Lem}\label{Lem:P-omega_alpha,alpha_term_1} For all $\alpha,\beta\in C^\infty_{1,\,1}(X,\,\C)$, we have \begin{eqnarray}\label{eqn:P-omega_alpha,alpha_term_1}\nonumber\langle\langle i\partial\bar\partial\alpha,\,\omega\wedge\beta\rangle\rangle & = & -\langle\langle\bar\partial\alpha,\,\bar\partial\beta\rangle\rangle + \langle\langle\Lambda_\omega(\bar\partial\alpha),\,\Lambda_\omega(\bar\partial\beta)\rangle\rangle - \langle\langle i\bar\partial\alpha,\,i\omega\wedge\bar\partial\Lambda_\omega\beta\rangle\rangle \\
 & + & \langle\langle i\bar\partial\alpha,\,\star(\beta\wedge\bar\partial\omega_{n-3})\rangle\rangle - \langle\langle i\bar\partial\alpha,\,(\Lambda_\omega\beta)\,\star(\bar\partial\omega_{n-2})\rangle\rangle.\end{eqnarray}

\end{Lem}

\noindent {\it Proof.} We have: \begin{eqnarray}\label{eqn:P-omega_alpha,alpha_term_1_proof_1}\langle\langle i\partial\bar\partial\alpha,\,\omega\wedge\beta\rangle\rangle = \langle\langle i\bar\partial\alpha,\,\partial^\star(\omega\wedge\beta)\rangle\rangle = -\langle\langle i\bar\partial\alpha,\,\star\bar\partial\star(\omega\wedge\beta)\rangle\rangle,\end{eqnarray} having used the standard formula $\partial^\star = -\star\bar\partial\star$.

\vspace{2ex}

$\bullet$ We now prove the following formula: \begin{eqnarray}\label{eqn:P-omega_alpha,alpha_term_1_proof_2}\star(\omega\wedge\alpha) = -\alpha\wedge\omega_{n-3} + (\Lambda_\omega\alpha)\,\omega_{n-2},  \hspace{5ex} \alpha\in\Lambda^{1,\,1}T^\star X.\end{eqnarray}

Let $\alpha\in\Lambda^{1,\,1}T^\star X$ be fixed. The Lefschetz decomposition of $\alpha$ spells: $\alpha = \alpha_{prim} + \frac{1}{n}\,(\Lambda_\omega\alpha)\,\omega$. Multiplying by $\omega$, we get $\omega\wedge\alpha = \omega\wedge\alpha_{prim} + \frac{1}{n}\,(\Lambda_\omega\alpha)\,\omega^2$ and then taking $\star$ we get: \begin{eqnarray*}\star(\omega\wedge\alpha) & = & \Lambda_\omega(\star\alpha_{prim}) + \frac{2}{n}\,(\Lambda_\omega\alpha)\,\omega_{n-2} = \Lambda_\omega(-\alpha_{prim}\wedge\omega_{n-2}) + \frac{2}{n}\,(\Lambda_\omega\alpha)\,\omega_{n-2} \\
  & = & -\frac{1}{(n-2)!}\,[\Lambda_\omega,\,L_\omega^{n-2}](\alpha_{prim}) + \frac{2}{n}\,(\Lambda_\omega\alpha)\,\omega_{n-2} = -\alpha_{prim}\wedge\omega_{n-3} + \frac{2}{n}\,(\Lambda_\omega\alpha)\,\omega_{n-2}.\end{eqnarray*} On the other hand, multiplying the Lefschetz decomposition of $\alpha$ by $\omega_{n-3}$, we get \begin{eqnarray*}\alpha_{prim}\wedge\omega_{n-3} = \alpha\wedge\omega_{n-3} - \frac{n-2}{n}\,(\Lambda_\omega\alpha)\,\omega_{n-2}.\end{eqnarray*} Plugging this expression of $\alpha_{prim}\wedge\omega_{n-3}$ into the above expression for $\star(\omega\wedge\alpha)$, we get (\ref{eqn:P-omega_alpha,alpha_term_1_proof_2}), as claimed.

\vspace{2ex}

$\bullet$ Next, we prove the following formulae: \begin{eqnarray}\label{eqn:P-omega_alpha,alpha_term_1_proof_3}\star(\omega\wedge\eta) & = & i\omega_{n-2}\wedge\eta, \hspace{5ex} \eta\in\Lambda^{0,\,1}T^\star X, \\
\label{eqn:P-omega_alpha,alpha_term_1_proof_4}\star(\omega_{n-3}\wedge\Gamma) & = & i\omega\wedge\Lambda_\omega\Gamma - i\Gamma, \hspace{5ex} \Gamma\in\Lambda^{1,\,2}T^\star X.\end{eqnarray} 

Let $\eta\in\Lambda^{0,\,1}T^\star X$. (iii) of (\ref{eqn:standard_comm_Lambda-L-powers}) gives the first equality below: \begin{eqnarray*}\star(\omega\wedge\eta) = \Lambda_\omega(\star\eta) = i\Lambda_\omega(\omega_{n-1}\wedge\eta) = \frac{i}{(n-1)!}\,[\Lambda_\omega,\,L_\omega^{n-1}](\eta) = i\,\omega_{n-2}\wedge\eta,\end{eqnarray*} where the second equality follows from the standard formula (\ref{eqn:prim-form-star-formula-gen}) applied to the primitive form $\eta$ of bidegree $(p,\,q)=(0,\,1)$, the third equality follows from $\Lambda_\omega\eta=0$ and the last one follows from (ii) of (\ref{eqn:standard_comm_Lambda-L-powers}) applied with $r=n-1$ and $k=1$. This proves (\ref{eqn:P-omega_alpha,alpha_term_1_proof_3}).

To prove (\ref{eqn:P-omega_alpha,alpha_term_1_proof_4}), let $\Gamma\in\Lambda^{1,\,2}T^\star X$. The Lefschetz decomposition spells: $\Gamma = \Gamma_{prim} + \omega\wedge\eta$ for a $(0,\,1)$-form $\eta$ and a primitive $(1,\,2)$-form $\Gamma_{prim}$. Taking $\star$ and applying the standard formula (\ref{eqn:prim-form-star-formula-gen}) to $\Gamma_{prim}$ and (\ref{eqn:P-omega_alpha,alpha_term_1_proof_3}) to $\eta$, we get: \begin{eqnarray*}\star\Gamma & = & -i\Gamma_{prim}\wedge\omega_{n-3} + i\omega_{n-2}\wedge\eta = -i\omega_{n-3}\wedge\bigg(\Gamma_{prim} - \frac{1}{n-2}\omega\wedge\eta\bigg) = -i\omega_{n-3}\wedge\bigg(\Gamma- \frac{n-1}{n-2}\,\omega\wedge\eta\bigg),\end{eqnarray*} where the last equality follows from $\Gamma_{prim} = \Gamma - \omega\wedge\eta$. Taking $\star$ and using the equality $\star\star\Gamma = -\Gamma$, we get the first equality below: \begin{eqnarray*}i\,\star(\omega_{n-3}\wedge\Gamma) = (n-1)\star(i\omega_{n-2}\wedge\eta) + \Gamma = -(n-1)\,\omega\wedge\eta + \Gamma,\end{eqnarray*} where the last equality follows by taking $\star$ in (\ref{eqn:P-omega_alpha,alpha_term_1_proof_3}) and using the equality $\star\star(\omega\wedge\eta) = -\omega\wedge\eta$. Hence \begin{eqnarray}\label{eqn:P-omega_alpha,alpha_term_1_proof_4_aux}\star(\omega_{n-3}\wedge\Gamma) = (n-1)\,i\omega\wedge\eta -i\Gamma.\end{eqnarray} Now, to compute $\eta$, we take $\Lambda_\omega$ in the Lefschetz decomposition $\Gamma = \Gamma_{prim} + \omega\wedge\eta$. We get: $\Lambda_\omega\Gamma = [\Lambda_\omega,\,L_\omega](\eta) = (n-1)\,\eta$ (see (i) of (\ref{eqn:standard_comm_Lambda-L-powers}) for the last equality). Together with (\ref{eqn:P-omega_alpha,alpha_term_1_proof_4_aux}), this proves (\ref{eqn:P-omega_alpha,alpha_term_1_proof_4}).

\vspace{2ex}

$\bullet$ {\it Computation of $\star\bar\partial\star(\omega\wedge\alpha)$.} Applying $\star\bar\partial$ to formula (\ref{eqn:P-omega_alpha,alpha_term_1_proof_2}), we get: \begin{eqnarray*}\star\bar\partial\star(\omega\wedge\alpha) & = & -\star\bar\partial(\alpha\wedge\omega_{n-3}) + \star\bar\partial\bigg((\Lambda_\omega\alpha)\,\omega_{n-2}\bigg) \\
 & = & -\star(\bar\partial\alpha\wedge\omega_{n-3}) -\star(\alpha\wedge\bar\partial\omega_{n-3}) + (\Lambda_\omega\alpha)\,\star\bar\partial\omega_{n-2} + \star(\omega_{n-2}\wedge\bar\partial\Lambda_\omega\alpha) \\
  & = & i\bar\partial\alpha -i\omega\wedge\Lambda_\omega(\bar\partial\alpha) -\star(\alpha\wedge\bar\partial\omega_{n-3}) + (\Lambda_\omega\alpha)\,\star\bar\partial\omega_{n-2} + i\omega\wedge\bar\partial\Lambda_\omega\alpha,\end{eqnarray*} where for the last equality we used (\ref{eqn:P-omega_alpha,alpha_term_1_proof_4}) and (\ref{eqn:P-omega_alpha,alpha_term_1_proof_3}).

\vspace{2ex}

$\bullet$ Taking the $L^2$-inner product of $-i\bar\partial\alpha$ against the above expression for $\star\bar\partial\star(\omega\wedge\beta)$ (in which we substitute $\beta$ for $\alpha$) and using (\ref{eqn:P-omega_alpha,alpha_term_1_proof_1}), we get (\ref{eqn:P-omega_alpha,alpha_term_1}). \hfill $\Box$

\vspace{3ex}

The computation of the second term on the right-hand side of (\ref{eqn:P-omega_alpha,alpha}) will lead to the following conclusion.

\begin{Lem}\label{Lem:P-omega_alpha,alpha_term_2} For all $\alpha,\beta\in C^\infty_{1,\,1}(X,\,\C)$, we have \begin{eqnarray}\label{eqn:P-omega_alpha,alpha_term_2}\nonumber\langle\langle\Lambda_\omega^2(i\partial\bar\partial\alpha),\,\Lambda_\omega\beta\rangle\rangle & = & -2\,\langle\langle i\bar\partial\alpha,\,i\omega\wedge\bar\partial\Lambda_\omega\beta\rangle\rangle -2\,\langle\langle i\bar\partial\alpha,\,(\Lambda_\omega\beta)\,\star(\bar\partial\omega_{n-2})\rangle\rangle.\end{eqnarray}

\end{Lem}

\noindent {\it Proof.} We have: \begin{eqnarray*}\label{eqn:P-omega_alpha,alpha_term_2_proof_1}\nonumber\langle\langle\Lambda_\omega^2(i\partial\bar\partial\alpha),\,\Lambda_\omega\beta\rangle\rangle & = & \langle\langle i\partial\bar\partial\alpha,\,(\Lambda_\omega\beta)\,\omega^2\rangle\rangle = \bigg\langle\bigg\langle i\bar\partial\alpha,\,\partial^\star\bigg((\Lambda_\omega\beta)\,\omega^2\bigg)\bigg\rangle\bigg\rangle \\
\nonumber & = & -\bigg\langle\bigg\langle i\bar\partial\alpha,\,\star\bar\partial\star\bigg((\Lambda_\omega\beta)\,\omega^2\bigg)\bigg\rangle\bigg\rangle = -2\,\bigg\langle\bigg\langle i\bar\partial\alpha,\,\star\bar\partial\bigg((\Lambda_\omega\beta)\,\omega_{n-2}\bigg)\bigg\rangle\bigg\rangle \\
& = & -2\,\bigg\langle\bigg\langle i\bar\partial\alpha,\,\star\bigg(\bar\partial\Lambda_\omega\beta\wedge\omega_{n-2}\bigg)\bigg\rangle\bigg\rangle - 2\,\bigg\langle\bigg\langle i\bar\partial\alpha,\,\star\bigg((\Lambda_\omega\beta)\,\bar\partial\omega_{n-2}\bigg)\bigg\rangle\bigg\rangle.\end{eqnarray*} Formula (\ref{eqn:P-omega_alpha,alpha_term_1_proof_3}) with $\eta=\bar\partial\Lambda_\omega\beta$ yields $\star(\bar\partial\Lambda_\omega\beta\wedge\omega_{n-2}) = i\omega\wedge(\bar\partial\Lambda_\omega\beta)$ and the contention follows.  \hfill $\Box$

\vspace{3ex}

The computation of the first term on the r.h.s. of the equality in Lemma \ref{Lem:P-omega_alpha,alpha_term_2}  leads to

\begin{Lem}\label{Lem:P-omega_alpha,alpha_further-term} For all $\alpha, \beta\in C^\infty_{1,\,1}(X,\,\C)$, we have \begin{eqnarray}\label{eqn:P-omega_alpha,alpha_further-term}\langle\langle i\bar\partial\alpha,\,i\omega\wedge\bar\partial\Lambda_\omega\beta\rangle\rangle = -\langle\langle i(\partial+\tau)\Lambda_\omega(\bar\partial\alpha),\,\beta\rangle\rangle + \langle\langle\Lambda_\omega\bar\partial\alpha,\,\Lambda_\omega\bar\partial\beta\rangle\rangle.\end{eqnarray}

\end{Lem}

\noindent {\it Proof.} Using the commutation relation (i) of (\ref{eqn:standard-comm-rel}), we get: \begin{eqnarray*}\langle\langle i\bar\partial\alpha,\,i\omega\wedge\bar\partial\Lambda_\omega\beta\rangle\rangle & = & \langle\langle i\Lambda_\omega(\bar\partial\alpha),\,i\bar\partial\Lambda_\omega\beta\rangle\rangle = \langle\langle i\Lambda_\omega(\bar\partial\alpha),\,i\,[\bar\partial,\,\Lambda_\omega](\beta)\rangle\rangle + \langle\langle\Lambda_\omega\bar\partial\alpha,\,\Lambda_\omega\bar\partial\beta\rangle\rangle \\
 & = & -\langle\langle i\Lambda_\omega(\bar\partial\alpha),\,(\partial+\tau)^\star(\beta)\rangle\rangle + \langle\langle\Lambda_\omega\bar\partial\alpha,\,\Lambda_\omega\bar\partial\beta\rangle\rangle.\end{eqnarray*} This proves the contention.  \hfill $\Box$

\vspace{3ex}

Putting (\ref{eqn:P-omega_alpha,alpha}) and Lemmas \ref{Lem:P-omega_alpha,alpha_term_1}, \ref{Lem:P-omega_alpha,alpha_term_2} and \ref{Lem:P-omega_alpha,alpha_further-term} together, we get the following preliminary conclusion.

\begin{Lem}\label{Lem:preliminary_P-omega_alpha,alpha} For all $\alpha, \beta\in C^\infty_{1,\,1}(X,\,\C)$, we have \begin{eqnarray}\label{eqn:preliminary_P-omega_alpha,alpha}\nonumber\langle\langle P_\omega(\alpha),\,\beta\rangle\rangle & = & -\langle\langle\bar\partial^\star\bar\partial\alpha,\,\beta\rangle\rangle + \frac{1}{n-1}\,\langle\langle\Lambda_\omega(\bar\partial\alpha),\,\Lambda_\omega(\bar\partial\beta)\rangle\rangle + \bigg(1-\frac{1}{n-1}\bigg)\,\langle\langle i\partial\Lambda_\omega(\bar\partial\alpha),\,\beta\rangle\rangle \\
\nonumber & + & \langle\langle i\bar\partial\alpha,\,\star(\beta\wedge\bar\partial\omega_{n-3})\rangle\rangle - \bigg(1-\frac{1}{n-1}\bigg)\, \langle\langle i\bar\partial\alpha,\,(\Lambda_\omega\beta)\,\star(\bar\partial\omega_{n-2})\rangle\rangle \\
& + & \bigg(1-\frac{1}{n-1}\bigg)\,\langle\langle i\tau\Lambda_\omega(\bar\partial\alpha),\,\beta\rangle\rangle.\end{eqnarray}

\end{Lem}

\subsubsection{The differential operator $R_\omega:C^\infty_{1,\,1}(X,\,\C)\longrightarrow C^\infty_{1,\,1}(X,\,\C)$}\label{subsubsection:R-omega}

We will now add terms to $P_\omega$ to make it elliptic, while preserving equality (\ref{eqn:star-rho-P_integral-link}) for the resulting operator. The first such term is described in the following

\begin{Def}\label{Def:R_omega_def} Let $X$ be an $n$-dimensional complex manifold equipped with a Hermitian metric $\omega$. We consider the following differential operator: \begin{eqnarray}\label{eqn:R_omega_def}R_\omega:C^\infty_{1,\,1}(X,\,\C)\longrightarrow C^\infty_{1,\,1}(X,\,\C), \hspace{5ex} R_\omega(\alpha) = (i\partial^\star\bar\partial^\star\alpha)\,\omega.\end{eqnarray}

\end{Def}

The first observation is that, if we replace $P_\omega$ with $P_\omega + R_\omega$, equality (\ref{eqn:star-rho-P_integral-link}) remains valid. Specifically, we have

\begin{Lem}\label{Lem:R_omega_property-integral} For every $\alpha\in C^\infty_{1,\,1}(X,\,\C)$, the following equality holds: \begin{eqnarray}\label{eqn:R_omega_property-integral}\int\limits_X R_\omega(\alpha)\wedge\omega_{n-1} = 0.\end{eqnarray}

\end{Lem}

\noindent {\it Proof.} We successively get: \begin{eqnarray*}\int\limits_X R_\omega(\alpha)\wedge\omega_{n-1} = n\int\limits_X (i\partial^\star\bar\partial^\star\alpha)\,\omega_n = n\int\limits_X (i\partial^\star\bar\partial^\star\alpha)\,\star_\omega 1 = n\,\langle\langle i\partial^\star\bar\partial^\star\alpha,\,1\rangle\rangle = n\,\langle\langle i\bar\partial^\star\alpha,\,\partial 1\rangle\rangle = 0.\end{eqnarray*} 

This proves the contention.  \hfill $\Box$

\vspace{2ex}

\begin{Lem}\label{Lem:R-omega_alpha,alpha} Let $(X,\,\omega)$ be a compact Hermitian manifold. For all $\alpha, \beta\in C^\infty_{1,\,1}(X,\,\C)$, we have \begin{eqnarray}\label{eqn:R-omega_alpha,alpha}\nonumber\langle\langle R_\omega(\alpha),\,\beta\rangle\rangle = -\langle\langle\bar\partial\bar\partial^\star\alpha,\,\beta\rangle\rangle + \langle\langle i\partial^\star(\omega\wedge\bar\partial^\star\alpha),\,\beta\rangle\rangle - \langle\langle\bar\partial^\star\alpha,\,\bar\tau^\star\beta\rangle\rangle,\end{eqnarray} where $\tau:\Lambda^{p,\,q}T^\star X\longrightarrow\Lambda^{p+1,\,q}T^\star X$ is the torsion operator definied pointwise by $\tau = \tau_\omega = [\Lambda_\omega,\,\partial\omega\wedge\cdot]$.

If $\alpha$ is real, we have $\langle\langle R_\omega(\alpha),\,\alpha\rangle\rangle = -||\bar\partial^\star\alpha||^2 + \langle\langle i\partial\Lambda_\omega(\bar\partial\alpha),\,\alpha\rangle\rangle - \langle\langle\bar\partial^\star\alpha,\,\bar\tau^\star\alpha\rangle\rangle.$

\end{Lem}

\noindent {\it Proof.} Computing, we get: \begin{eqnarray*}\langle\langle R_\omega(\alpha),\,\beta\rangle\rangle & = & \langle\langle i\bar\partial^\star\alpha,\,\partial\Lambda_\omega\beta\rangle\rangle = \langle\langle i\bar\partial^\star\alpha,\,[\partial,\,\Lambda_\omega]\beta\rangle\rangle + i\,\langle\langle \bar\partial^\star\alpha,\,\Lambda_\omega(\partial\beta)\rangle\rangle \\
& = & \langle\langle i\bar\partial^\star\alpha,\,-i(\bar\partial^\star + \bar\tau^\star)\beta\rangle\rangle + i\,\langle\langle\partial^\star(\omega\wedge\bar\partial^\star\alpha),\,\beta\rangle\rangle,\end{eqnarray*} where for the last equality we used (ii) of (\ref{eqn:standard-comm-rel}). This proves the first claim.

The second claim follows from this after we notice that $\langle\langle\partial^\star(\omega\wedge\bar\partial^\star\alpha),\,\alpha\rangle\rangle = \langle\langle\alpha,\,\bar\partial\Lambda_\omega(\partial\alpha)\rangle\rangle = \overline{\overline{\langle\langle\alpha,\,\bar\partial\Lambda_\omega(\partial\alpha)\rangle\rangle}} = \overline{\langle\langle\bar\partial\Lambda_\omega(\partial\alpha),\,\alpha\rangle\rangle} = \langle\langle\partial\Lambda_\omega(\bar\partial\alpha),\,\alpha\rangle\rangle$ for every real $(1,\,1)$-form $\alpha$.  \hfill $\Box$

\vspace{2ex}

The term $\langle\langle i\partial\Lambda_\omega(\bar\partial\alpha),\,\beta\rangle\rangle$ featuring in (\ref{eqn:preliminary_P-omega_alpha,alpha}) can be transformed as follows.

\begin{Lem}\label{Lem:del-lambda-dbar-alpha_alpha} Let $(X,\,\omega)$ be a compact Hermitian manifold. For all $\alpha, \beta\in C^\infty_{1,\,1}(X,\,\C)$, we have \begin{eqnarray*}\label{eqn:del-lambda-dbar-alpha_alpha}\langle\langle i\partial\Lambda_\omega(\bar\partial\alpha),\,\beta\rangle\rangle = \langle\langle\Lambda_\omega(\bar\partial\alpha),\,\Lambda_\omega(\bar\partial\beta)\rangle\rangle - \langle\langle\bar\partial^\star\Lambda_\omega(\bar\partial\alpha)\,\omega,\,\beta\rangle\rangle - i\,\langle\langle\Lambda_\omega(\bar\partial\alpha),\,\tau^\star\beta\rangle\rangle.\end{eqnarray*}

\end{Lem}

\noindent {\it Proof.} We have: \begin{eqnarray*}\langle\langle i\partial\Lambda_\omega(\bar\partial\alpha),\,\beta\rangle\rangle & = & \langle\langle i\Lambda_\omega(\bar\partial\alpha),\,\partial^\star\beta\rangle\rangle = \langle\langle i\Lambda_\omega(\bar\partial\alpha),\,i[\Lambda_\omega,\bar\partial]\beta - \tau^\star\beta\rangle\rangle \\
 & = & \langle\langle\Lambda_\omega(\bar\partial\alpha),\,\Lambda_\omega(\bar\partial\beta)\rangle\rangle - \langle\langle\Lambda_\omega(\bar\partial\alpha),\,\bar\partial\Lambda_\omega\beta\rangle\rangle - i\,\langle\langle\Lambda_\omega(\bar\partial\alpha),\,\tau^\star\beta\rangle\rangle.\end{eqnarray*} The contention follows from this after we notice that $\langle\langle\Lambda_\omega(\bar\partial\alpha),\,\bar\partial\Lambda_\omega\beta\rangle\rangle = \langle\langle\omega\wedge\bar\partial^\star\Lambda_\omega(\bar\partial\alpha),\,\beta\rangle\rangle$. \hfill $\Box$

\vspace{2ex}

Putting together Lemmas \ref{Lem:preliminary_P-omega_alpha,alpha}, \ref{Lem:R-omega_alpha,alpha} and \ref{Lem:del-lambda-dbar-alpha_alpha}, we get

\begin{Lem}\label{Lem:P-omega_alpha,alpha_1} Let $(X,\omega)$ be a compact Hermitian manifold. For all $\alpha, \beta\in C^\infty_{1,\,1}(X,\,\C)$, we have \begin{eqnarray}\label{eqn:P-omega_alpha,alpha_1}\nonumber\langle\langle (P_\omega + R_\omega)(\alpha),\,\beta\rangle\rangle & = & -\langle\langle\Delta''\alpha,\,\beta\rangle\rangle \\
\nonumber & + & \langle\langle i\partial\Lambda_\omega(\bar\partial\alpha),\,\beta\rangle\rangle + \langle\langle i\partial^\star(\omega\wedge\bar\partial^\star\alpha),\,\beta\rangle\rangle + \frac{1}{n-1}\,\langle\langle\bar\partial^\star\Lambda_\omega(\bar\partial\alpha)\,\omega,\,\beta\rangle\rangle \\
\nonumber & + & \langle\langle i\bar\partial\alpha,\,\star(\beta\wedge\bar\partial\omega_{n-3})\rangle\rangle - \bigg(1-\frac{1}{n-1}\bigg)\, \langle\langle i\bar\partial\alpha,\,(\Lambda_\omega\beta)\,\star(\bar\partial\omega_{n-2})\rangle\rangle \\
\nonumber & + & \bigg(1-\frac{1}{n-1}\bigg)\,\langle\langle i\tau\Lambda_\omega(\bar\partial\alpha),\,\beta\rangle\rangle + \frac{1}{n-1}\,i\,\langle\langle\Lambda_\omega(\bar\partial\alpha),\,\tau^\star\beta\rangle\rangle \\
& - & \langle\langle\bar\partial^\star\alpha,\,\bar\tau^\star\beta\rangle\rangle,\end{eqnarray} where $\Delta'':=\bar\partial\bar\partial^\star + \bar\partial^\star \bar\partial$ is the $\bar\partial$-Laplacian induced by the metric $\omega$.

\end{Lem}

\noindent {\it Proof.} The conclusion is straightforward after splitting the term $(1-(1/(n-1)))\,\langle\langle i\partial\Lambda_\omega(\bar\partial\alpha),\,\beta\rangle\rangle$ on the first line of (\ref{eqn:preliminary_P-omega_alpha,alpha}) into the difference between $\langle\langle i\partial\Lambda_\omega(\bar\partial\alpha),\,\beta\rangle\rangle$ and $(1/(n-1))\,\langle\langle i\partial\Lambda_\omega(\bar\partial\alpha),\,\beta\rangle\rangle$ and expressing the latter part using Lemma \ref{Lem:del-lambda-dbar-alpha_alpha}.  \hfill $\Box$

\vspace{2ex}

\begin{Cor}\label{Cor:Q_omega_elliptic} Let $(X,\omega)$ be a compact Hermitian manifold. The second-order differential operator $Q_\omega:C^\infty_{1,\,1}(X,\,\C)\longrightarrow C^\infty_{1,\,1}(X,\,\C)$ defined by \begin{eqnarray}\label{eqn:Q_omega_def}Q_\omega(\alpha): = P_\omega(\alpha) + R_\omega(\alpha) - i\partial\Lambda_\omega(\bar\partial\alpha) - i\partial^\star(\omega\wedge\bar\partial^\star\alpha) - \frac{1}{n-1}\,\bar\partial^\star\Lambda_\omega(\bar\partial\alpha)\,\omega,\end{eqnarray} is {\bf elliptic}.

When $\omega$ is {\bf K\"ahler}, we even have $Q_\omega = -\Delta''_\omega$.

\end{Cor}

\noindent {\it Proof.} Except for $\Delta''=\Delta''_\omega $, the only terms on the r.h.s. of (\ref{eqn:P-omega_alpha,alpha_1}) that involve expressions of order two in $\alpha$ are the three terms on the second line. They have been incorporated into $Q_\omega$. Thus, \begin{eqnarray}\label{eqn:Q_omega_obs_1}Q_\omega = -\Delta''_\omega + l.o.t.,\end{eqnarray} where l.o.t. stands for terms of order $\leq 1$. In other words, $Q_\omega$ has the same principal part as the elliptic operator $-\Delta''_\omega$. It is therefore elliptic.

 We also note that all the terms (all of which are of order $\leq 1$) on lines $2$-$4$ on the r.h.s. of (\ref{eqn:P-omega_alpha,alpha_1}) vanish if $\omega$ is K\"ahler. This implies the last statement of the corollary in the K\"ahler case.  \hfill $\Box$

\vspace{2ex}

An immediate consequence of (\ref{eqn:Q_omega_def}) is that the restrictions to $\ker\Delta''_\omega$ of $P_\omega$ and $Q_\omega$ coincide.

\begin{Cor}\label{Cor:P_omega_Q_omega_restrictions} Let $(X,\omega)$ be a compact Hermitian manifold. Then \begin{eqnarray*}Q_\omega(\alpha) = P_\omega(\alpha) \end{eqnarray*} for every $\alpha\in\ker\Delta''_\omega$.

\end{Cor}

\noindent {\it Proof.} This follows at once from (\ref{eqn:Q_omega_def}) and from the well-known fact that $\ker\Delta'' = \ker\bar\partial\cap\ker\bar\partial^\star$.  \hfill $\Box$

\vspace{2ex}

The next observation is that equality (\ref{eqn:star-rho-P_integral-link}) remains valid under certain conditions if we replace $P_\omega$ with $Q_\omega$ since the expressions $R_\omega(\alpha)$ (for which we have Lemma \ref{Lem:R_omega_property-integral}), $i\partial\Lambda_\omega(\bar\partial\alpha)$, $i\partial^\star(\omega\wedge\bar\partial^\star\alpha)$ and $\bar\partial^\star\Lambda_\omega(\bar\partial\alpha)\,\omega$ have trivial contributions to the integral on the r.h.s. of (\ref{eqn:star-rho-P_integral-link}). Specifically, we have

\begin{Lem}\label{Lem:extra-terms_order-two_integral} Let $(X,\,\omega)$ be a compact Hermitian manifold. For every $\alpha\in C^\infty_{1,\,1}(X,\,\C)$, we have \begin{eqnarray}\label{eqn:extra-terms_order-two_integral_1}\int\limits_X\bar\partial^\star\Lambda_\omega(\bar\partial\alpha)\,\omega\wedge\omega_{n-1} = 0.\end{eqnarray} If $\omega$ is {\bf balanced}, we also have \begin{eqnarray}\label{eqn:extra-terms_order-two_integral_2}\int\limits_X i\partial\Lambda_\omega(\bar\partial\alpha)\wedge\omega_{n-1} = 0 \hspace{3ex}\mbox{and}\hspace{3ex} \int\limits_X i\partial^\star(\omega\wedge\bar\partial^\star\alpha)\wedge\omega_{n-1} = 0.\end{eqnarray}

\end{Lem}

\noindent {\it Proof.} To prove (\ref{eqn:extra-terms_order-two_integral_1}), we note that \begin{eqnarray*}\int\limits_X\bar\partial^\star\Lambda_\omega(\bar\partial\alpha)\,\omega\wedge\omega_{n-1} & = & n\int\limits_X\bar\partial^\star\Lambda_\omega(\bar\partial\alpha)\,\omega_n = n\int\limits_X\bar\partial^\star\Lambda_\omega(\bar\partial\alpha)\,\star_\omega 1 = n\,\langle\langle\bar\partial^\star\Lambda_\omega(\bar\partial\alpha),\, 1\rangle\rangle \\
 & = & n\,\langle\langle\Lambda_\omega(\bar\partial\alpha),\, \bar\partial(1)\rangle\rangle = 0.\end{eqnarray*}

To prove the first part of (\ref{eqn:extra-terms_order-two_integral_2}), we start by using the balanced assumption on $\omega$ to get the first equality below: \begin{eqnarray*}\int\limits_X i\partial\Lambda_\omega(\bar\partial\alpha)\wedge\omega_{n-1} = i\int\limits_X\partial\bigg(\Lambda_\omega(\bar\partial\alpha)\wedge\omega_{n-1}\bigg) = 0,\end{eqnarray*} where the second equality follows from the Stokes theorem.  

To prove the second part of (\ref{eqn:extra-terms_order-two_integral_2}), we use the equality $\omega_{n-1} = \star\omega$ to get the first equality below: \begin{eqnarray*}\int\limits_X i\partial^\star(\omega\wedge\bar\partial^\star\alpha)\wedge\omega_{n-1} = \int\limits_X i\partial^\star(\omega\wedge\bar\partial^\star\alpha)\wedge\star\omega = \langle\langle i\partial^\star(\omega\wedge\bar\partial^\star\alpha),\,\omega\rangle\rangle = \langle\langle i\bar\partial^\star\alpha),\,\Lambda_\omega(\partial\omega)\rangle\rangle = 0,\end{eqnarray*} where the last equality follows from the balanced assumption on $\omega$ which is equivalent to $\partial\omega$ being $\omega$-primitive (see e.g. proof of Proposition \ref{Prop:astheno-balanced-K}, meaning that $ \Lambda_\omega(\partial\omega) = 0$. \hfill $\Box$

\vspace{3ex}

One of the conclusions of this discussion is the following version of Lemma \ref{Lem:star-rho-P_integral-link} where $P_\omega$ is replaced with the elliptic operator $Q_\omega$. 

\begin{Lem}\label{Lem:star-rho-Q_integral-link} Let $X$ be an $n$-dimensional compact complex manifold and let $\omega$, $\gamma$ be two Hermitian metrics on $X$. Suppose that $\omega$ is {\bf balanced}.

  Then, for every smooth $(1,\,1)$-form $\eta$ on $X$, we have: \begin{eqnarray}\label{eqn:star-rho-Q_integral-link}\int\limits_X\eta\wedge\star_\gamma\rho_{\omega,\,\gamma} = \frac{n-1}{n-2}\int\limits_X Q_\omega\bigg(T_\gamma(\eta)\bigg)\wedge\omega_{n-1},\end{eqnarray} where $\rho_{\omega,\,\gamma} = (\gamma_{n-2}\wedge\cdot)^{-1}(i\partial\bar\partial\omega_{n-2})$.

\end{Lem}

\noindent {\it Proof.} Corollary \ref{Cor:Q_omega_elliptic} and Lemmas \ref{Lem:R_omega_property-integral} and \ref{Lem:extra-terms_order-two_integral} show that \begin{eqnarray*} \int\limits_X\bigg(Q_\omega(\alpha) - P_\omega(\alpha)\bigg)\wedge\omega_{n-1} = 0\end{eqnarray*} for every $(1,\,1)$-form $\alpha$ whenever the metric $\omega$ is balanced. The contention follows from this and from Lemma \ref{Lem:star-rho-P_integral-link}.  \hfill $\Box$

  \vspace{2ex}

  It goes without saying that Corollaries \ref{Cor:pluriclosed-star-split_P_omega} and \ref{Cor:pluriclosed-star-split_P_omega_vanishing} have obvious analogues when $P_\omega$ is replaced with $Q_\omega$ and $\omega$ is balanced.

\vspace{3ex}

We now notice that formula (\ref{eqn:Q_omega_def}) simplifies when $\alpha=\omega$.

\begin{Cor}\label{Cor:Q_omega_omega} Let $(X,\omega)$ be a compact Hermitian manifold. The following equality holds: \begin{eqnarray}\label{eqn:Q_omega_omega}Q_\omega(\omega) = P_\omega(\omega) + \frac{n}{n-1}\,R_\omega(\omega) + \partial\partial^\star\omega - i\partial^\star(\omega\wedge\bar\partial^\star\omega).\end{eqnarray} 

If $\omega$ is {\bf balanced}, then $Q_\omega(\omega) = P_\omega(\omega)$.

\end{Cor}

\noindent {\it Proof.} We transform, one by one, the third and the fifth terms in (\ref{eqn:Q_omega_def}) when $\alpha =\omega$. We get: \begin{eqnarray*}i\partial\Lambda_\omega(\bar\partial\omega) = \partial\bigg(i\,[\Lambda_\omega,\,\bar\partial]\,\omega\bigg) = \partial\partial^\star\omega + \partial\tau^\star\omega = -\partial\partial^\star\omega,\end{eqnarray*} where the first equality follows from $\Lambda_\omega(\omega)=n$ (which implies $\bar\partial\Lambda_\omega(\omega)=0$), the second equality follows from the commutation relation (i) of (\ref{eqn:standard-comm-rel}) and the third equality follows from (i) of (\ref{eqn:Michelson_torsion-form}).

Meanwhile, we get: \begin{eqnarray*}\bar\partial^\star\Lambda_\omega(\bar\partial\omega) = \bar\partial^\star\bigg([\Lambda_\omega,\,\bar\partial]\,\omega\bigg) = -i\bar\partial^\star\partial^\star\omega - i\bar\partial^\star\tau^\star\omega = i\bar\partial^\star\partial^\star\omega,\end{eqnarray*} where the arguments given above were repeated to transform the quantity $[\Lambda_\omega,\,\bar\partial]\,\omega$.

Putting these pieces of information together and using (\ref{eqn:Q_omega_def}) combined with the definition (\ref{eqn:R_omega_def}) of $R_\omega$, we get (\ref{eqn:Q_omega_omega}). 

If $\omega$ is balanced, $\partial^\star\omega = 0$ and $\bar\partial^\star\omega = 0$, hence also $R_\omega(\omega) = 0$. Together with (\ref{eqn:Q_omega_omega}), this proves the last claim. \hfill $\Box$

\subsection{Link between the function $f_\omega$ and the operator $P_\omega$}\label{subsection:f_omega_P-omega} Let $X$ be an $n$-dimensional complex manifold on which an arbitrary Hermitian metric $\omega$ has been fixed. We will compute the associated $C^\infty$ function $f_\omega:X\to\R$ defined in (\ref{eqn:function_f_omega}). We have: \begin{eqnarray}\label{eqn:function_f_omega_computation_1}\nonumber f_\omega & = & \frac{\omega\wedge i\partial\bar\partial\omega_{n-2}}{\omega_n} = \frac{\star(i\partial\bar\partial\omega_{n-2})\wedge\omega_{n-1}}{\omega_n}  \\ 
& = & \frac{\star\bigg(i\partial(\omega_{n-3}\wedge\bar\partial\omega)\bigg)\wedge\omega_{n-1}}{\omega_n} = \frac{\star(i\partial\bar\partial\omega\wedge\omega_{n-3})\wedge\omega_{n-1}}{\omega_n} + \frac{\star(i\partial\omega\wedge\bar\partial\omega\wedge\omega_{n-4})\wedge\omega_{n-1}}{\omega_n},\end{eqnarray} where the second equality follows from Lemma \ref{Lem:star_prod_star}. 

Since $i\partial\bar\partial\omega$ is a $(2,\,2)$-form and $i\partial\omega\wedge\bar\partial\omega$ is a $(3,\,3)$-form, we will need the results of the following two computations.

\begin{Lem}\label{Lem:star-Gamma_22-omega_n-3} For any $(2,\,2)$-form $\Gamma$ on an $n$-dimensional Hermitian manifold $(X,\,\omega)$, we have: \begin{eqnarray}\label{eqn:star-Gamma_22-omega_n-3}\star(\Gamma\wedge\omega_{n-3}) = -\Lambda_\omega\Gamma + \frac{1}{2}\,(\Lambda_\omega^2\Gamma)\,\omega\end{eqnarray} if $n\geq 3$.

\end{Lem} 

\noindent {\it Proof.} We saw in the proof of Lemma \ref{Lem:P-omega_trace-trace-square} that if we use the Lefschetz decomposition of $\Gamma$ to write \begin{eqnarray*}\Gamma = \Gamma_{prim} + \omega\wedge V_{prim} + f\omega_2 = \Gamma_{prim} + \omega\wedge V,\end{eqnarray*} where $\Gamma_{prim}$ is an $\omega$-primitive $(2,\,2)$-form if $n\geq 4$ and $\Gamma_{prim} = 0$ if $n=3$, $V_{prim}$ is  an $\omega$-primitive $(1,\,1)$-form, $f$ is function and we put $V:=V_{prim} + (1/2)f\omega$, then we have: \begin{eqnarray}\label{eqn:star-Gamma_22-omega_n-3_proof_1}V = \frac{1}{n-2}\,\Lambda_\omega\Gamma - \frac{1}{2(n-1)(n-2)}\,(\Lambda_\omega^2\Gamma)\,\omega.\end{eqnarray}

On the other hand, we have: \begin{eqnarray*}\Gamma\wedge\omega_{n-3} = (n-2)\,V\wedge\omega_{n-2} = (n-2)\,V_{prim}\wedge\omega_{n-2} + \frac{(n-2)(n-1)}{n}\,(\Lambda_\omega V)\,\omega_{n-1},\end{eqnarray*} where $V = V_{prim} + (1/n)\,(\Lambda_\omega V)\,\omega$ is the Lefschetz decomposition of $V$. (We used the fact that $\Gamma_{prim}\wedge\omega_{n-3} = 0$, a consequence of $\Gamma_{prim}$ being a primitive $4$-form.) Since $\star V_{prim} = -V_{prim}\wedge\omega_{n-2}$ (see the standard formula (\ref{eqn:prim-form-star-formula-gen})), we infer that $\star(V_{prim}\wedge\omega_{n-2}) = -V_{prim}$, hence \begin{eqnarray}\label{eqn:star-Gamma_22-omega_n-3_proof_2}\nonumber\star(\Gamma\wedge\omega_{n-3}) & = & -(n-2)\,V_{prim} + \frac{(n-2)(n-1)}{n}\,(\Lambda_\omega V)\,\omega \\
\nonumber & = & -(n-2)\,\bigg(V_{prim} + \frac{1}{n}\,(\Lambda_\omega V)\,\omega\bigg) + (n-2)\,(\Lambda_\omega V)\,\omega \\
 & = & -(n-2)\,\bigg(V - (\Lambda_\omega V)\,\omega\bigg).\end{eqnarray}

Now, taking $\Lambda_\omega$ in (\ref{eqn:star-Gamma_22-omega_n-3_proof_1}), we get $\Lambda_\omega V = \frac{1}{2(n-1)}\,\Lambda_\omega^2\Gamma$. Plugging this and (\ref{eqn:star-Gamma_22-omega_n-3_proof_1}) into (\ref{eqn:star-Gamma_22-omega_n-3_proof_2}), we get (\ref{eqn:star-Gamma_22-omega_n-3}), as claimed.  \hfill $\Box$

\begin{Lem}\label{Lem:star-Omega_33-omega_n-4} For any $(3,\,3)$-form $\Omega$ on an $n$-dimensional Hermitian manifold $(X,\,\omega)$, we have: \begin{eqnarray}\label{eqn:star-Omega_33-omega_n-4}\star(\Omega\wedge\omega_{n-4}) = -\frac{1}{2!}\Lambda_\omega^2\Omega + \frac{1}{3!}\,(\Lambda_\omega^3\Omega)\,\omega\end{eqnarray} if $n\geq 4$.

\end{Lem} 

\noindent {\it Proof.} If $n\geq 6$, let $\Omega = \Omega_{prim} + \Omega_{1,\,prim}\wedge\omega + \Omega_{2,\,prim}\wedge\omega_2 + f\omega_3$ be the Lefschetz decomposition of $\Omega$. Thus, $\Omega_{1,\,prim}$ and $\Omega_{2,\,prim}$ are primitive forms of respective bidegrees $(2,\,2)$ and $(1,\,1)$, while $f$ is a function. The same decomposition of $\Omega$ holds with $\Omega_{prim} = 0$ when $n\in\{4,\,5\}$. Indeed, the pointwise multiplication map $\omega\wedge\cdot:\Lambda^{2,\,2}T^\star X\longrightarrow\Lambda^{3,\,3}T^\star X$ is surjective and non-injective when $n=4$ (so the decomposition of $\Omega$ is not unique), while it is bijective (yielding a unique decomposition of $\Omega$) when $n=5$. This means that $\Omega$ can be divided by $\omega$ when $n=4$ or $n=5$. The Lefschetz decomposition can then be applied to the resulting quotient $(2,\,2)$-form. 

Multiplying by $\omega_{n-4}$ and then taking $\star$, we get \begin{eqnarray}\label{eqn:star-Omega_33-omega_n-4_proof_1}\nonumber\star(\Omega\wedge\omega_{n-4}) & = & \star\bigg(\frac{(n-3)(n-2)}{2}\,\Omega_{2,\,prim}\wedge\omega_{n-2} + \frac{(n-3)(n-2)(n-1)}{6}\,f\omega_{n-1}\bigg)  \\
 & = & -\frac{(n-3)(n-2)}{2}\,\Omega_{2,\,prim} + \frac{(n-3)(n-2)(n-1)}{6}\,f\omega,\end{eqnarray} where the standard formula (\ref{eqn:prim-form-star-formula-gen}) was used to get: $\star\Omega_{2,\,prim} = -\Omega_{2,\,prim}\wedge\omega_{n-2}$, hence $\star(\Omega_{2,\,prim}\wedge\omega_{n-2}) = - \Omega_{2,\,prim}$.

\vspace{1ex}

On the other hand, we can take $\Lambda_\omega$ successively in the Lefschetz decomposition of $\Omega$ to compute $\Omega_{2,\,prim}$ and $f$. After one application of $\Lambda_\omega$, we get: \begin{eqnarray*}\Lambda_\omega\Omega & = & [\Lambda_\omega,\,L_\omega](\Omega_{1,\,prim}) + \frac{1}{2}\,[\Lambda_\omega,\,L_\omega^2](\Omega_{2,\,prim}) + \frac{1}{3!}\,[\Lambda_\omega,\,L_\omega^3](f) \\
& = & (n-4)\,\Omega_{1,\,prim} + (n-3)\,\omega\wedge\Omega_{2,\,prim} + \frac{(n-2)}{2}f\omega^2,\end{eqnarray*} where the last equality follows from (i) and (ii) of (\ref{eqn:standard_comm_Lambda-L-powers}).

Taking $\Lambda_\omega$ again, we get: \begin{eqnarray}\label{eqn:star-Omega_33-omega_n-4_proof_2}\nonumber\Lambda_\omega^2\Omega & = & (n-3)\,[\Lambda_\omega,\,L_\omega](\Omega_{2,\,prim}) + \frac{(n-2)}{2}\,[\Lambda_\omega,\,L_\omega^2](f) \\
 & = & (n-3)(n-2)\,\Omega_{2,\,prim} + (n-2)(n-1)\,f\omega,\end{eqnarray} where the last equality follows from (i) and (ii) of (\ref{eqn:standard_comm_Lambda-L-powers}).

Taking $\Lambda_\omega$ one final time, we get: $\Lambda_\omega^3\Omega = (n-2)(n-1)n\,f$, hence \begin{eqnarray}\label{eqn:star-Omega_33-omega_n-4_proof_3}f = \frac{1}{(n-2)(n-1)n}\,\Lambda_\omega^3\Omega.\end{eqnarray} Plugging this into (\ref{eqn:star-Omega_33-omega_n-4_proof_2}), we get \begin{eqnarray}\label{eqn:star-Omega_33-omega_n-4_proof_4}\Omega_{2,\,prim} = \frac{1}{(n-3)(n-2)}\,\Lambda_\omega^2\Omega - \frac{1}{(n-3)(n-2)n}\,(\Lambda_\omega^3\Omega)\,\omega.\end{eqnarray}

\vspace{1ex}

Putting together (\ref{eqn:star-Omega_33-omega_n-4_proof_1}), (\ref{eqn:star-Omega_33-omega_n-4_proof_3}) and (\ref{eqn:star-Omega_33-omega_n-4_proof_4}), we get (\ref{eqn:star-Omega_33-omega_n-4}), as claimed.  \hfill $\Box$

\vspace{3ex}

We are now in a position to conclude the following

\begin{Prop}\label{Prop:formula_function_f_omega} Let $(X,\,\omega)$ be an $n$-dimensional Hermitian manifold. The $C^\infty$ function $f_\omega:X\to\R$ associated with $\omega$ as defined in (\ref{eqn:function_f_omega}) is given by the formulae: \begin{eqnarray}\label{eqn:formula_function_f_omega_1}f_\omega & = & \frac{n-2}{2!}\,\Lambda_\omega^2(i\partial\bar\partial\omega) + \frac{n-3}{3!}\,\Lambda_\omega^3(i\partial\omega\wedge\bar\partial\omega) \\
\label{eqn:formula_function_f_omega_2} & = & (n-1)\,\Lambda_\omega(P_\omega(\omega)) + \frac{n-3}{3!}\,\Lambda_\omega^3(i\partial\omega\wedge\bar\partial\omega).\end{eqnarray}

\end{Prop}

\noindent {\it Proof.} Thanks to (\ref{eqn:P-omega_trace-trace-square_Cor}), (\ref{eqn:formula_function_f_omega_2}) is an immediate consequence of (\ref{eqn:formula_function_f_omega_1}). 

To prove (\ref{eqn:formula_function_f_omega_1}), we take $\Gamma= i\partial\bar\partial\omega$ in (\ref{eqn:star-Gamma_22-omega_n-3}) and $\Omega=i\partial\omega\wedge\bar\partial\omega$ in (\ref{eqn:star-Omega_33-omega_n-4}) and use formula (\ref{eqn:function_f_omega_computation_1}) obtained earlier for $f_\omega$. \hfill $\Box$

\subsection{Link between the $(1,\,1)$-form $\rho_\omega$ and the operator $P_\omega$}\label{subsection:rho_omega_P-omega} Let $X$ be an $n$-dimensional complex manifold on which an arbitrary Hermitian metric $\omega$ has been fixed. We will compute the associated $C^\infty$ $(1,\,1)$-form $\rho_\omega$ defined in (\ref{eqn:rho_omega_def}). We already know that it is given by formula (\ref{eqn:rho_omega_formula}), which can be rewritten as \begin{eqnarray}\label{eqn:rho_omega_formula_bis}\rho_\omega = \frac{1}{n-1}\,f_\omega\,\omega + i\bar\partial^\star\partial^\star\omega_2.\end{eqnarray} Since $f_\omega$ was computed in (\ref{eqn:formula_function_f_omega_2}), we implicitly get a formula for $\rho_\omega$. However, we will prove a different formula better adapted to our purposes.

\begin{Prop}\label{Prop:formula_rho_omega_P-omega} Let $(X,\,\omega)$ be an $n$-dimensional Hermitian manifold. The $C^\infty$ $(1,\,1)$-form $\rho_\omega$ associated with $\omega$ as defined in (\ref{eqn:rho_omega_def}) is given by the formula: \begin{eqnarray}\label{eqn:formula_rho_omega_P-omega}\rho_\omega = P_\omega(\omega) + \frac{1}{2}\,\Lambda_\omega^2(i\partial\omega\wedge\bar\partial\omega) - \frac{1}{3(n-1)}\,\Lambda_\omega^3(i\partial\omega\wedge\bar\partial\omega)\,\omega.\end{eqnarray}

\end{Prop}

\noindent {\it Proof.} Definition (\ref{eqn:rho_omega_def}) gives the first equality below: \begin{eqnarray}\label{eqn:formula_rho_omega_P-omega_proof_1}\nonumber\rho_\omega & = & (\omega_{n-2}\wedge\cdot)^{-1}(i\partial\bar\partial\omega_{n-2}) = (\omega_{n-2}\wedge\cdot)^{-1}(i\partial\bar\partial\omega\wedge\omega_{n-3}) + (\omega_{n-2}\wedge\cdot)^{-1}(i\partial\omega\wedge\bar\partial\omega\wedge\omega_{n-4}) \\
& = & P_\omega(\omega) + (\omega_{n-2}\wedge\cdot)^{-1}(i\partial\omega\wedge\bar\partial\omega\wedge\omega_{n-4}).\end{eqnarray} 

Since $i\partial\omega\wedge\bar\partial\omega$ is a $(3,\,3)$-form, we will need the following analogue in this bidegree of the pointwise formula (\ref{eqn:Gamma-omega_n-3_trace-trace-square}).

\begin{Lem}\label{Lem:Omega-omega_n-4_trace-trace-square-cube} For every form $\Omega\in\Lambda^{3,\,3}T^\star X$ on an $n$-dimensional complex manifold $X$ with $n\geq 4$, the pointwise formula holds: \begin{eqnarray}\label{eqn:Omega-omega_n-4_trace-trace-square-cube}(\omega_{n-2}\wedge\cdot)^{-1}(\Omega\wedge\omega_{n-4}) = \frac{1}{2}\Lambda_\omega^2(\Omega) - \frac{1}{3(n-1)}\,\Lambda_\omega^3(\Omega)\,\omega.\end{eqnarray}

\end{Lem}

\noindent {\it Proof.} Let $\Omega$ be a $(3,\,3)$-form. When $n\geq 6$, from its Lefschetz decomposition we get \begin{eqnarray*}\Omega = \Omega_{prim} + \omega\wedge\Gamma,\end{eqnarray*} where $\Omega_{prim}$ is an $\omega$-primitive $(3,\,3)$-form and $\Gamma$ is a (not necessarily primitive) $(2,\,2)$-form. As explained in the proof of Lemma \ref{Lem:star-Omega_33-omega_n-4}, the same decomposition of $\Omega$ holds with $\Omega_{prim} = 0$ when $n\in\{4,\,5\}$. 

Then, $\Omega_{prim}\wedge\omega_{n-5} = 0$, hence also $\Omega_{prim}\wedge\omega_{n-4} = 0$, so $\Omega\wedge\omega_{n-4} = (n-3)\,\Gamma\wedge\omega_{n-3}$. Thus, \begin{eqnarray}\label{eqn:Omega-omega_n-4_trace-trace-square-cube_proof_1}(\omega_{n-2}\wedge\cdot)^{-1}(\Omega\wedge\omega_{n-4}) = (n-3)\,(\omega_{n-2}\wedge\cdot)^{-1}(\Gamma\wedge\omega_{n-3}) = (n-3)\,\Lambda_\omega(\Gamma) - \frac{n-3}{2(n-1)}\,\Lambda_\omega^2(\Gamma)\,\omega,\end{eqnarray} where the last equality follows from (\ref{eqn:Gamma-omega_n-3_trace-trace-square}).

We are thus reduced to computing $\Lambda_\omega(\Gamma)$ and $\Lambda_\omega^2(\Gamma)$. Taking $\Lambda_\omega$ in the above formula defining $\Gamma$ and using $\Lambda_\omega(\Omega_{prim}) = 0$, we get $\Lambda_\omega(\Omega) = [\Lambda_\omega,\,L_\omega](\Gamma) + \omega\wedge\Lambda_\omega(\Gamma) = (n-4)\,\Gamma + \omega\wedge\Lambda_\omega(\Gamma)$, where the last equality follows from (i) of (\ref{eqn:standard_comm_Lambda-L-powers}). Taking $\Lambda_\omega$ again, we get \begin{eqnarray*}\Lambda_\omega^2(\Omega) = (n-4)\,\Lambda_\omega(\Gamma) + [\Lambda_\omega,\,L_\omega](\Lambda_\omega(\Gamma)) + \omega\wedge\Lambda_\omega^2(\Gamma) = 2(n-3)\,\Lambda_\omega(\Gamma) + \Lambda_\omega^2(\Gamma)\,\omega,\end{eqnarray*} where the last equality follows again from (i) of (\ref{eqn:standard_comm_Lambda-L-powers}). Taking $\Lambda_\omega$ again, we get \begin{eqnarray*}\Lambda_\omega^3(\Omega)= 2(n-3)\,\Lambda_\omega^2(\Gamma) + n\,\Lambda_\omega^2(\Gamma) = 3(n-2)\,\Lambda_\omega^2(\Gamma).\end{eqnarray*}

This yields \begin{eqnarray}\label{eqn:Omega-omega_n-4_trace-trace-square-cube_proof_2}\Lambda_\omega^2(\Gamma) = \frac{1}{3(n-2)}\,\Lambda_\omega^3(\Omega).\end{eqnarray} Using this and another of the above formulae, we get: \begin{eqnarray}\label{eqn:Omega-omega_n-4_trace-trace-square-cube_proof_3}\Lambda_\omega(\Gamma) = \frac{1}{2(n-3)}\,\bigg(\Lambda_\omega^2(\Omega) - \frac{1}{3(n-2)}\,\Lambda_\omega^3(\Omega)\,\omega\bigg).\end{eqnarray}

Putting (\ref{eqn:Omega-omega_n-4_trace-trace-square-cube_proof_1}), (\ref{eqn:Omega-omega_n-4_trace-trace-square-cube_proof_2}) and (\ref{eqn:Omega-omega_n-4_trace-trace-square-cube_proof_3}) together, we get (\ref{eqn:Omega-omega_n-4_trace-trace-square-cube}), as claimed.  \hfill $\Box$

\vspace{2ex}

\noindent {\it End of proof of Proposition \ref{Prop:formula_rho_omega_P-omega}.} Taking $\Omega=i\partial\omega\wedge\bar\partial\omega$ in (\ref{eqn:Omega-omega_n-4_trace-trace-square-cube}) and using (\ref{eqn:formula_rho_omega_P-omega_proof_1}), we get (\ref{eqn:formula_rho_omega_P-omega}). \hfill $\Box$

\section{Appendix: Commutation relations}\label{section:Appendix} We briefly recall here some standard formulae that were used throughout the paper.

 Let $(X,\,\omega)$ be a compact complex Hermitian manifold. Recall the following standard Hermitian commutation relations ([Dem84], see also [Dem97, VII, $\S.1$]): \begin{eqnarray}\label{eqn:standard-comm-rel}\nonumber &  & (i)\,\,(\partial + \tau)^{\star} = i\,[\Lambda,\,\bar\partial];  \hspace{3ex} (ii)\,\,(\bar\partial + \bar\tau)^{\star} = - i\,[\Lambda,\,\partial]; \\
&  & (iii)\,\, \partial + \tau = -i\,[\bar\partial^{\star},\,L]; \hspace{3ex} (iv)\,\,
\bar\partial + \bar\tau = i\,[\partial^{\star},\,L],\end{eqnarray}

\noindent where the upper symbol $\star$ stands for the formal adjoint w.r.t. the $L^2$ inner product induced by $\omega$, $L=L_{\omega}:=\omega\wedge\cdot$ is the Lefschetz operator of multiplication by $\omega$, $\Lambda=\Lambda_{\omega}:=L^{\star}$ and $\tau_\omega =\tau:=[\Lambda,\,\partial\omega\wedge\cdot]$ is the torsion operator (of order zero and type $(1,\,0)$) associated with the metric $\omega$.

\vspace{3ex}

Other standard formulae (see e.g. [Voi02]) are the following: \begin{eqnarray}\label{eqn:standard_comm_Lambda-L-powers}\nonumber &  & (i)\,\, [\Lambda,\,L] = (n-k)\,\mbox{Id} \hspace{5ex} \mbox{on}\hspace{1ex} k\mbox{-forms, for every non-negative integer}\hspace{1ex} k; \\
\nonumber & & (ii)\,\, [L^r,\,\Lambda] = r(k-n+r-1)\,L^{r-1} \hspace{5ex} \mbox{on}\hspace{1ex} k\mbox{-forms, for all integers}\hspace{1ex} k\geq 0, r\geq 2; \\
 & & (iii)\,\,\star\,L = \Lambda\,\star \hspace{3ex}\mbox{and}\hspace{3ex} \star\,\Lambda = L\,\star.\end{eqnarray}

\vspace{3ex}

We also used the following result involving again the torsion operator $\tau_\omega$. 

\begin{Lem}\label{Lem:Michelson_torsion-form} Let $(X,\,\omega)$ be a compact complex Hermitian manifold with $\mbox{dim}_\C X=n$. The following identities hold: \begin{equation}\label{eqn:Michelson_torsion-form}-\frac{1}{2}\,\bar\tau_\omega^\star\omega \stackrel{(i)}{=} \bar\partial^\star_\omega\omega \stackrel{(ii)}{=} i\,\Lambda_\omega(\partial\omega).\end{equation} 

  In particular, $\omega$ is {\bf balanced} if and only if $\bar\tau_\omega^\star\omega=0$.

\end{Lem}

\noindent {\it Proof.} $\bullet$ To prove identity (i) in (\ref{eqn:Michelson_torsion-form}), we will show that the multiplication operators by the $(1,\,0)$-forms $\bar\tau^\star\omega$ and $-2\bar\partial^\star\omega$ acting on functions, namely $$\bar\tau^\star\omega\wedge\cdot, \hspace{1ex} -2\bar\partial^\star\omega\wedge\cdot:C^\infty_{0,\,0}(X,\,\C)\longrightarrow C^\infty_{1,\,0}(X,\,\C),$$ coincide by showing that their adjoints $$(\bar\tau^\star\omega\wedge\cdot)^\star, \hspace{1ex} (-2\bar\partial^\star\omega\wedge\cdot)^\star:C^\infty_{1,\,0}(X,\,\C)\longrightarrow C^\infty_{0,\,0}(X,\,\C)$$ coincide.

Let $\alpha\in C^\infty_{1,\,0}(X,\,\C)$ and $g\in C^\infty_{0,\,0}(X,\,\C)$ be arbitrary. We have:

\begin{eqnarray}\label{eqn:Michelson_torsion-form_proof_1}\langle\langle(\bar\partial^\star\omega\wedge\cdot)^\star\alpha,\,g \rangle\rangle = \langle\langle\bar{g}\alpha,\,\bar\partial^\star\omega\rangle\rangle = \langle\langle\bar\partial(\bar{g}\alpha),\,\omega\rangle\rangle = \int\limits_X\bar\partial(\bar{g}\alpha)\wedge\star\omega = \int\limits_X\bar{g}\alpha\wedge\bar\partial\omega_{n-1},\end{eqnarray} where we put $\omega_{n-1}:=\omega^{n-1}/(n-1)!$ and we used the standard identity $\star\omega = \omega_{n-1}$.

Meanwhile, we have: \begin{eqnarray}\label{eqn:Michelson_torsion-form_proof_2}\nonumber\langle\langle(\bar\tau^\star\omega\wedge\cdot)^\star\alpha,\,g \rangle\rangle & = &  \langle\langle\bar{g}\alpha,\,\bar\tau^\star\omega\rangle\rangle = \langle\langle\bar{g}\,\bar\tau(\alpha),\,\omega\rangle\rangle = \langle\langle\bar{g}\,\Lambda(\bar\partial\omega\wedge\alpha),\,\omega\rangle\rangle \\
  \nonumber & = & \langle\langle\bar\partial\omega\wedge\alpha,\,g\,\omega^2\rangle\rangle = \int\limits_X\bar\partial\omega\wedge\alpha\wedge\star(\bar{g}\omega^2) = -2\,\int\limits_X\bar{g}\alpha\wedge \bar\partial\omega\wedge\omega_{n-2} \\
  & = & -2\,\int\limits_X\bar{g}\alpha\wedge\bar\partial\omega_{n-1},\end{eqnarray} where for the third identity on the first line we used the definition $\bar\tau = [\Lambda,\,\bar\partial\omega\wedge\cdot]$ of $\bar\tau$ and the fact that $\Lambda(\alpha)=0$ for bidegree reasons, while for the third identity on the second line we used the standard identity $\star\omega_2 = \omega_{n-2}$, where $\omega_2:=\omega^2/2!$.

Comparing (\ref{eqn:Michelson_torsion-form_proof_1}) and (\ref{eqn:Michelson_torsion-form_proof_2}), we get $\langle\langle(\bar\tau^\star\omega\wedge\cdot)^\star\alpha,\,g \rangle\rangle = -2\,\langle\langle(\bar\partial^\star\omega\wedge\cdot)^\star\alpha,\,g \rangle\rangle$ for all $\alpha$ and $g$. Hence $(\bar\tau^\star\omega\wedge\cdot)^\star = -2\,(\bar\partial^\star\omega\wedge\cdot)^\star$, which proves (i) of (\ref{eqn:Michelson_torsion-form}).

\vspace{1ex}

$\bullet$ To prove identity (ii) in (\ref{eqn:Michelson_torsion-form}), we start from the Hermitian commutation relation (ii) in (\ref{eqn:standard-comm-rel}): $$[\Lambda,\,\partial] = i\,(\bar\partial^\star+\bar\tau^\star)$$ that we apply to $\omega$. We get the equivalent identities: \begin{eqnarray*}[\Lambda,\,\partial]\,\omega = i\bar\partial^\star\omega + i\bar\tau^\star\omega \iff \Lambda(\partial\omega) - \partial(\Lambda\omega) = -i\bar\partial^\star\omega  \iff \Lambda(\partial\omega) = -i\bar\partial^\star\omega,\end{eqnarray*} the last of which is (ii) of (\ref{eqn:Michelson_torsion-form}), where for the first equivalence we used the identity $\bar\tau^\star\omega = -2\,\bar\partial^\star\omega$ proved above as (i) in (\ref{eqn:Michelson_torsion-form}), while for the second equivalence we used the fact that $\Lambda\omega = n$, hence $\partial(\Lambda\omega)=0$. \hfill $\Box$

\vspace{3ex}

\noindent {\bf References.} \\

\vspace{1ex}

\noindent [Ang11]\, D. Angella --- {\it The Cohomologies of the Iwasawa Manifold and of Its Small Deformations} --- J. Geom. Anal. (2013), no. 3, 1355–1378.

\vspace{1ex}

\noindent [CE53]\, E. Calabi, B. Eckmann --- {\it A Class of Compact, Complex Manifolds Which Are Not Algebraic} --- Ann. of Math. {\bf 58} (1953) 494-500.

\vspace{1ex}

\noindent [Dem 84]\, J.-P. Demailly --- {\it Sur l'identit\'e de Bochner-Kodaira-Nakano en g\'eom\'etrie hermitienne} --- S\'eminaire d'analyse P. Lelong, P. Dolbeault, H. Skoda (editors) 1983/1984, Lecture Notes in Math., no. {\bf 1198}, Springer Verlag (1986), 88-97.

\vspace{1ex}

\noindent [Dem 97]\, J.-P. Demailly --- {\it Complex Analytic and Algebraic Geometry}---http://www-fourier.ujf-grenoble.fr/~demailly/books.html

\vspace{1ex}

\noindent [DP21]\, S.Dinew, D. Popovici --- {\it A Generalised Volume Invariant for Aeppli Cohomology Classes of Hermitian-Symplectic Metrics} --- Adv. Math. (2021), https://doi.org/10.1016/j.aim.2021.108056

\vspace{1ex}

\noindent [DP22]\, S. Dinew, D. Popovici --- {\it A Variational Approach to SKT and Balanced Metrics} --- arXiv e-print DG 2209.12813v1

\vspace{1ex}

\noindent [FU13]\, A. Fino, L. Ugarte --- {\it On Generalized Gauduchon Metrics} --- Proc. Edinburgh Math. Soc. {\bf 56} (2013), 733-753.

\vspace{1ex}

\noindent [FV15]\, A. Fino, L. Vezzoni --- {\it Special Hermitian metrics on compact solvmanifolds} --- J. Geom. Phys. {\bf91} (2015), 40-53.

\vspace{1ex}

\noindent [FWW13]\, J. Fu, Z. Wang, D. Wu --- {\it Semilinear Equations, the $\gamma_k$ Function and Generalized Gauduchon Metrics} --- J. Eur. Math. Soc. {\bf 15} (2013), 659-680.

\vspace{1ex}

\noindent [Gau77a]\, P. Gauduchon --- {\it Le th\'eor\`eme de l'excentricit\'e nulle} --- C. R. Acad. Sci. Paris, S\'er. A, {\bf 285} (1977), 387-390.

\vspace{1ex}

\noindent [Gau77b]\, P. Gauduchon --- {\it Fibr\'es hermitiens \`a endomorphisme de Ricci non n\'egatif} --- Bull. Soc. Math. France {\bf 105} (1977) 113-140.

\vspace{1ex}

\noindent [HP96]\, P. S. Howe, G. Papadopoulos --- {\it Twistor Spaces for HKT Manifolds} --- Phys. Lett. B {\bf 379} (1996) 80.

\vspace{1ex}

\noindent [IP13]\, S. Ivanov, G. Papadopoulos --- {\it Vanishing Theorems on $(l/k)$-strong K\"ahler Manifolds with Torsion} ---  Adv.Math. {\bf 237} (2013) 147-164.

\vspace{1ex}

\noindent [JY93]\, J. Jost, S.-T. Yau --- {\it A Nonlinear Elliptic System for Maps from Hermitian to Riemannian Manifolds and Rigidity Theorems in Hermitian Geometry} --- Acta Math. {\bf 170} (2) (1993) 221–254. Correction in: “A Nonlinear Elliptic System for Maps from Hermitian to Riemannian Manifolds and Rigidity Theorems in Hermitian Geometry”, Acta Math. {\bf 173} (2) (1994) 307.

\vspace{1ex}

\noindent [MT01]\, K. Matsuo, T. Takahashi --- {\it On Compact Astheno-K\"ahler Manifolds} --- Colloq. Math. {\bf 89} (2001), 213–221.

\vspace{1ex}

\noindent [Nak75]\, I. Nakamura --- {\it Complex parallelisable manifolds and their small deformations} --- J. Diff. Geom. {\bf 10} (1975), 85-112.

\vspace{1ex}

\noindent [Pop15]\, D. Popovici --- {\it Aeppli Cohomology Classes Associated with Gauduchon Metrics on Compact Complex Manifolds} --- Bull. Soc. Math. Fr. {\bf 143} (3) (2015) 1-37.

\vspace{1ex}

\noindent [Siu80]\, Y.-T. Siu --- {\it The Complex-Analyticity of Harmonic Maps and the Strong Rigidity of Compact K\"ahler Manifolds} --- Ann. Math. {\bf 112} (1), (1980), 73-111. 

\vspace{1ex}

\noindent [TT17]\, N. Tardini, A. Tomassini --- {\it On geometric Bott-Chern formality and deformations} --- Annali di Matematica  {\bf 196} (1), (2017), doi: 10.1007/s10231-016-0575-6

\vspace{1ex}

\noindent [Voi02]\, C. Voisin --- {\it Hodge Theory and Complex Algebraic Geometry. I.} --- Cambridge Studies in Advanced Mathematics, 76, Cambridge University Press, Cambridge, 2002.

\vspace{6ex}

\noindent Universit\'e Paul Sabatier, Institut de Math\'ematiques de Toulouse

\noindent 118, route de Narbonne, 31062, Toulouse Cedex 9, France

\noindent Email: popovici@math.univ-toulouse.fr

\end{document}